\documentclass{article} 
\usepackage{amsmath,amssymb,amsthm} 
\usepackage[OT2,T1]{fontenc}
\usepackage[utf8]{inputenc}
\usepackage{enumitem}
\usepackage{graphicx}
\usepackage{hyperref}
\usepackage{tikz}
\usetikzlibrary{cd}
\usepackage{tikz-cd}
\usepackage{xcolor}
\usepackage{ stmaryrd }
\usepackage{mathrsfs}
\usepackage{titling}
\usepackage{fullpage}
\usepackage{mathtools}
\usepackage[english]{babel}
\usepackage{setspace}
\usepackage[normalem]{ulem}

\usepackage[alphabetic,lite]{amsrefs}

\usepackage{arydshln}

\usetikzlibrary{arrows,positioning,decorations.pathmorphing,
	decorations.markings
}
\tikzset{inner sep=0pt,
	root/.style={circle,draw,minimum size=7pt,thick},
	fatroot/.style={circle,draw,minimum size=10pt,thick},
	short root/.style={circle,fill,minimum size=7pt},
	doublearrow/.style={postaction={decorate},
		decoration={markings,mark=at position .7
			with {\arrow{angle 60}}},double distance=3pt,thick}
}

\DeclareUnicodeCharacter{00A0}{ }

\newtheorem{proposition}{Proposition}[section]
\newtheorem{definition}[proposition]{Definition}
\newtheorem{remark}[proposition]{Remark}

\newtheorem{theorem}[proposition]{Theorem}
\newtheorem{lemma}[proposition]{Lemma}
\newtheorem{corollary}[proposition]{Corollary}
\newtheorem{example}[proposition]{Example}

\newtheorem{conjecture}[proposition]{Conjecture}
\numberwithin{equation}{subsection}
\newtheorem{convention}[proposition]{Convention}

\newcommand{\G}{\mathbb{G}}
\DeclareMathOperator{\GL}{GL}
\DeclareMathOperator{\PGL}{PGL}
\DeclareMathOperator{\SL}{SL}

\DeclareMathOperator{\PSp}{PSp}

\DeclareMathOperator{\SO}{SO}
\DeclareMathOperator{\PSO}{PSO}
\DeclareMathOperator{\Spin}{Spin}

\DeclareMathOperator{\liesl}{\mathfrak{sl}}

\DeclareMathOperator{\so}{\mathfrak{so}}

\DeclareMathOperator{\lieg}{\mathfrak{g}}
\DeclareMathOperator{\lieh}{\mathfrak{h}}
\DeclareMathOperator{\liet}{\mathfrak{t}}
\DeclareMathOperator{\Ad}{Ad}
\DeclareMathOperator{\Lie}{Lie}

\DeclareMathOperator{\Hom}{Hom}
\DeclareMathOperator{\Aut}{Aut}

\DeclareMathOperator{\Gal}{Gal}

\DeclareMathOperator{\End}{End}
\DeclareMathOperator{\Sym}{Sym}
\DeclareMathOperator{\vol}{vol}

\DeclareMathOperator{\rank}{rank}
\DeclareMathOperator{\image}{image}

\DeclareMathOperator{\Isom}{Isom}

\DeclareMathOperator{\Spec}{Spec}

\DeclareMathOperator{\Frac}{Frac}

\DeclareMathOperator{\Pic}{Pic}

\newcommand{\sh}[1]{\mathscr{#1}}

\newcommand{\A}{\mathbb{A}}
\renewcommand{\P}{\mathbb{P}}

\renewcommand{\O}{\mathcal{O}}
\newcommand{\GIT}{\mathbin{/\mkern-6mu/}}
\newcommand{\HH}{\mathrm{H}}

\newcommand{\Real}{\mathbb{R}}
\newcommand{\Q}{\mathbb{Q}}
\newcommand{\Z}{\mathbb{Z}}
\newcommand{\F}{\mathbb{F}}

\DeclareMathOperator{\Sel}{Sel}
\DeclareMathOperator{\rk}{rk}

\DeclareSymbolFont{cyrletters}{OT2}{wncyr}{m}{n}
\DeclareMathSymbol{\Sha}{\mathalpha}{cyrletters}{"58}

\newcommand{\extp}{\@ifnextchar^\@extp{\@extp^{\,}}}
\def\@extp^#1{\mathop{\bigwedge\nolimits^{\!#1}}}

\DeclareUnicodeCharacter{00A0}{ }
\newcommand{\define}[1]{{\fontfamily{cmss}\selectfont{#1}}}

\newcommand{\intbigG}{\underline{G}}
\newcommand{\intbigH}{\underline{H}}
\newcommand{\intbigV}{\underline{V}}
\newcommand{\intbigB}{\underline{B}}
\newcommand{\intbigT}{\underline{T}}
\newcommand{\intbigt}{\underline{\mathfrak{t}}}

\newcommand{\bigH}{H}
\newcommand{\bigh}{\lieh}
\newcommand{\bigG}{G}
\newcommand{\bigB}{B}
\newcommand{\bigP}{P}
\newcommand{\bigPopp}{\bar{P}}
\newcommand{\bigN}{N}
\newcommand{\bigNopp}{\bar{N}}
\newcommand{\bigT}{T}

\renewcommand{\bigg}{\mathfrak{g}}
\newcommand{\bigV}{V}
\newcommand{\bigtheta}{\theta}

\newcommand{\bigkappa}{\kappa}

\newcommand{\bigpi}{\pi}

\newcommand{\Jac}{J}
\newcommand{\CJac}{\bar{J}}

\newcommand{\intbigcurve}{\underline{C}}
\newcommand{\Bun}[1]{\mathsf{B}#1}

\newcommand{\height}{\mathrm{ht}}

\DeclareMathOperator{\disc}{disc}
\DeclareMathOperator{\ord}{ord}
\DeclareMathOperator{\rs}{rs}
\DeclareMathOperator{\reg}{reg}
\newcommand{\Siegel}{\mathfrak{S}}
\DeclareMathOperator{\Pff}{Pff}

\tikzset{
  symbol/.style={
    draw=none,
    every to/.append style={
      edge node={node [sloped, allow upside down, auto=false]{$#1$}}}
  }
}

\usepackage{microtype}

\title{Graded Lie Algebras, Compactified Jacobians and Arithmetic Statistics}
\author{Jef Laga}

\begin{document}

\maketitle

\begin{abstract}
A simply laced Dynkin diagram gives rise to a family of curves over $\Q$ and a coregular representation, using deformations of simple singularities and Vinberg theory respectively.
Thorne has conjectured and partially proven a strong link between the arithmetic of these curves and the rational orbits of these representations.

In this paper, we complete Thorne's picture and show that $2$-Selmer elements of the Jacobians of the smooth curves in each family can be parametrised by integral orbits of the corresponding representation. 
Using geometry-of-numbers techniques, we deduce statistical results on the arithmetic of these curves. 
We prove these results in a uniform manner.
This recovers and generalises results of Bhargava, Gross, Ho, Shankar, Shankar and Wang. 

The main innovations are: an analysis of torsors on affine spaces using results of Colliot-Th\'el\`ene and the Grothendieck--Serre conjecture, a study of geometric properties of compactified Jacobians using the Białynicki-Birula decomposition, and a general construction of integral orbit representatives.
\end{abstract}

\tableofcontents

\section{Introduction}

\subsection{Context}


This paper is a contribution to arithmetic statistics of algebraic curves: given a family of curves $\mathcal{F}$ over $\Q$, what can be said about the rational points of $C$ (or related objects) as $C$ varies in $\mathcal{F}$?
Over the last twenty years, Bhargava and his collaborators have made spectacular progress in this direction.
One of their key ideas is that many arithmetic objects can be parametrised by rational or integral orbits of a representation $(G,V)$. 
When the representation is coregular, meaning that the ring of invariants $\Q[V]^G$ is a polynomial ring, they have developed powerful geometry-of-numbers techniques to count integral orbits of $V$. 
Combining orbit parametrisations with these counting techniques has led to many striking results; see \cite{BS-2selmerellcurves,BS-3Selmer, BS-4Selmer, BS-5Selmer, Bhargava-Gross-hyperellcurves, BhargavaGrossWang-positiveproportionnopoints} for some highlights and \cite{Ho-howmanyrationalpoints, Bhargava-ICMreport} for surveys of these results.

This raises the question: how does one find such orbit parametrisations? 
Typically they arise from classical algebro-geometric constructions. 
For example, elements of $\Sel_2E$ for an elliptic curve $E/\Q$ correspond to locally soluble genus-$1$ curves $C$ that are double covers of $\P^1_{\Q}$ \cite[\S1.3]{CremonaFisherStoll}, so give rise to $\PGL_2(\Q)$-orbits of binary quartic forms \cite[Theorem 3.5]{BS-2selmerellcurves}.
This example goes back to Birch and Swinnerton-Dyer \cite{BirchSwinnertonDyer-notesonellcurves1} (building on ideas of Mordell) and has been used by Bhargava and Shankar to compute the average size of the $2$-Selmer group of elliptic curves \cite{BS-2selmerellcurves}.
See \cite{BhargavaHo-coregularspacesgenusone} for an exhaustive list of orbit parametrisations of genus-$1$ curves which are obtained using similar (but more difficult) algebro-geometric constructions.
See also \cite{Bhargava-Gross-hyperellcurves} for an orbit parametrisation of $2$-Selmer groups of odd hyperelliptic curves using the geometry of pencils of quadrics \cite{Wang-phdthesis}.
Even though these considerations have been hugely successful, Wei Ho writes that `Finding appropriate groups $G$ and vector spaces $V$ related to the Selmer elements is still a relatively ad hoc process' \cite[p 45]{Ho-howmanyrationalpoints}.

Gross \cite{Gross-BhargavasrepresentationsandVinberg} observed that most coregular representations employed in arithmetic statistics arise from \emph{Vinberg theory}, that is the theory of graded Lie algebras. 
This suggests the possibility to take Vinberg theory as a starting point, and to attempt to naturally construct families of curves in this setting.
This is exactly the perspective taken in Thorne's PhD thesis \cite{Thorne-thesis} in the case of $2$-Selmer groups.
Given a simply laced Dynkin diagram of type $A, D,E$, he \emph{canonically} constructs a family of curves and a coregular representation whose rational orbits should be related to the arithmetic of the curves in the family. 
This canonical construction unifies many orbit parametrisations in the literature and has already produced new results in arithmetic statistics; see \cite{Thorne-E6paper, Romano-Thorne-ArithmeticofsingularitiestypeE, Thorne-averagesizeelliptictwomarkedfunctionfields}. 
However, to obtain all the expected consequences it remained to be shown that all elements of the $2$-Selmer group give rise to rational orbits \cite[Conjecture 4.16]{Thorne-thesis} and that such rational orbits admit integral representatives.

The main goal of this paper is to resolve both these questions, and to do so in a uniform manner for all the ADE-families considered. 
By using geometry-of-numbers techniques developed by Bhargava and his collaborators, we obtain an upper bound on the average size of the $2$-Selmer group of the Jacobians of the smooth curves in each family. 
This has consequences for the ranks of the Jacobians and the rational points of the curves in these families.

\subsection{Statement of results}

Let $\mathsf{D}$ be a Dynkin diagram of type $A_n$, $D_n$ or $E_n$ and let $C\rightarrow B$ be the family of projective curves over $\Q$ with affine equation given by Table \ref{table: introduction different cases}.
For example, if $\mathsf{D}= A_{2g}$, then $B = \Spec \Q[p_2,\dots,p_{2g+1}]$ and $C\rightarrow B$ is the family of all monic odd hyperelliptic curves of genus $g$.
If $\mathsf{D} = E_7$, then $C\rightarrow B$ is the family of all plane quartic curves with a marked rational flex point.
The family $C\rightarrow B$ is a semi-universal deformation of its central fibre (by setting all coefficients $p_i$ equal to zero), which is a simple singularity of type $\mathsf{D}$. (See Proposition \ref{proposition: first properties of the family of curves}.)
We exclude the case $\mathsf{D} = A_1$.

\begin{table}
\centering
\begin{tabular}{|l | l | c | c | c| }
	\hline
	   Type & Equation & $m$  \\
	\hline       
	$A_{2g}$ & $y^2 = x^{2g+1} + p_2x^{2g-1} + \dots +p_{2g+1}$& $1$  \\
	$A_{2g+1}$ & $y^2 = x^{2g+2} + p_2x^{2g} + \dots +p_{2g+2}$& $2$  \\
	$D_{2g+2}\, (g\geq 1)$ & $y(xy+p_{2g+2}) = x^{2g+1}+p_2x^{2g}+p_4x^{2g-1} + \dots + p_{4g+2} $ & $3$  \\
	$D_{2g+1}\, (g\geq 2)$ & $y(xy+p_{2g+1}) = x^{2g}+p_2x^{2g-1}+p_4x^{2g-2} + \dots + p_{4g} $ & $2$ \\
	$E_6$ & $y^3 = x^4+(p_2x^2+p_5x+p_8)y + (p_6x^2+p_9x+p_{12})$ & $1$ \\
	$E_7$ & $y^3 = x^3y +p_{10}x^2 +x(p_2y^2 + p_8y + p_{14}) + p_6y^2 + p_{12}y + p_{18}$ & $2$ \\
	$E_8$ & $y^3 = x^5+ (p_2x^3+p_8x^2+p_{14}x+p_{20})y+(p_{12}x^3+p_{18}x^2+p_{24}x+p_{30}) $ & $1$ \\
	\hline

\end{tabular}
\caption{Families of curves}
\label{table: introduction different cases}
\end{table}

Write $B^{\rs} \subset B$ for the locus above which $C\rightarrow B$ is smooth, the complement of a discriminant hypersurface.
For every field $k/\Q$ and $b\in B^{\rs}(k)$, write $J_b$ for the Jacobian of the smooth projective curve $C_b$, an abelian variety over $k$ of dimension equal to the genus of $C_b$. 
Our first main theorem is an orbit parametrisation for elements of $J_b(k)/2J_b(k)$. 

To each diagram $\mathsf{D}$ one may canonically associate a representation $V$ of a reductive group $G/\Q$. This construction, due to Thorne \cite{Thorne-thesis}, is recalled in  \S\ref{subsection: a stable Z/2Z-grading} and is based on Vinberg's theory of graded Lie algebras.
See \S\ref{subsection: explicit determination G and V} for an explicit description of $G$ and $V$, although we will almost never use this description.
The geometric quotient $V \GIT G = \Spec \Q[V]^G$ (parametrising $G$-invariant polynomials of $V$) turns out to be isomorphic to $B$. 
For every field $k/\Q$ and $b\in B(k)$, write $V_b$ for the subset of elements of $V$ which map to $b$ under the map $V\rightarrow V \GIT G \simeq B$.

\begin{theorem}[Theorem \ref{theorem: inject 2-descent orbits}]  \label{theorem: intro inject 2-descent orbits}
For every field $k/\Q$ and element $b\in B^{\rs}(k)$, there exists an injection $\eta_b\colon J_b(k)/2J_b(k) \hookrightarrow G(k)\backslash V_b(k)$ compatible with base change. 
\end{theorem}

See Theorem \ref{theorem: inject 2-descent orbits} for a more precise formulation and an explicit construction of this injection.
Using a local-global principle for $G$, one can also embed the $2$-Selmer group of $J_b$ inside the $G(\Q)$-orbits of $V(\Q)$. 
Recall that the $2$-Selmer group of an abelian variety $A/\Q$ is a finite dimensional $\F_2$-vector space $\Sel_2A$ defined by local conditions and fitting inside an exact sequence 
$$
0 \rightarrow A(\Q)/2A(\Q) \rightarrow \Sel_{2}A \rightarrow \Sha(A/\Q)[2]\rightarrow 0.
$$ 

\begin{theorem}[Corollary \ref{corollary: inject 2-Selmer orbits}]  \label{theorem: intro inject 2-selmer orbits}
For every $b\in B^{\rs}(\Q)$, the injection $\eta_b$ extends to an injection $\Sel_2 J_b \hookrightarrow G(\Q) \backslash V_b(\Q)$. 
\end{theorem}

If $\mathsf{D}$ is of type $A_2$, $C\rightarrow B$ is the family of elliptic curves in short Weierstrass form and we essentially recover the orbit-parametrisation of Birch--Swinnerton-Dyer \cite{BirchSwinnertonDyer-notesonellcurves1} used by Bhargava--Shankar \cite{BS-2selmerellcurves}.
If $\mathsf{D}$ is of type $A_{2g}$, we recover the orbit parametrisation of Bhargava and Gross \cite{Bhargava-Gross-hyperellcurves}.

Crucially, we additionally show that the $G(\Q)$-orbits corresponding to $\Sel_2 J_b$ using Theorem \ref{theorem: intro inject 2-selmer orbits} have integral representatives away from small primes, see Corollary \ref{corollary: weak global integral representatives}.
Using geometry-of-numbers techniques to count integral orbits of $V$, Theorem \ref{theorem: intro inject 2-selmer orbits} may thus be used to give an upper bound on the average size of the $2$-Selmer group of $\Sel_2J_b$. 
Using the identification $B = \Spec \Q[p_{d_1}, \dots,p_{d_r}]$ from Table \ref{table: introduction different cases}, let $\mathcal{F}$ be the subset of elements $b = (p_{d_1}(b),\dots,p_{d_r}(b)) \in \Z^r$ with $b\in B^{\rs}(\Q)$. 
We define the height of $b\in \mathcal{F}$ by the formula
\begin{align*}
\height(b) \coloneqq \max\left(|p_{d_1}(b)|^{1/d_1}, \dots,|p_{d_r}(b)|^{1/d_r}\right).
\end{align*}
Note that for every $X\in \Real_{>0}$, the set $\{b\in \mathcal{F} \mid \height(b) <X \}$ is finite. 
To state the next theorem, note that each curve $C_b$ has points at infinity not lying in the affine patch of Table \ref{table: introduction different cases}, and we call those points the \emph{marked points}. 
Their cardinality is displayed in Table \ref{table: introduction different cases}.

\begin{theorem}[Theorem \ref{theorem: main theorem}]
\label{theorem: intro average size 2-Selmer}
Let $m$ be the number of marked points.
Then when ordered by height, the average size of the $2$-Selmer group of $J_b$ for $b\in \mathcal{F}$ is bounded above by $3\cdot 2^{m-1}$.
More precisely, we have 
	\begin{equation*}
	\limsup_{X\rightarrow +\infty} \frac{ \sum_{b\in \mathcal{F},\; \height(b)<X }\# \Sel_2 J_b    }{\#  \{b \in \mathcal{F}\mid \height(b) < X\}}	\leq 3 \cdot 2^{m-1}.
	\end{equation*}
\end{theorem}

The same result holds true even if we impose finitely many congruence conditions on $\mathcal{F}$.
Assuming a certain plausible uniformity estimate (Conjecture \ref{conjecture: uniformity estimate}), we show that the limit exists and the bound $3\cdot 2^{m-1}$ is sharp, see \S\ref{subsection: a conditional lower bound}.
The constant $3\cdot 2^{m-1}$ is consistent with the heuristics of Poonen and Rains \cite{PoonenRains-maximalisotropic}, once we incorporate the fact that the $m$ marked points give rise to a `trivial' subgroup of $\Sel_2J_b$ which has size $2^{m-1}$ most of the time (Proposition \ref{proposition: average size trivial selmer group}).
See \S\ref{subsection: relation to other works} for a comparison of this theorem with previously obtained results.

Just like in previous cases, Theorem \ref{theorem: intro average size 2-Selmer} implies a bound on the average of the Mordell--Weil rank $\rk(\Jac_b)$ of $\Jac_b$, the rank of the finitely generated abelian group $\Jac_b(\Q)$.
Let $\Sel_2^{triv}J_b\subset \Sel_2 J_b$ be the subgroup generated by the differences of the marked points. 
The chain of inequalities $$\rk J_b \leq \dim_{\F_2}(\Sel_2 J_b/\Sel^{triv}_2J_b)+(m-1)\leq \frac{1}{2}(\#(\Sel_2J_b/\Sel_2^{triv}J_b))+(m-1)$$
combined with Theorem \ref{theorem: intro average size 2-Selmer} and the fact that $\Sel_2^{triv} J_b$ is of size $2^{m-1}$ for $100\%$ of curves, we obtain:
\begin{corollary}
Let $m$ be the number of marked points of the family $C\rightarrow B$.
Then when ordered by height, the average rank $\rk(\Jac_b)$ where $b\in \mathcal{F}$ is bounded above by $m+1/2$.
\end{corollary}

\subsection{Relation to other works}\label{subsection: relation to other works}

Theorem \ref{theorem: intro average size 2-Selmer} has been previously obtained for many $\mathsf{D}$:
\begin{itemize}
\item $A_2$: Bhargava--Shankar \cite{BS-2selmerellcurves}, who prove in this case that the average is exactly $3$.
\item $A_{2g}$: Bhargava--Gross \cite{Bhargava-Gross-hyperellcurves}.
\item $A_{2g+1}, \, g\geq 2$: Shankar--Wang \cite{ShankarWang-hypermarkednonweierstrass}.
\item $D_{2g+1}, \, g\geq 2$: Shankar \cite{Shankar-2selmerhypermarkedpoints}.
\item $A_3, D_4$: Bhargava--Ho \cite[Theorem 1.1(c),(g)]{BhargavaHo-2Selmergroupsofsomefamilies}. 
\item $E_6$: Laga \cite{Laga-E6paper}.
\end{itemize}
All these works combine geometry-of-numbers techniques with the orbit parametrisation of Theorem \ref{theorem: intro inject 2-selmer orbits} to obtain Theorem \ref{theorem: intro average size 2-Selmer} in their specific case, just as we do here.
However, their construction of orbits (in other words, the proof of Theorem \ref{theorem: intro inject 2-selmer orbits}) and analysis of the representation $(G,V)$ requires specific arguments in each case. 
One of the main points of this paper is that we are able to prove Theorem \ref{theorem: intro average size 2-Selmer} in a uniform way.
Inspecting the above list, we see that the only cases not previously considered in the literature are $\mathsf{D} = D_{2g+2}$ with $g\geq 2$, $E_7$ and $E_8$.
Concretely, this concerns the universal family of hyperelliptic genus $g$ curves with two non-conjugate non-Weierstrass points ($D_{2g+2}$), the universal family of plane quartic curves with a marked flex point ($E_7$) and the universal family of trigonal genus $4$ curves with a marked triple ramification point ($E_8$). 

Theorem \ref{theorem: intro average size 2-Selmer} has a number of interesting consequences for the rational points of $C_b$, typically using various forms of the Chabauty--Coleman method. See for example \cite[Corollary 1.4]{Bhargava-Gross-hyperellcurves} and \cite{PoonenStoll-Mosthyperellipticnorational} for such results in the case $\mathsf{D} = A_{2g}$, and \cite[Corollary 1.3 \& Theorem 1.4]{Laga-E6paper} in the case $\mathsf{D} = E_6$.
Theorem \ref{theorem: intro average size 2-Selmer} should give similar such consequences for $\mathsf{D} = D_{2g+2}$, $E_7$ and $E_8$, but we have not pursued this in this paper.

We describe what is new in this paper compared with Thorne's work. 
He has shown the analogue of Theorem \ref{theorem: intro inject 2-descent orbits} for the subset of $J_b(k)/2J_b(k)$ lying in the image of the Abel--Jacobi map $C_b(k) \rightarrow J_b(k)/2J_b(k)$ with respect to a fixed marked point \cite[Theorem 4.15]{Thorne-thesis}.
This allowed him and Beth Romano to deduce arithmetic statistical results on the \emph{$2$-Selmer set} of the curve $C_b$ (a pointed subset of $\Sel_2 J_b$) and the integral points of the affine curve $C_b^{\circ}$, see \cite{Thorne-E6paper, Romano-Thorne-ArithmeticofsingularitiestypeE}.
The first main innovation of this work is the construction of orbits associated with \emph{all} elements of $J_b(k)/2J_b(k)$, as was conjectured in \cite[Conjecture 4.16]{Thorne-thesis}.
The second main innovation is an integral study of the representations $(G,V)$. 
In particular, we show that orbits arising from Theorem \ref{theorem: intro inject 2-selmer orbits} admit integral representatives away from small primes (Theorem \ref{theorem: integral representatives exist}). 
This technical result is essential for applying orbit-counting methods and allows us to obtain new results on the arithmetic of the curves $C_b$. 


\subsection{Method of proof}\label{subsection: method of proof}

We briefly describe the proof of Theorem \ref{theorem: intro inject 2-descent orbits}, which is the first main novelty of this paper.
Thorne has shown that the stabiliser $Z_G(v)$ of an arbitrary element $v\in V_b(k)$ is canonically isomorphic to $J_b[2]$, the $2$-torsion subgroup of the Jacobian of $C_b$.
In fact, there always exists a distinguished orbit $\kappa_b \in G(k)\backslash V_b(k)$ and a well-known lemma in arithmetic invariant theory (Lemma \ref{lemma: AIT full generality}) shows that by twisting $\kappa_b$ the set $G(k) \backslash V_b(k)$ can be identified with the pointed kernel of the map on Galois cohomology $\HH^1(k,Z_G(\kappa_b)) \rightarrow \HH^1(k,G)$. 

To prove Theorem \ref{theorem: intro inject 2-descent orbits}, it therefore suffices to prove that the composition
\begin{align}\label{equation: intro composite 2-descent with orbit classifying map}
J_b(k)/2J_b(k) \xrightarrow{\delta} \HH^1(k,J_b[2]) \simeq \HH^1(k,Z_G(\kappa_b)) \rightarrow \HH^1(k,G),
\end{align}
where $\delta$ is the $2$-descent map of $J_b$, is trivial.
We solve this problem by considering it universally. 
More precisely, a `categorified' version of \eqref{equation: intro composite 2-descent with orbit classifying map} associates to every element $P\in J_b(k)$ a $G$-torsor $T_P \rightarrow \Spec k$ such that its isomorphism class $[T_P]\in \HH^1(k,G)$ equals the image of $P$ under \eqref{equation: intro composite 2-descent with orbit classifying map}.
This process can be carried out in a relative setting: let $J^{\rs}$ be the relative Jacobian of the family of smooth curves $C|_{B^{\rs}} \rightarrow B^{\rs}$. 
Then we may construct a $G$-torsor $T\rightarrow J^{\rs}$ whose pullback along a point $P\colon \Spec k \rightarrow J^{\rs}$ is isomorphic to $T_P$. 
The crucial observation is that the geometry of the total space $J^{\rs}$ is very simple, despite the fibres of $J^{\rs} \rightarrow B^{\rs}$ being abelian varieties so arguably not so simple. 
For example, $J^{\rs}$ is a rational variety. This fact (or rather a similar, more precise statement), together with an analysis of $G$-torsors on affine spaces and progress on the Grothendieck--Serre conjecture, allows us to prove $J^{\rs}$ admits a Zariski open cover above which $T$ is trivial. This implies that each $T_P$ is trivial, proving the theorem.

To analyse the geometry of $J^{\rs}$, we introduce a compactification of $J^{\rs}$ over the whole of $B$: there exists a projective scheme $\bar{J} \rightarrow B$ restricting to $J^{\rs}$ over $B^{\rs}$ called the \emph{compactified Jacobian} of $C\rightarrow B$.
The scheme $\bar{J}$ parametrises rank-$1$ torsion-free sheaves following Altman--Kleiman \cite{AltmanKleiman-CompactifyingThePicardScheme},
with the caveat that in the reducible fibres of $C \rightarrow B$ we have to impose a stability condition in the sense of Esteves \cite{Esteves-compactifyingrelativejacobian} to obtain a well-behaved moduli problem.
The main selling point of this paper can be summarised as follows: geometric properties of $\bar{J}$ are very useful in the construction of orbits associated with elements of $J^{\rs}$.
For example, we show that even though the fibres of $\bar{J}\rightarrow B$ might be highly singular, $\bar{J}$ is a smooth and geometrically integral variety.
Moreover, the Białynicki-Birula decomposition from geometric representation theory shows that $\bar{J}$ has a decomposition into affine cells, so has a very transparent geometry.
The consequences for the geometry of $J^{\rs}$ are strong enough to carry out the strategy of the previous paragraph and consequently prove Theorem \ref{theorem: intro inject 2-descent orbits}.
The smoothness of $\bar{J}$ is essential and follows from the fact that $C\rightarrow B$ is a semi-universal deformation of its central fibre.

The second main innovation of this paper is a uniform construction of integral representatives, and we again exploit the geometry of the compactified Jacobian.
We achieve this by deforming to the case of square-free discriminant and using a general result on extending reductive group schemes over open dense subschemes of regular arithmetic surfaces (Lemma \ref{lemma: purity for the stack M}).
We are able to deform to this case using Bertini theorems over $\Q_p$ and $\F_p$ and (again!) the smoothness of $\bar{J}$. 
Many of the ideas were already present in our earlier work \cite{Laga-E6paper} treating the $E_6$ case, but we use stack-theoretic language to streamline it significantly.
Constructing integral orbits has often been a subtle point in the past, and we expect that our methods will have applications to settings different to the one considered here.

To deduce Theorem \ref{theorem: intro average size 2-Selmer} from Theorem \ref{theorem: intro inject 2-selmer orbits}, we use the robust geometry-of-numbers methods developed by Bhargava and his collaborators to count integral orbits coregular representations, about which we make two remarks.
Firstly, the only step in the counting argument that cannot be carried out uniformly is controlling orbits lying in the cuspidal region of the fundamental domain, the so-called `cutting of the cusp'. 
For every ADE diagram, this relies on combinatorial calculations in the associated root system. These calculations have appeared in the literature (on which we rely) except for the $D_{2g+2}$ case, which we handle explicitly in an appendix chapter.
This is the only part of the paper where we rely on the previous works listed at the beginning of \S\ref{subsection: relation to other works}.
It would be very interesting to find a less computational or even uniform proof for these calculations.
This remark extends to other representations employed in arithmetic statistics, for example the ones used in \cite{BS-4Selmer, BS-5Selmer}. 
See \cite[Table 2]{BS-5Selmer} for an example of the intricacies involved.
Secondly, the reason that we only obtain an upper bound in Theorem \ref{theorem: intro average size 2-Selmer} is our failure to prove the uniformity estimate of Conjecture \ref{conjecture: uniformity estimate}, which is the crucial ingredient for the square-free sieve needed to obtain the lower bound. 
For the representations considered in this paper this conjecture has only been solved in the $A_2$ case \cite[Theorem 2.13]{BS-2selmerellcurves}. 
Establishing more cases of this conjecture would be very interesting but seems difficult at present.


\subsection{Other coregular representations}



There are at least two ways in which coregular representations can arise from Vinberg theory that are not treated in this paper.
In both cases, we expect that our methods go a long way towards proving analogous results to Theorem \ref{theorem: intro average size 2-Selmer}.

Firstly, one may try to incorporate gradings on nonsimply laced Lie algebras (so of type $B,C, F,G$) into the picture. 
Again there will be families of curves, but the relevant Selmer groups may arise from a general isogeny, not just multiplication by an integer.
Such gradings have already appeared in the literature, implicitly and explicitly: a $\Z/2\Z$-grading on $G_2$ has been used to study $2$-Selmer groups of elliptic curves with a marked $3$-torsion point \cite[Theorem 1.1(f)]{BhargavaHo-2Selmergroupsofsomefamilies}; a $\Z/3\Z$-grading on $G_2$ has been used to study $3$-isogeny Selmer groups of the curves $y^2 = x^3+k$ \cite{BhargavaElkiesShnidman}; a $\Z/2\Z$-grading on $F_4$ has been used to study $2$-Selmer groups of a family of Prym surfaces \cite{Laga-F4paper}. 

Secondly, although their occurrence is more sporadic, there are also interesting $\Z/m\Z$-gradings on simple Lie algebras for $m\geq 3$; see for example \cite{Thorne-Romano-E8}, where the authors calculate the average size of the $3$-Selmer group of the family of odd genus-$2$ curves using a $\Z/3\Z$-grading on $E_8$. The representations used by Bhargava--Shankar for $3$-, $4$- and $5$-Selmer groups of elliptic curves \cite{BS-3Selmer, BS-4Selmer, BS-5Selmer} can also be interpreted this way.

To the best of our knowledge, every coregular representation appearing in the literature on arithmetic statistics of algebraic curves arises from Vinberg theory (that is, from a $\Z/m\Z$-grading on a semisimple Lie algebra), \emph{except} for one family of notable examples: the representation of $\SL_n$ acting on pairs of symmetric matrices $\Sym^2(n) \oplus \Sym^2(n)$.
This representation is used in \cite{BhargavaGrossWang-positiveproportionnopoints} to show that a positive proportion of locally soluble hyperelliptic curves over $\Q$ of fixed genus have no points over any odd degree extension. 
One feature that distinguishes their setting from ours is that their representation lacks a `Kostant section', which is related to the fact that the curves they study do not come with specified marked points.
We wonder if one can still interpret this representation in terms of Lie theory and study its arithmetic from this perspective.


\subsection{Organisation}

We now summarise the chapters of this paper.
In \S\ref{section: background} we recall some background results in Vinberg theory and arithmetic invariant theory.
In \S\ref{section: recollections of Thorne's thesis} we recall the constructions and main results of Thorne's thesis and introduce the Vinberg representation $(G,V)$ and family of curves $C\rightarrow B$.
In \S\ref{section: the mildly singular locus}, we extend the results of Thorne's thesis from the smooth fibres of $C\rightarrow B$ to those fibres admitting at most one singular nodal point.
In \S\ref{section: the compactified jacobian}, we introduce and study compactified Jacobians of the family $C\rightarrow B$.
In \S\ref{section: constructing orbits} we analyse torsors on affine spaces and use this and our results from the previous chapters to prove Theorems \ref{theorem: intro inject 2-descent orbits} and \ref{theorem: intro inject 2-selmer orbits} on the construction of orbits.
In \S\ref{integral representatives} we prove that such orbits admit integral representatives away from small primes.
In \S\ref{section: geometry-of-numbers}, we employ Bhargava's orbit-counting techniques and count integral orbits of the representation $(G,V)$.
We combine all the results from the previous chapters in \S\ref{section: the average size of the 2-Selmer group} to obtain Theorem \ref{theorem: intro average size 2-Selmer}.
In an appendix chapter, we perform some combinatorial calculations in the root system of type $D_{2n}$ to complete the proof of Proposition \ref{proposition: cutting off cusp} (cutting off the cusp) in this case.
We note that the main novel contributions of this paper lie in Chapters \ref{section: the mildly singular locus}, \ref{section: the compactified jacobian}, \ref{section: constructing orbits} and \ref{integral representatives}. 


\subsection{Acknowledgements}

This paper is a revised version of the author's PhD thesis, written under the supervision of Jack Thorne. I would like to thank him for many useful suggestions, stimulating conversations and for carefully reading an early draft of this manuscript.
I would also like to thank the anonymous referee for their careful reading of the paper. 
Finally, it is a pleasure to thank Alex Bartel, Jean-Louis Colliot-Th\'el\`ene, Tom Fisher, Jesse Leo Kass and Beth Romano for interesting discussions related to the contents of this paper.

\subsection{Notation}\label{subsection: notation}

\subsubsection{General}

For a field $k$ we write $k^s$ for a fixed separable closure and $\Gamma_k = \Gal(k^s/k)$ for its absolute Galois group. 

If $X$ is a scheme over $S$ and $T\rightarrow S$ a morphism we write $X_T$ for the base change of $X$ to $T$. If $T = \Spec A$ is an affine scheme we also write $X_A$ for $X_T$. 

If $G$ is a smooth group scheme over $S$ then we write $\HH^1(S,G)$ for the set of isomorphism classes of \'etale sheaf torsors under $G$ over $S$, which is a pointed set coming from nonabelian \v{C}ech cohomology. If $S = \Spec R$ we write $\HH^1(R,G)$ for the same object.
If $k$ is a field then $\HH^1(k,G)$ coincides with the first nonabelian Galois cohomology set of $G(k^s)$.

If $G\rightarrow S$ is a group scheme acting on $X\rightarrow S$ and $x \in X(T)$ is a $T$-valued point, we write $Z_G(x) \rightarrow T$ for the centraliser of $x$ in $G$. It is defined by the following pullback square:
\begin{center}
	\begin{tikzcd}
		 Z_G(x) \arrow[d] \arrow[r ] & T \arrow[d] \\
		G\times_S X \arrow[r]   & X\times_S X                        
	\end{tikzcd}
	\end{center}
Here $G\times_S X \rightarrow X \times_S X$ denotes the map $(g,x) \mapsto (g\cdot x,x)$  and $T \rightarrow X\times_S X$ denotes the composition of $x$ with the diagonal $X\rightarrow X \times_S X$. 

If $V$ is a vector space over a field $k$ we write $k[V]$ for the graded algebra $\Sym(V^{\vee})$. Then $V$ is naturally identified with the $k$-points of the scheme $\Spec k[V]$, and we call this latter scheme $V$ as well.
If $G$ is a group scheme over $k$ acting on $V$ we write $V \GIT G\coloneqq \Spec k[V]^G$ for the \define{GIT quotient} of $V$ by $G$.

\subsubsection{Root lattices}\label{subsubsection: root lattices}

We define a \define{lattice} to be a finitely generated free $\Z$-module $\Lambda$ together with a symmetric and positive-definite bilinear form $(\cdot,\cdot)\colon \Lambda\times  \Lambda \rightarrow \Z$. We write $\Lambda^{\vee}\coloneqq \{\lambda\in \Lambda \otimes \Q \mid (\lambda, \Lambda) \subset \Z\}$ for the \define{dual lattice} of $\Lambda$, which is naturally identified with $\Hom(\Lambda,\Z)$. 
We say $\Lambda$ is a \define{root lattice} if $(\lambda,\lambda)$ is an even integer for all $\lambda\in \Lambda$ and the set 
$$\{ \alpha \in \Lambda \mid (\alpha,\alpha) =2 \}$$
generates $\Lambda$. If $\Phi \subset \Real^n$ is a simply laced root system then $\Lambda= \Z\Phi$ is a root lattice. In that case we define the type of $\Lambda$ to be the Dynkin type of $\Phi$.

If $S$ is a scheme, an \define{\'etale sheaf of root lattices} $\Lambda$ over $S$ is defined as a locally constant \'etale sheaf of finite free $\Z$-modules together with a bilinear pairing $\Lambda\times \Lambda \rightarrow \Z$ (where $\Z$ denotes the constant \'etale sheaf on $S$) such that for every geometric point $\bar{s}$ of $S$ the stalk $\Lambda_{\bar{s}}$ is a root lattice. In that case $\Aut(\Lambda)$ is a finite \'etale $S$-group.

\subsubsection{Reductive groups and Lie algebras}

A reductive group scheme over $S$ is a smooth $S$-affine group scheme $G\rightarrow S$ whose geometric fibres are connected reductive groups. 
See \cite{Springer-linearalgebraicgroups} for the basics of reductive groups over a field and \cite{Conrad-reductivegroupschemes} for reductive group schemes over a general base.
A reductive group is assumed to be connected.

If $G,H,\dots$ are algebraic groups then we will use gothic letters $\lieg,\lieh,\dots$ to denote their Lie algebras.
If $G$ is a reductive group with split maximal torus $T\subset G$, we shall write $\Phi_{\liet}\subset X^*(T)$ for the set of roots of $T$ in $\lieg$, and $\Phi^{\vee}_{\liet}\subset X_*(T)$ for its set of coroots.
The map $\alpha \in \Phi_{\liet} \mapsto d\alpha\in \Hom(\liet,k)$ identifies $\Phi_{\liet}$ with the set of roots of $\liet$ in $\lieg$, and we will use this identification without further comment.

If $x$ is an element of a Lie algebra $\lieg$ then we write $\mathfrak{z}_{\lieg}(x)$ for the centraliser of $x$ in $\lieg$, a subalgebra of $\lieg$. 
We note that if $G$ is an algebraic group over a field $k$ and $x\in \lieg$ any element, then the inclusion $\Lie Z_G(x) \subset \mathfrak{z}_{\lieg}(x)$ is an equality if the characteristic of $k$ is zero or if $x$ is semisimple \cite[Proposition 1.10]{Humphreys-conjugacyclassesalgebraic}.

\begin{table}
\centering
\begin{tabular}{|c | c | c |}
	\hline
	   Symbol & Definition & Reference in paper \\
	\hline       
		$H$ & Split adjoint group of type $ADE$ & \S\ref{subsection: a stable Z/2Z-grading} \\
		$(T,P,\{X_{\alpha}\})$ & Pinning of $H$ & \S\ref{subsection: a stable Z/2Z-grading} \\
		$\theta$ & Split stable involution of $H$ & \S\ref{subsection: a stable Z/2Z-grading} \\
		$G$ & Fixed points of $\theta$ on $H$ & \S\ref{subsection: a stable Z/2Z-grading} \\
		$V$ & $(-1)$-part of action of $\theta$ on $\mathfrak{h}$ & \S\ref{subsection: a stable Z/2Z-grading}\\
		$B$ & GIT quotient $V\GIT G$ & \S\ref{subsection: a stable Z/2Z-grading}\\
		$\pi\colon V \rightarrow B$ & Invariant map & \S\ref{subsection: a stable Z/2Z-grading} \\
		$\kappa_E \colon B \rightarrow V$ & Kostant section associated to $E$ & \S\ref{subsection: Kostant sections} \\
		$C^{\circ} \rightarrow B$ & Family of affine curves & \S\ref{subsection: A family of curves}  \\
		$C \rightarrow B$ & Family of projective curves & \S\ref{subsection: A family of curves} \\
		$\infty_1,\dots,\infty_m$ & Marked points of $C\rightarrow B$ & \S\ref{subsection: A family of curves} \\
		$p_{d_1},\dots ,p_{d_r}$ & Invariant polynomials of $G$-action on $V$ &  \S\ref{subsection: A family of curves} \\
		$\kappa\colon B\rightarrow V$ & Fixed choice of Kostant section & \S\ref{subsection: The universal centraliser} \\
		$A, Z \rightarrow B^{\rs}$ & Centraliser of $\kappa|_{B^{\rs}}$ in $H$ and $G$ & \S\ref{subsection: The universal centraliser}\\
		$\Lambda \rightarrow B^{\rs}$ & Character group scheme of the torus $A\rightarrow B^{\rs}$& \S\ref{subsection: The universal centraliser}\\
		$N_{\Lambda}$ & $\image\left(\Lambda/2\Lambda\rightarrow \Lambda^{\vee}/2\Lambda^{\vee}\right)$ & \S\ref{subsection: The universal centraliser} \\
		$J^{\rs} \rightarrow B^{\rs}$ & Jacobian variety of $C^{\rs} \rightarrow B^{\rs}$ & \S\ref{subsection: A family of curves} \\
		$D$ & Zero locus of discriminant $\Delta\in \Q[B]$ & \S\ref{subsection: the discriminant locus} \\
		$B^1$ & Complement of singular locus of $D$ in $B$ & \S\ref{subsection: the discriminant locus} \\
		$\bar{J}\rightarrow B$ & Compactified Jacobian & \S\ref{subsection: the definition} \\
		$\intbigH, \intbigG, \intbigV$ & Extensions of above objects over $\Z$ &\S\ref{subsection: integral structures} \\
		$N$ & Sufficiently large integer &\S\ref{subsection: spreading out} \\
		$S$ & $\Z[1/N]$ &\S\ref{subsection: spreading out} \\
	\hline
\end{tabular}
\caption{Notation used throughout the paper}
\label{table: notation}
\end{table}

\section{Background}\label{section: background}

\subsection{The adjoint quotient of a Lie algebra}\label{subsection: the chevalley morphism}

To motivate the results in Vinberg theory, we first recall some classical results in the invariant theory of Lie algebras.

Let $H$ be a connected reductive group over a field $k$ of characteristic zero with Lie algebra $\lieh$.
The group $H$ acts on $\lieh$ via the adjoint representation. 
Let $p \colon \lieh \rightarrow \lieh \GIT H = \Spec k[\lieh]^H$ be the so-called \define{adjoint quotient} induced by the inclusion $k[\lieh]^H\subset k[\lieh]$.
We interpret $\lieh \GIT H$ as the space of invariants of the $H$-action on $\lieh$ and $p$ as the morphism of taking invariants.
Recall that an element $x\in \lieh$ is said to be \define{regular} if $\dim \mathfrak{z}_{\lieh}(x)$ is minimal among elements of $\lieh$; this minimal value equals the rank of $H$.
The subset of regular elements defines an open subscheme $\lieh^{\reg}\subset \lieh$.
The following classical proposition summarises the invariant theory of $\lieh$.

\begin{proposition}\label{proposition: classical invariant theory of H on lieh}
\begin{itemize}
\item Every semisimple element of $\lieh$ is contained in a Cartan subalgebra, and if $k$ is algebraically closed every two Cartan subalgebras are $H(k)$-conjugate.
\item Let $\mathfrak{c} \subset \lieh$ be a Cartan subalgebra and let $W = N_H(\mathfrak{c})/Z_H(\mathfrak{c})$. Then the inclusion $\mathfrak{c} \subset \lieh$ induces an isomorphism (the Chevalley isomorphism)
\begin{align*}
\mathfrak{c} \GIT W \simeq \lieh \GIT H.
\end{align*}
Since $W$ is a finite reflection group, this quotient is isomorphic to affine space.
\item If $k$ is algebraically closed and  $b\in (\lieh\GIT H)(k)$, the fibre $p^{-1}(b)$ contains a unique open $H(k)$-orbit (consisting of the regular elements with invariants $b$) and a unique closed $H(k)$-orbit (consisting of the semisimple elements with invariants $b$).
\end{itemize} 
\end{proposition}



We will often use induction arguments to reduce a statement for $\lieh$ to a reductive Lie algebra of smaller rank. 
To this end, the following lemma will be helpful.
We suppose for this lemma that $H$ is split and $T\subset H$ is a split maximal torus.
This determines a root datum $(X^*(T),\Phi_{\liet},X_*(T),\Phi^{\vee}_{\liet})$ (in the sense of \cite[\S7.4]{Springer-linearalgebraicgroups}) and a Weyl group $W = N_G(T)/T$. 
\begin{lemma}\label{lemma: centraliser semisimple element}
Let $x\in \liet$ be a semisimple element. Then the centraliser $Z_H(x)$ is a (connected) reductive group. 
Moreover, let
\begin{align*}
\Phi_{\liet}(x)  = \{ \alpha\in \Phi_{\liet} \mid \alpha(x) = 0\} \quad \text{and} \quad \Phi_{\liet}^{\vee}(x)  = \{ \alpha^{\vee} \in \Phi_{\liet}^{\vee} \mid \alpha \in \Phi_{\liet}(x)\}.
\end{align*}
Let $W_x = Z_{W}(x)$. Then the root datum of $Z_H(x)$ is $(X^*(T),\Phi_{\liet}(x),X_*(T),\Phi^{\vee}_{\liet}(x))$, and the Weyl group of $Z_H(x)$ with respect to $T$ is isomorphic to $W_x$. 
\end{lemma}
\begin{proof}
The centraliser $Z_H(x)$ is connected by \cite[Theorem 3.14]{Steinberg-Torsioninreductivegroups}. 
The fact that it is reductive and has the above root datum follows from \cite[Lemma 3.7]{Steinberg-Torsioninreductivegroups}.
The claim about the Weyl group of $Z_H(x)$ follows from \cite[Lemma 3.7(c)]{Steinberg-Torsioninreductivegroups}, again using the fact that $Z_H(x)$ is connected.
%
\end{proof}

Let $\mathfrak{c} \subset \lieh$ a Cartan subalgebra. 
The \define{discriminant polynomial} $\Delta \in k[\lieh]^H$ is the image of the product of all the roots $\prod \alpha \in k[\mathfrak{c}]^W$ with respect to $\mathfrak{c}$ under the Chevalley isomorphism $k[\mathfrak{c}]^W \xrightarrow{\sim} k[\lieh]^H$; it is independent of the choice of $\mathfrak{c}$.
For $x\in \lieh$ we have $\Delta(x) \neq 0$ if and only if $x$ is regular semisimple.
The \define{discriminant locus} (or discriminant divisor) $D\subset \lieh\GIT H$ is the zero locus of $\Delta$.
This subscheme will play a fundamental role later in this paper (in particular in Chapter \ref{section: the mildly singular locus}).

The next lemma says that the \'etale local structure of $\lieh \GIT H$ and $D$ near a point is determined by the centraliser of a semisimple lift of that point.
\begin{lemma}\label{lemma: etale local structure GIT quotient}
Let $x\in \lieh$ be a semisimple element with centraliser $\mathfrak{z}_{\lieh}(x)$. 
Let $\mathfrak{c} \subset \lieh$ be a Cartan subalgebra containing $x$.
Let $W$ and $W_x$ be the respective Weyl groups of $\lieh$ and $\mathfrak{z}_{\lieh}(x)$ with respect to $\mathfrak{c}$.
Consider the diagram:
 \begin{center}
    \begin{tikzcd}
    \mathfrak{c} \arrow[d, "\phi_x"'] \arrow[rd,"\phi"] &  \\
     \mathfrak{c} \GIT W_x \arrow[r, "\psi"]           &  \mathfrak{c} \GIT W               
    \end{tikzcd}
  \end{center}
Then $\psi$ is \'etale at $\phi_x(x)$. 
Moreover, if $D$ and $D_x$ denote the discriminant divisors of $\lieh$ and $\mathfrak{z}_{\lieh}(x)$ respectively, then $\psi^*D = D_x + R$, where $R$ is a divisor of $\mathfrak{c} \GIT W_x$ not containing $\phi_x(x)$ in its support. 
\end{lemma}
\begin{proof}
We may assume that $k$ is algebraically closed and hence that $\lieh$ is split.
Since $\phi$ and $\phi_x$ are finite and faithfully flat (they are even Galois with Galois group $W$ and $W_x$), the map $\psi$ is finite and faithfully flat.
The fact that $\psi$ is \'etale at $\phi_x(x)$ follows from the fact that the stabiliser of the $W$-action on $x$ is precisely $W_x$ (Lemma \ref{lemma: centraliser semisimple element}). 

To prove the claim about the discriminant divisors, let $\Delta = \prod_{\alpha \in \Phi_{\mathfrak{c}}} \alpha \in k[\mathfrak{c}]^W$ and $\Delta_x  = \prod_{\alpha \in \Phi_{\mathfrak{c}}(x)} \alpha \in k[\mathfrak{c}]^{W_x}$ denote the respective discriminant polynomials.
By definition of $\Phi_{\mathfrak{c}}(x)$, $\Delta = \Delta_x \cdot \mathcal{R}$ as elements of $k[\mathfrak{c}]^{W_x}$, where $\mathcal{R} \in k[\mathfrak{c}]^{W_x}$ is a polynomial that does not vanish at $x$. 
Since $D$ and $D_x$ are the zero loci of $\Delta$ and $\Delta_x$ respectively, this proves the claim.
\end{proof}


\subsection{Vinberg theory}\label{subsection: Vinberg theory}

We keep the notations from \S\ref{subsection: the chevalley morphism}.
Let $m \geq 1$ be an integer.
A \define{$\Z/m\Z$-grading} on $\lieh$ is, by definition, a direct sum decomposition
\begin{align*}
\lieh = \bigoplus_{i \in \Z/m\Z} \lieh(i)
\end{align*}
into linear subspaces satisfying $[\lieh(i),\lieh(j) ] \subset \lieh(i+j)$.
Given a $\Z/m\Z$-grading on $\lieh$, let $\lieg \coloneqq \lieh(0)$ and $V \coloneqq \lieh(1)$. 
Then $\lieg$ is a subalgebra of $\lieh$ and the restriction of the adjoint representation of $\lieh$ induces an action of $\lieg$ on $V$. 
If $\zeta \in k$ is a primitive $m$-th root of unity, giving a $\Z/m\Z$-grading amounts to giving, by considering $\zeta^i$-eigenspaces, an automorphism $\theta$ of $\lieh$ of order dividing $m$. 
In general when no such $\zeta$ exists or is fixed, giving a $\Z/m\Z$-grading amounts to giving a homomorphism $\mu_m \rightarrow \Aut(\lieh)$ of group schemes over $k$.

Let $\mu_m \rightarrow \Aut(H)$ be a morphism of group schemes. The composition $\mu_m \rightarrow \Aut(H) \rightarrow \Aut(\lieh)$ determines a $\Z/m\Z$-grading on $\lieh$. If $G$ is the identity component of the centraliser of $\mu_m$ in $H$, then $G$ has Lie algebra $\lieg$ and acts on $V = \lieh(1)$ by restriction of the adjoint action. The pair $(G,V)$ is called a \define{Vinberg representation}, and its study is dubbed \define{Vinberg theory} \cite{Vinberg-theweylgroupofgraded}.
If $\lieh$ is a semisimple $\Z/m\Z$-graded Lie algebra, a natural choice for $H$ is the adjoint group $\Aut(\lieh)^{\circ}$ of $\lieh$: this is the unique (up to nonunique isomorphism) connected semisimple group with trivial centre and Lie algebra $\lieh$. 

We now summarise some of the highlights of Vinberg theory, referring to \cite{Panyushev-Invarianttheorythetagroups,  Levy-Vinbergtheoryposchar} for proofs.
We call an element $x\in V$ semisimple, nilpotent or regular if it is so when considered as an element of $\lieh$. 
We call a subspace $\mathfrak{c} \subset V$ that consists of semisimple elements, that satisfies $[\mathfrak{c},\mathfrak{c}]=0$, and is maximal with these properties (among subspaces of $V$) a \define{Cartan subspace}.

\begin{lemma}
If $x\in V$ has Jordan decomposition $x = x_s + x_n$ where $x_s, x_n$ are commuting elements that are semisimple and nilpotent respectively, then $x_s,x_n \in V$.
\end{lemma}

\begin{proposition}
Every semisimple $x\in V$ is contained in a Cartan subspace of $V$. 
Every two Cartan subspaces are $G(k^s)$-conjugate. 
\end{proposition} 

We call a triple $(e,h,f)$ an \define{$\liesl_2$-triple} of $\lieh$ if $e,h,f$ are nonzero elements of $\lieh$ satisfying the following relations:
\begin{equation*}
[h,e] = 2e , \quad [h,f] = -2f ,\quad [e,f] = h .
\end{equation*}
The classical Jacobson--Mozorov lemma states that every nilpotent element in $\lieh$ can be completed to an $\liesl_2$-triple. 
If $\lieh$ is $\Z/m\Z$-graded, we say an $\liesl_2$-triple is \define{normal} if $e \in \lieh(1)$, $h \in \lieh(0)$ and $f\in \lieh(-1)$.

\begin{lemma}[Graded Jacobson--Mozorov]\label{lemma: graded Jacobson-Mozorov}
Every nilpotent $e\in \lieh(1)$ is contained in a normal $\liesl_2$-triple $(e,h,f)$.
\end{lemma}
\begin{proof}
    See \cite[Proposition 4]{KostantRallis-Orbitsrepresentationssymmetrisspaces}, which only treats the case $m=2$ but whose proof works for any $m$.
    (We will only need the $m=2$ case in this paper.)
\end{proof}

The next proposition describes the basic geometric invariant theory of the representation $(G,V)$.
Let $\pi\colon V\rightarrow V \GIT G = \Spec k[V]^G$ be the graded analogue of the adjoint quotient from \S\ref{subsection: the chevalley morphism}.
\begin{proposition}\label{proposition: basic invariant theory of vinberg representation full generality}
\begin{itemize}
\item Let $\mathfrak{c} \subset V$ be a Cartan subspace and $W(\mathfrak{c})  \coloneqq N_G(\mathfrak{c})/Z_G(\mathfrak{c})$. Then the inclusion $\mathfrak{c} \subset V$ induces an isomorphism 
\begin{align*}
\mathfrak{c} \GIT W(\mathfrak{c}) \simeq V\GIT G. 
\end{align*}
The group $W(\mathfrak{c})$ is a finite pseudo-reflection group, so the quotient is isomorphic to affine space.

\item If $k$ is algebraically closed and $b\in (V \GIT G)(k)$, then the fibre $\pi^{-1}(b)$ contains a unique closed $G(k)$-orbit, the set of semisimple elements with invariants $b$.
\end{itemize}
\end{proposition}

\begin{remark}
In contrast to the $m= 1$ case, it is not true in general that two regular elements $x,y \in V(k)$ with the same invariants (that is, with the same image in $V \GIT G$) are $G(k^s)$-conjugate. Indeed, it follows from Proposition \ref{proposition: H^theta acts simply transitively regular nilpotents} below that the $\Z/2\Z$-gradings introduced in \S\ref{subsection: a stable Z/2Z-grading} can have multiple $G(k^s)$-orbits of regular nilpotent elements.
\end{remark}

\subsection{Stable gradings}

We keep the notations of \S\ref{subsection: the chevalley morphism} and \S\ref{subsection: Vinberg theory}.
Of particular interest to us are the stable gradings.

\begin{definition}\label{definition: stable gradings}
Suppose that $k$ is algebraically closed.
We say a vector $v\in V$ is \define{stable} (in the sense of geometric invariant theory) if the $G$-orbit of $v$ is closed and its stabiliser $Z_G(v)$ is finite.
We say $(G,V)$ is \define{stable} if $V$ contains stable vectors.
If $k$ is not necessarily algebraically closed, we say $(G,V)$ is stable if $(G_{k^s},V_{k^s})$ is.
\end{definition}

Stable gradings of simple Lie algebras over an algebraically closed field of characteristic zero have been classified \cite[\S7.1, \S7.2]{GrossLevyReederYu-GradingsPosRank} in terms of regular elliptic conjugacy classes of (twisted) Weyl groups.
In the case of involutions (in other words, $\Z/2\Z$-gradings), this classification takes the simple form of Lemma \ref{lemma: stable involutions over alg closed field are conjugate}, see \cite[Lemma 2.6]{Thorne-thesis} for a proof. We say two involutions $\theta,\theta' \colon H\rightarrow H$ are $H(k)$-conjugate if there exists an $h\in H(k)$ such that $\theta' = \Ad(h)\circ  \theta \circ \Ad(h)^{-1}$.

\begin{lemma}\label{lemma: stable involutions over alg closed field are conjugate}
Suppose that $k$ is algebraically closed.
Then there exists a unique $H(k)$-conjugacy class of stable involutions.
\end{lemma}

For example, if $H$ is a torus then the only stable involution is given by the inversion map $h \mapsto h^{-1}$.

One of the main advantages of stable gradings is that they have a particularly good invariant theory. 
The next proposition describes this more precisely in the case of $\Z/2\Z$-gradings.
In particular, it shows that regular semisimple orbits over algebraically closed fields are well understood.
We refer to \cite[\S2]{Thorne-thesis} for precise references.

\begin{proposition} \label{proposition: invariant theory (G,V) stable involution}
	Suppose that $\theta\colon H\rightarrow H$ is a stable involution, with associated Vinberg representation $(G,V)$.
	Then the following properties are satisfied:
	\begin{enumerate}
		\item Let $\mathfrak{c} \subset V$ be a Cartan subspace and $W(\mathfrak{c})  \coloneqq N_G(\mathfrak{c})/Z_G(\mathfrak{c})$. Then $\mathfrak{c}$ is a Cartan subalgebra of $\lieh$ and the map $N_{\bigG}(\mathfrak{c}) \rightarrow W_{\mathfrak{c}} \coloneqq N_{\bigH}(\mathfrak{c})/Z_{\bigH}(\mathfrak{c})$ is surjective.
		Consequently, the inclusions $\mathfrak{c} \subset V \subset \lieh$ induce isomorphisms
		$$\mathfrak{c}\GIT W_{\mathfrak{c}} \simeq \bigV \GIT \bigG \simeq \lieh \GIT \bigH  .$$
		In particular, the quotient is isomorphic to affine space. 
		\item Suppose that $k$ is algebraically closed and let $x,y\in \bigV(k)$ be regular semisimple elements. Then $x$ is $\bigG(k)$-conjugate to $y$ if and only if $x,y$ have the same image in $\bigV\GIT \bigG$. 
		\item Let $\Delta \in \Q[\bigV]^{\bigG}$ be the restriction of the Lie algebra discriminant of $\bigh$ to the subspace $\bigV$ and suppose that $k$ is algebraically closed. Then for all $x\in \bigV(k)$, $x$ is regular semisimple if and only if $\Delta(x) \neq 0$, if and only if $x$ is stable in the sense of Definition \ref{definition: stable gradings}. 
	\end{enumerate}
\end{proposition}

\subsection{Arithmetic Invariant Theory}\label{subsection: arithmetic invariant theory}

Let $k$ be a field with separable closure $k^s$.
Let $G/k$ be a smooth algebraic group acting on a $k$-vector space $V$.
In general, a fixed $G(k^s)$-orbit in $V(k^s)$ might break up into multiple $G(k)$-orbits, and the study of this phenomenon is referred to as arithmetic invariant theory \cite{BhargavaGross-AIT}.
We recall its relation to Galois cohomology, which lies at the basis of the orbit parametrisations in this paper.

\begin{lemma}\label{lemma: AIT baby case field bhargava gross}
    Suppose that $G$ acts on a $k$-scheme $X$. Suppose that the $k$-point $e\in X(k)$ has smooth stabiliser $Z_{G}(e)$ and that the action of $G(k^s)$ on $X(k^s)$ is transitive.
    Then there is a natural bijection
    \begin{align*}
    G(k) \backslash X(k) \xleftrightarrow{1:1} \ker(\HH^1(k,Z_{G}(e))\rightarrow \HH^1(k,G)).
    \end{align*}
\end{lemma}
\begin{proof}
This is \cite[Proposition 1]{BhargavaGross-AIT}.
The bijection is explicitly constructed as follows: if $x \in X(k)$, transitivity ensures that there exists an element $g
\in G(k^s)$ with $x = g\cdot e$.
For every element $\sigma\in \Gamma_k = \Gal(k^s/k)$, we again have $x = \sigma(g) \cdot e$, so the map $\sigma \mapsto g^{-1} \sigma(g)$ defines a $1$-cocycle with values in $Z_G(e)$ which is trivial in $\HH^1(k,G)$.
\end{proof}

In fact, we will need a relative version of Lemma \ref{lemma: AIT baby case field bhargava gross} which is valid over any base scheme.

\begin{lemma}\label{lemma: AIT full generality}
Let $G\rightarrow S$ be a smooth affine group scheme acting on an $S$-scheme $X$.
Let $e\in X(S)$ be an $S$-point and suppose that the action map $m\colon G \rightarrow X, g\mapsto g\cdot e$ is smooth and surjective.
Then the assignment $x\mapsto$ `isomorphism class of the $Z_G(e)$-torsor $m^{-1}(x)$' induces a bijection between the set of $G(S)$-orbits of $X(S)$ and the kernel of the map of pointed sets $\HH^1(S,Z_{G}(e)) \rightarrow \HH^1(S,G)$.  
\end{lemma}
\begin{proof}
	This is \cite[Exercise 2.4.11]{Conrad-reductivegroupschemes}:
	the conditions imply that $X\simeq G/Z_{G}(e)$ and since $G$ and $Z_{G}(e)$ (the fibre above $e$ of the smooth map $m$) are $S$-smooth we can replace fppf cohomology by \'etale cohomology. 
\end{proof}

\subsection{The Grothendieck--Serre conjecture}\label{subsection: the grothendieck-serre conjecture}

We discuss some general results concerning principal bundles over reductive group schemes which will be useful in \S\ref{section: constructing orbits}.
Recall from \cite[Definition 3.1.1]{Conrad-reductivegroupschemes} that a reductive group scheme over $S$ is a smooth $S$-affine group scheme $G\rightarrow S$ whose geometric fibres are connected reductive groups.

\begin{definition}\label{definition: Grothendieck Serre conjecture}
Let $R$ be a regular local ring with fraction field $K$ and let $G \rightarrow \Spec R$ be a reductive group scheme.
We say that the \define{Grothendieck--Serre conjecture holds for $R$ and $G$} if the restriction map
\begin{align*}
\HH^1(R,G) \rightarrow \HH^1(K,G)
\end{align*}
is injective.
\end{definition}

Note that the injectivity $\HH^1(R,G)\rightarrow \HH^1(K,G)$ is stronger than requiring that this map has trivial kernel, since this is merely a map of pointed sets.
It is conjectured that the Grothendieck-Serre conjecture holds for every reductive group scheme over every regular local ring; see \cite{Panin-Grothendieck-Serreconjecturesurvey} for a survey and \cite[\S1.4]{Cesnavicius-GrothendieckSerreunramifiedcase} for a short summary of known results.
Below we will single out the known cases that we will need.

\begin{lemma} \label{lemma: Grothendieck Serre conjecture implies spreading out result}
Let $X$ be a regular integral scheme with function field $K$.
Let $G$ be a reductive $X$-group scheme. 
Suppose that the Grothendieck-Serre conjecture holds for all local rings of $X$ and $G$.
Then every two $G$-torsors over $X$ that are generically isomorphic (that is, isomorphic over $K$) are Zariski locally isomorphic (that is, isomorphic after restricting to a Zariski open cover). 
\end{lemma}
\begin{proof}
Let $T,T'$ be two $G$-torsors over $X$ which are generically isomorphic and let $x\in X$. 
We need to prove that $x$ has an open neighbourhood over which $T$ and $T'$ are isomorphic. 
Since the Grothendieck-Serre conjecture holds for the local ring $\O_{X,x}$, the torsors $T$ and $T'$ are isomorphic when restricted to $\Spec \O_{X,x}$. 
The result follows from spreading out this isomorphism.
\end{proof}

\begin{proposition}\label{proposition: known cases Grothendieck Serre conjecture}
Let $R$ be a regular local ring and $G$ a reductive $R$-group. Suppose that at least one of the following is satisfied:
\begin{itemize}
\item $R$ is a discrete valuation ring;
\item $R$ contains an infinite field.
\end{itemize}
Then the Grothendieck-Serre conjecture holds for $R$ and $G$.
\end{proposition}
\begin{proof}
The case of a discrete valuation ring was proved by Nisnevich \cite{Nisnevich-Espaceshomogenesprincipaux}, with corrections by Guo \cite{Guo-GrothendieckSerresemilocalDedekind}.
The case where $R$ contains an infinite field was proved by Fedorov and Panin \cite{FedorovPanin-GrothendieckSerreconjinftinitefield}.
\end{proof}

The conjecture is known in many other cases; see \cite{Cesnavicius-GrothendieckSerreunramifiedcase} for a recent general result when $R$ is of mixed characteristic and \cite{Panin-GrothendieckSerrefinitefield} for the case where $R$ contains a finite field.

\begin{corollary}\label{corollary: known cases Grothendieck Serre conjecture spreading out version}
Let $X$ be a regular integral scheme and $G$ a reductive $X$-group. Suppose that at least one of the following conditions is satisfied:
\begin{itemize}
\item $X$ is a Dedekind scheme;
\item $X$ has a map to the spectrum of an infinite field.
\end{itemize}
Then every two $G$-torsors over $X$ that are generically isomorphic are Zariski locally isomorphic. 
\end{corollary}
\begin{proof}
Combine Lemma \ref{lemma: Grothendieck Serre conjecture implies spreading out result} and Proposition \ref{proposition: known cases Grothendieck Serre conjecture}.
\end{proof}

\section{Around Thorne's thesis}\label{section: recollections of Thorne's thesis}

In the remainder of this paper we will focus on a particular Vinberg representation. Given a Dynkin diagram of type $A,D,E$, we will canonically construct a stable $\Z/2\Z$-grading on the Lie algebra of the corresponding type following Thorne's thesis \cite[\S2]{Thorne-thesis}.
We then recall and extend some of its basic properties in \S\ref{subsection: explicit determination G and V}--\ref{subsection: Kostant sections}.
In \S\ref{subsection: A family of curves} we introduce the corresponding family of curves $C\rightarrow B$ and in \S\ref{subsection: The universal centraliser} we recall the relation between stabilisers in $V$ and the $2$-torsion in the Jacobians of smooth fibres of $C\rightarrow B$.
We do not claim any originality in this chapter, except maybe for some of the calculations in \S\ref{subsection: monodromy}.

\subsection{A split stable \texorpdfstring{$\Z/2\Z$}{Z/2Z}-grading}\label{subsection: a stable Z/2Z-grading}

Let $\bigH$ be a split adjoint simple group of type $A,D,E$ over $\Q$ with Dynkin diagram $\mathsf{D}$. 
We have an exact sequence
\begin{align}\label{equation: exact sequence automorphism group H}
1\rightarrow \bigH\rightarrow \Aut(\bigH) \rightarrow \Aut(\mathsf{D}) \rightarrow 1.
\end{align}
Assume that $\bigH$ is equipped with a pinning $(\bigT,\bigP,\{X_{\alpha}\})$.
By definition, this means that $\bigT \subset \bigH$ is a split maximal torus (which determines a root system $\Phi_H \coloneqq \Phi_{\liet}$), $\bigP\subset \bigH$ is a Borel subgroup containing $\bigT$ (which determines a root basis $S_{\bigH} \subset \Phi_{\bigH}$), and $X_{\alpha}$ is a generator for each root space $\bigh_{\alpha}$ for $\alpha \in S_{\bigH}$. 
The subgroup $\Aut((\bigH, \bigT,\bigP,\{X_{\alpha}\}))\subset \Aut(\bigH)$ of elements preserving the pinning determines a splitting of the sequence \eqref{equation: exact sequence automorphism group H}.

On the other hand, if $W = N_{\bigH}(\bigT)/\bigT$ denotes the Weyl group of $\Phi_{\bigH}$, we have an exact sequence
\begin{align}\label{equation: exact sequence automorphism group Phi}
1\rightarrow W \rightarrow \Aut(\Phi_{\bigH}) \rightarrow \Aut(\mathsf{D}) \rightarrow 1.
\end{align}
We define $\vartheta \in \Aut(\bigH)$ as the unique element of $\Aut((\bigH, \bigT,\bigP,\{X_{\alpha}\}))$ whose image in $\Aut(\mathsf{D})$ under \eqref{equation: exact sequence automorphism group H} coincides with the image of $-1\in \Aut(\Phi_{\bigH})$ in $\Aut(\mathsf{D})$ under \eqref{equation: exact sequence automorphism group Phi}.
Note that $\vartheta = 1$ if and only if $-1 \in W$.

Write $\check{\rho} \in X_*(\bigT)$ for the sum of the fundamental coweights with respect to $S_{\bigH}$, characterised by the property that $(\alpha\circ \check{\rho})(t) = t$ for all $\alpha \in S_{\bigH}$. 
Let $$\bigtheta \coloneqq \vartheta \circ \Ad(\check{\rho}(-1)) = \Ad(\check{\rho}(-1)) \circ\vartheta.$$ Then $\theta$ defines an involution of $\lieh$ and thus by considering $(\pm1)$-eigenspaces it determines a $\Z/2\Z$-grading 
$$\bigh = \bigh(0) \oplus \bigh(1).$$
Let $\bigG \coloneqq (\bigH^{\bigtheta})^{\circ}$ be the identity component of the centraliser of $\bigtheta$ in $\bigH$ and let $\bigV\coloneqq \bigh(1)$. 
The space $V$ defines a representation of $\bigG$ by restricting the adjoint representation. If we write $\bigg$ for the Lie algebra of $\bigG$ then $\bigV$ is a Lie algebra representation of $\bigg = \bigh(0)$. The Vinberg representation $(\bigG,\bigV)$ is a central object of study of this paper.
Its relevance can be summarised in the following proposition, proved in \cite[Proposition 1.9]{thorne-planequarticsAIT}.

\begin{proposition}\label{proposition: characterization split stable involutions}
Up to $\bigH(\Q)$-conjugacy, $\theta$ is the unique involution of $\bigH$ with the property that $\theta$ is stable (Definition \ref{definition: stable gradings}) and the reductive group $G$ is split over $\Q$.
\end{proposition}

The first property of Proposition \ref{proposition: characterization split stable involutions} is geometric: it characterises the $H(\bar{\Q})$-conjugacy class of $\theta$.
The second property is arithmetic, and it is equivalent to requiring the existence of a regular nilpotent in $V(\Q)$.
(For the last claim, see \cite[Corollary 2.15]{Thorne-thesis}.)
Note that in our construction of $\theta$ the element $E = \sum_{\alpha \in S_H} X_{\alpha}$ is a regular nilpotent in $V(\Q)$.

Write $\bigB \coloneqq \bigV\GIT \bigG = \Spec \Q[\bigV]^{\bigG}$ and $\bigpi\colon \bigV \rightarrow \bigB$ for the natural quotient map.
We have a $\G_m$-action on $V$ given by $\lambda \cdot v = \lambda v$ and there is a unique $\G_m$-action on $B$ such that $\pi$ is $\G_m$-equivariant. 
The invariant theory of the pair $(G,V)$ is summarised in Proposition \ref{proposition: invariant theory (G,V) stable involution}. 
We additionally record the following fact concerning the smooth locus of $\pi$, which follows from the proof of \cite[Proposition 3.10]{Thorne-thesis}.

\begin{lemma}\label{lemma: smooth locus of invariant map}
Let $x\in V$ and let $d\pi_x$ be the induced map on tangent spaces $T_x V \rightarrow T_{\pi (x)}B$. 
Then $d\pi_x$ is surjective if and only if $x$ is regular. 
Consequently, the smooth locus of $\pi$ coincides with $V^{\reg}\coloneqq V\cap \lieh^{\reg}$. 
\end{lemma}
%
%

\subsection{Explicit determination of \texorpdfstring{$(G,V)$}{(G,V)}}\label{subsection: explicit determination G and V}

Using the results of \cite{Reeder-torsion} applied to the Kac diagram of $\bigtheta$ \cite[\S7.1, \S7.2]{GrossLevyReederYu-GradingsPosRank} (or an explicit description of $\theta$ in the case of classical groups), one may calculate the isomorphism class of the split group $G$ and the representation $V$ explicitly.
These results are summarised in Table \ref{table: ADE}, where we have used the following notation:
\begin{itemize}
    \item If $G$ is defined as a subgroup of $\GL_n$, then $(n)$ denotes the representation of $G$ corresponding to this embedding.
    \item In case $A_n$, $\Sym^2_0(n)$ denotes the unique codimension one $G$-subrepresentation of $\Sym^2(n)$. Equivalently, $\Sym^2_0(n)$ can be viewed as the set of self-adjoint linear maps $(n) \rightarrow (n)$ of trace zero.
    \item In case $D_{2r}$, $\Delta(\mu_2)$ denotes the image of $\mu_2$ diagonally embedded in the centre $\mu_2 \times \mu_2$ of $\SO_{2r}\times \SO_{2r}$.
    \item In case $E_6$, $\wedge^4_0(8)$ denotes the unique $42$-dimensional subrepresentation of the $\PSp_8$-representation $\wedge^4(8)$.
    \item In case $E_8$, $\Spin_{16}/\mu_2$ denotes a $\mu_2$-quotient of $\Spin_{16}$ that is not isomorphic to $\SO_{16}$; it does not seem to have a more succinct name.
\end{itemize}
We will only need these explicit identifications in the proof of Proposition \ref{proposition: cutting off cusp}.
Moreover we will calculate the component group of $H^{\theta}$ and the centre of $G$ more uniformly in \S\ref{subsection: the component group of Htheta} and \S\ref{subsection: the fundamental group of G}.

We treat the $E_7$ case as an example. 
The extended Dynkin diagram is given by:
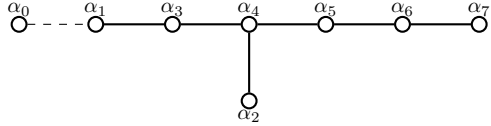
\begin{center}
\begin{tikzpicture}[transform shape, scale=.8]
\node[root] (a) {}; 
\node[root] (b) [right=of a] {}; 
\node[root] (c) [right=of b] {};
\node[root] (d) [right=of c] {};
\node[root] (e) [right=of d] {};
\node[root] (f) [right=of e] {};
\node[root] (g) [right=of f] {};
\node[root] (h) [below=of d] {};
\node [above] at (a.north) {$\alpha_0$};
\node [above] at (b.north) {$\alpha_1$};
\node [above] at (c.north) {$\alpha_3$};
\node [above] at (d.north) {$\alpha_4$};
\node [above] at (e.north) {$\alpha_5$};
\node [above] at (f.north) {$\alpha_6$};
\node [above] at (g.north) {$\alpha_7$};
\node [below] at (h.south) {$\alpha_2$};
\draw[thick] (b) -- (c) -- (d) -- (e) -- (f)--(g);
\draw[thick] (d) -- (h);
\draw[dashed] (a) -- (b);
\end{tikzpicture} 
\end{center} 
The normalised Kac coordinates of $\theta$ (given in \cite[\S7.1, Table 4]{GrossLevyReederYu-GradingsPosRank}) are everywhere zero, except at the bottom node $\alpha_2$, which has coordinate $1$. 
We may now apply the results of \cite[\S2.4]{Reeder-torsion}.
Since the Kac coordinates are invariant under the automorphism of the extended diagram, the component group of $H^{\theta}$ is of order $2$.
Since the highest root has coordinate $2$ at $\alpha_2$, the centre of $G$ is of order $2$.
If we delete the node $\alpha_2$, we obtain a diagram of type $A_7$, so $G$ semisimple is of type $A_7$.
Since $G$ is split, it follows that $G\simeq \SL_8/\mu_4$.
Moreover, the representation $V$ has highest weight the fundamental weight corresponding to $\alpha_4$, so is isomorphic to $\wedge^4(8)$, where $(8)$ denotes the defining representation of $\SL_8$.

\begin{table}
\centering
\begin{tabular}{| l | c | c | c| }
	\hline
	   Type &  $G$  & $V$  &  $\pi_0(H^{\theta})$   \\
	\hline       
	$A_{2r}$ & $\SO_{2r+1}$. & $\Sym^2_0(2r+1)$& $1$\\
	$A_{2r+1}$ &$\PSO_{2r+2}$& $\Sym^2_0(2r+2)$ & $\Z/2\Z$ \\
	$D_{2r}$ & $(\SO_{2r} \times \SO_{2r})/\Delta(\mu_2)$ & $(2r)\boxtimes (2r)$& $\Z/2\Z\times \Z/2\Z$ \\  
	$D_{2r+1}$ & $\SO_{2r+1}\times \SO_{2r+1}$ & $(2r+1)\boxtimes (2r+1) $&   $\Z/2\Z$ \\
	$E_6$ &  $\PSp_8$& $\wedge^4_0(8) $ &$1$ \\
	$E_7$ &  $\SL_8/\mu_4$& $\wedge^4(8)$ &  $\Z/2\Z$ \\
	$E_8$ &  $\Spin_{16}/\mu_2$& $\text{half spin}$ &$1$ \\
	\hline

\end{tabular}
\caption{Short description of each representation}
\label{table: ADE}
\end{table}

\subsection{The component group of \texorpdfstring{$H^{\theta}$}{Htheta}}\label{subsection: the component group of Htheta}

The group $H^{\theta}$ is typically disconnected, and we have a tautological exact sequence 
\begin{align*}
    1\rightarrow G\rightarrow H^{\theta} \rightarrow \pi_0(H^{\theta}) \rightarrow 1.
\end{align*}
The component group $\pi_0(H^{\theta})$ is a finite \'etale group scheme over $\Q$.
We will show that $\pi_0(H^{\theta})$ is split and describe it in two different ways, which will be useful in the proof of Proposition \ref{proposition: orbits corresponding to marked points are reducible}.

Firstly, we use Weyl groups.
Recall that $W_H = N_H(T)/T$ denotes the Weyl group of $H$.
We know that $T^{\theta} = T^{\vartheta}$ is a maximal torus of $G$, and moreover the centraliser $Z_H(T^{\theta})$ of $T^{\theta}$ equals $T$; these claims can be verified explicitly or follow from \cite[Lemmas 5.1 and 5.3]{Richardson-orbitsinvariantsrepresentationsinvolutions}.
It follows that $N_H(T^{\theta}) \subset N_H(T)$, so $W_{H^{\theta}} \coloneqq N_{H^{\theta}}(T^{\theta})/T^{\theta}$ is naturally a subgroup of $W_H$.
Let $W_G \coloneqq N_G(T^{\theta})/T^{\theta}$ be the Weyl group of $G$, a normal subgroup of $W_{H^{\theta}}$.
We have inclusions $W_G \subset W_{H^{\theta}} \subset W_H$.

\begin{lemma}\label{lemma: component group using weyl groups}
The inclusion $N_{H^{\theta}}(T^{\theta}) \subset H^{\theta}$ induces an isomorphism $W_{H^{\theta}}/W_G \simeq \pi_0(H^{\theta})$.
\end{lemma}
\begin{proof}
This is implicit in the proof of \cite[Lemma 3.9]{Reeder-torsion}; we sketch the details.
It suffices to prove that $H^{\theta} = G \cdot N_{H^{\theta}}(T^{\theta})$.
This can be checked on geometric points, so let $k/\Q$ be an algebraically closed field and $h\in H^{\theta}(k)$. 
The conjugate subgroup $\Ad(h)\cdot T^{\theta}$ is a maximal torus of $G_k$. 
Since $G_k$ is reductive, $G(k)$ acts transitively on its maximal tori, so $\Ad(h)\cdot T^{\theta} = \Ad(g) \cdot T^{\theta}$ for some $g\in G$. 
We see that $g^{-1}h\in N_{H^{\theta}}(T^{\theta})$, as claimed.
\end{proof}

\begin{corollary}\label{corollary: pi0 is constant}
The finite \'etale $\Q$-group $\pi_0(H^{\theta})$ is constant (in other words, has trivial Galois action) and the map $H^{\theta}(\Q)\rightarrow \pi_0(H^{\theta})$ is surjective.
\end{corollary}
\begin{proof}
It suffices to prove the latter claim.
Since $T$ is a maximal torus of $H$, $W_H$ is a constant group scheme, so its subgroup $W_{H^{\theta}}$ is constant too.
By Lemma \ref{lemma: component group using weyl groups} it suffices to show that $N_{H^{\theta}}(T^{\theta})(\Q)\rightarrow W_{H^{\theta}}$ is surjective.
This follows from Hilbert's theorem 90 since the torus $T^{\theta}$ is $\Q$-split.
\end{proof}


For the second description, choose a Cartan subspace $\mathfrak{c} \subset V$ and let $C\subset H$ be the maximal torus with Lie algebra $\mathfrak{c}$.
Since $\theta$ acts as $-1$ on $\mathfrak{c}$ it acts via inversion on $C$ hence $C[2]\subset H^{\theta}$. 
The next lemma is \cite[Proposition 1]{KostantRallis-Orbitsrepresentationssymmetrisspaces}:

\begin{lemma}\label{lemma: component group meets Cartan maximal torus}
We have $H^{\theta} = G \cdot C[2]$. In other words, the inclusion $C[2] \subset H^{\theta}$ induces a surjection $C[2] \twoheadrightarrow \pi_0(H^{\theta})$.
\end{lemma}

Lemma \ref{lemma: component group meets Cartan maximal torus} allows us to give an explicit description of $\pi_0(H^{\theta})$. 

\begin{corollary}\label{corollary: explicit description component group in terms of fundamental group}
Let $H_{sc} \rightarrow H$ be the simply connected cover of $H$ and let $\pi_1(H)$ denote the centre of $H_{sc}$. 
Then there is an isomorphism $\pi_0(H^{\theta}) \simeq \pi_1(H)/2\pi_1(H)$.
\end{corollary}
\begin{proof}
Let $C\subset H$ be a maximal torus whose Lie algebra is a Cartan subspace of $V$. (Such a torus certainly exists: take the centraliser of a regular semisimple element of $V$.)
Let $C_{sc}\subset H_{sc}$ be its preimage in $H_{sc}$. 
We have an exact sequence
\begin{align}\label{equation: sequence comparison component group with fundamental group}
1\rightarrow \pi_1(H) \rightarrow C_{sc} \rightarrow C \rightarrow 1.
\end{align}
Examining the long exact sequence associated with the $2$-torsion of \eqref{equation: sequence comparison component group with fundamental group} shows that
\begin{align*}
\frac{C[2]}{\image(C_{sc}[2] \rightarrow C[2])} \simeq  \pi_1(H)/2\pi_1(H).
\end{align*}
We claim that the left-hand-side is isomorphic to $\pi_0(H^{\theta})$. 
Indeed, the involution $\theta\colon H\rightarrow H$ uniquely extends to an involution of $H_{sc}$, still denoted by $\theta$, and a theorem of Steinberg \cite[Theorem 8.1]{Steinberg-Endomorphismsofalgebraicgroups} shows that $H_{sc}^{\theta}$ is connected. 
It follows that the induced map $H_{sc}^{\theta} \rightarrow H^{\theta}$ surjects onto $G$. 
Therefore the kernel of the natural map $C[2] \rightarrow \pi_0(H^{\theta})$ (which is surjective by Lemma \ref{lemma: component group meets Cartan maximal torus}) agrees with the image of the map $C_{sc}[2] \rightarrow C[2]$, as claimed.
\end{proof}

\subsection{The fundamental group of  \texorpdfstring{$G$}{G}}\label{subsection: the fundamental group of G}

\begin{proposition}\label{proposition: simply connected cover of $G$}
The group $G$ is semisimple and its fundamental group has order $2 \# \pi_0(H^{\theta})$.
\end{proposition}
\begin{proof}
Let $H_{sc} \rightarrow H$ be the simply connected cover of $H$ and let $\pi_1(H)$ denote the centre of $H_{sc}$.
By a previously invoked theorem of Steinberg \cite[Theorem 8.1]{Steinberg-Endomorphismsofalgebraicgroups}, $H_{sc}^{\theta}$ is connected. 
Therefore the induced map $H_{sc}^{\theta} \rightarrow G$ is surjective with kernel $\pi_1(H)[2]$. 
Moreover $\pi_1(H)[2]$ has cardinality $ \# \pi_0(H^{\theta})$ by Corollary \ref{corollary: explicit description component group in terms of fundamental group}.
Hence it suffices to prove that $H_{sc}^{\theta}$ is semisimple and its fundamental group is of order $2$.
This is a result of Kaletha, see \cite[Proposition A.1]{thorne-planequarticsAIT}.
\end{proof}

\subsection{Regular nilpotent elements in \texorpdfstring{$V$}{V}}\label{subsection: regular nilpotent elements in V}

The next proposition describes the set of regular nilpotent elements in $V$ \cite[Lemma 2.14]{Thorne-thesis}.

\begin{proposition}\label{proposition: H^theta acts simply transitively regular nilpotents}
For every field $k/\Q$, the group $H^{\theta}(k)$ acts simply transitively on the set of regular nilpotent elements of $V(k)$.
\end{proposition}

\begin{corollary}\label{corollary: bijection between pi0 and regular nilpotent orbits}
Let $k/\Q$ be a field and $E\in V(k)$ a regular nilpotent. (For example, $E = \sum_{\alpha\in S_H} X_{\alpha}$.)
Then the map $h\mapsto h\cdot E$ induces a bijection between $\pi_0(H^{\theta})$ and the set of $G(k)$-orbits of regular nilpotent elements in $V(k)$.
\end{corollary}
\begin{proof}
Follows from Proposition \ref{proposition: H^theta acts simply transitively regular nilpotents} and the fact that $H^{\theta}(k) \rightarrow \pi_0(H^{\theta})$ is surjective (Corollary \ref{corollary: pi0 is constant}).
\end{proof}

We see in particular that if $H^{\theta}$ is disconnected then there are multiple $G$-orbits of regular nilpotent elements in $V$. To state the next result, recall from \S\ref{subsection: Vinberg theory} the notion of a normal $\liesl_2$-triple.

\begin{corollary}\label{corollary: regular nilpotents can be uniquely completed to sl2 triple}
Let $k/\Q$ be a field and $E\in V(k)$ a regular nilpotent element. 
Then $E$ is contained in a unique normal $\liesl_2$-triple.
\end{corollary}
\begin{proof}
Proposition \ref{proposition: H^theta acts simply transitively regular nilpotents} shows that the stabiliser $Z_G(E)$ is trivial. 
Therefore the corollary follows from \cite[Lemma 2.17]{Thorne-thesis}.
\end{proof}

\subsection{Kostant sections}\label{subsection: Kostant sections}

We describe sections of the GIT quotient $\pi\colon \bigV \rightarrow \bigB$ whose remarkable construction is originally due to Kostant. 
Let $E \in V(\Q)$ be a regular nilpotent element and let $(E,X,F)$ be the unique normal $\liesl_2$-triple containing $E$ using Corollary \ref{corollary: regular nilpotents can be uniquely completed to sl2 triple}.
We define the affine linear subspace $\bigkappa_E \coloneqq \left(E +\mathfrak{z}_{\bigh}(F) \right) \cap \bigV \subset \bigV$.
We call $\kappa_E$ the \define{Kostant section} associated with $E$, or simply a Kostant section.
\begin{proposition}\label{proposition: Kostant section properties}
    \begin{enumerate}
        \item The composition $\bigkappa_E \hookrightarrow \bigV\rightarrow \bigB$ is an isomorphism. 
        \item $\bigkappa_E$ is contained in the open subscheme of regular elements of $\bigV$.
        \item The morphism $\bigG \times \kappa_E \rightarrow \bigV, (g,v) \mapsto g\cdot v$ is \'etale.
    \end{enumerate}
\end{proposition}
\begin{proof}
    Parts 1 and 2 are \cite[Lemma 3.5]{Thorne-thesis}; the last part is \cite[Proposition 3.4]{Thorne-thesis}, together with the fact that $G\times \kappa_E$ and $V$ have the same dimension (apply \cite[Lemma 2.21]{Thorne-thesis} to $x=0$).
\end{proof}

Every Kostant section $\kappa_E$ determines a morphism $B\rightarrow V$ that is a section of the quotient map $\pi\colon V\rightarrow B$, and we denote this section by $\kappa_E$ too.
For any $b\in B(k)$ we write $\kappa_{E,b}$ for the fibre of $\kappa_E$ over $b$.

\begin{definition}\label{definition: k-reducible orbits}
Let $k/\Q$ be a field and $v\in V(k)$.
We say $v$ is \define{$k$-reducible} if $v$ is not regular semisimple or $v$ is $G(k)$-conjugate to $\kappa_{E,b}$ for some Kostant section $\kappa_E$ and where $b = \pi(v)$. 
Otherwise, we call $v$ \define{$k$-irreducible}.
\end{definition}
If $k$ is algebraically closed, then every element of $V(k)$ is $k$-reducible by Proposition \ref{proposition: invariant theory (G,V) stable involution}.

\subsection{A family of curves}\label{subsection: A family of curves}

If $k/\Q$ is a field, an element $v\in \lieh(k)$ is called \define{subregular} if $\dim \mathfrak{z}_{\lieh}(x) = r+2$, where $r$ is the rank of $\lieh$.
By \cite[Proposition 2.27]{Thorne-thesis}, the vector space $V$ contains subregular nilpotent elements; let $e\in V(\Q)$ be such an element and fix $(e,x,f)$ a normal $\liesl_2$-triple extending it, using Lemma \ref{lemma: graded Jacobson-Mozorov}.

Slodowy \cite{Slodowy-simplesingularitiesalggroups} has shown that the restriction of the invariant map $(e+ \mathfrak{z}_{\lieh}(f)) \rightarrow B$ is a family of surfaces.
Moreover, he has shown that this family is a semi-universal deformation of its central fibre, which is a simple surface singularity of the type corresponding to that of $H$.
Proposition \ref{proposition: first properties of the family of curves} is a $\Z/2\Z$-graded analogue of Slodowy's result, due to Thorne.
Define  $C^{\circ} \coloneqq (e+ \mathfrak{z}_{\lieh}(f))\cap V$.
Restricting the invariant map $\pi\colon V\rightarrow B$ to $C^{\circ}$ defines a morphism $\varphi\colon C^{\circ}\rightarrow B$.

\begin{proposition}\label{proposition: first properties of the family of curves}
    \begin{enumerate}
    \item  The geometric fibres of $\varphi$ are reduced connected curves.
    \item  The central fibre $C_0^{\circ} = \varphi^{-1}(0)$ has a unique singular point which is a simple singularity of type $A_r, D_r, E_r$, corresponding to that of $H$.
    \item  We can choose coordinates $p_{d_1},\dots,p_{d_r}$ on $B$ where $p_{d_i}$ is homogeneous of degree $d_i$ and coordinates $(x,y,p_{d_1},\dots,p_{d_{r-1}})$ on $C^{\circ}$ such that $C^{\circ}\rightarrow B$ is given by the affine equation of Table \ref{table: introduction different cases}.
    \item  The formal completion of $C^{\circ}\rightarrow B$ along its central fibre defines a morphism of formal schemes $\widehat{C^{\circ}}\rightarrow \widehat{B}$ which is a semi-universal deformation of its central fibre.
    \item  The morphism $\varphi$ is faithfully flat. It is smooth at $x\in C^{\circ}$ if and only if $x$ is a regular element of $V$.
    \item  The action map $G\times C^{\circ} \rightarrow V, (g,x) \mapsto g\cdot x$ is smooth.
    \end{enumerate}
\end{proposition}
\begin{proof}
This is proved in Thorne's thesis. 
The first three parts are \cite[Theorem 3.8]{Thorne-thesis}; for the definition of a simple curve singularity, see \cite[End of \S2]{Thorne-actualthesis}.
The fourth part follows from the fact that the semi-universal deformation of an isolated hypersurface singularity can be explicitly computed \cite[\S2.4]{Slodowy-simplesingularitiesalggroups} and agrees with the equations given in Table \ref{table: introduction different cases}.
The last two parts are contained in \cite[Proposition 3.4 and Proposition 3.10]{Thorne-thesis}.
\end{proof}

The next lemma describes the singularities of the fibres of $C^{\circ} \rightarrow B$ very precisely; see \cite[Corollary 3.16]{Thorne-thesis} for its proof.

\begin{lemma}\label{lemma: singularities of C_b using dynkin diagram}
Let $k/\Q$ be a field, $b\in B(k)$ and $v \in V_b(k)$ a semisimple element. 
Then there is a bijection between the connected components of the Dynkin diagram of $Z_H(v)$ and the singularities of $C_b^{\circ}$, which takes each (connected, simply laced) Dynkin diagram to a singularity of the corresponding type. 
\end{lemma}

We compactify the flat affine family of curves $C^{\circ} \rightarrow B$ to a flat projective family of curves $C\rightarrow B$ as described in \cite[Lemma 4.9]{Thorne-thesis}.
That lemma implies that the complement $C \setminus C^{\circ}$ is a disjoint union of sections $\infty_1,\dots,\infty_m \colon B\rightarrow C$ and $C\rightarrow B$ is smooth in a Zariski open neighbourhood of these sections.
For every field $k/\Q$ and $b\in B(k)$, the curve $C_b$ has $k$-rational points $\infty_{1,b}, \dots, \infty_{m,b}\in C_b(k)$; we call these the \define{marked points of $C_b$}.

\begin{lemma}\label{lemma: bijection between marked points, irreducible components and certain nilpotents}
There are natural bijections between:
\begin{enumerate}
    \item The sections $\infty_1,\dots,\infty_m$ of $C\rightarrow B$;
    \item Irreducible components of $C_0$;
    \item $G$-orbits of regular nilpotent elements of $V$ whose closure contains $e$. 
\end{enumerate}
The bijections are given as follows: given a section $\infty_i$, map it to the irreducible component containing $\infty_{i,0}\in C_0$; given an irreducible component of $C_0$, map it to the $G$-orbit of any point on its smooth locus.
\end{lemma}
\begin{proof}
See \cite[Lemma 4.14]{Thorne-thesis} and its proof.
\end{proof}

For the remainder of this paper, we fix a section $\infty_1 = \infty$ of $C\rightarrow B$ and a regular nilpotent element $E\in V(\Q)$ whose $G$-orbit corresponds to $\infty$ under Lemma \ref{lemma: bijection between marked points, irreducible components and certain nilpotents}.
Moreover, we fix a choice of polynomials $p_{d_1},\dots,p_{d_r} \in \Q[V]^{G}$ and coordinates $x,y$ of $C^{\circ}$ satisfying the conclusions of Proposition \ref{proposition: first properties of the family of curves}. Recall that we have defined a $\G_m$-action on $B$ which satisfies $\lambda \cdot p_{d_i} = \lambda^{d_i}p_{d_i}$. 
There exist unique positive integers $a,b$ such that $\lambda \cdot (x,y,p_{d_1},\dots,p_{d_{r-1}}) := (\lambda^a x,\lambda^b y,\lambda^{2d_1}p_{d_1},\dots,\lambda^{2d_{r-1}}p_{d_{r-1}})$ defines a $\G_m$-action on $C$ and such that the morphism $C\rightarrow B$ is $\G_m$-equivariant with respect to the square of the usual $\G_m$-action on $B$. (The integers $(a,b)$ are given by $(w_r,w_{r+1})$ in the table of \cite[Proposition 3.6]{Thorne-thesis}. These weights can also be defined Lie theoretically, but we will not need this fact in what follows.)

\subsection{Universal centralisers}\label{subsection: The universal centraliser}

Recall from the last paragraph of \S\ref{subsection: A family of curves} that we have fixed a regular nilpotent $E\in V(\Q)$; let $\kappa\colon B\rightarrow V$ be the Kostant section corresponding to $E$ constructed in \S\ref{subsection: Kostant sections}.
Recall from our conventions in \S\ref{subsection: notation} that if $v\colon S\rightarrow V$ is an $S$-point of $V$ then $Z_G(v)\rightarrow S$ denotes the centraliser of $v$ in $G$. 

\begin{definition}\label{definition: universal centralisers}
Let $Z\rightarrow B$ be the centraliser $Z_G(\kappa)$ of the Kostant section $\kappa\colon B\rightarrow V$ with respect to the $G$-action on $V$.
Similarly, let $A\rightarrow B$ be the centraliser $Z_H(\kappa)$ of $\kappa\colon B \rightarrow \lieh$ with respect to the $H$-action on $\lieh$.
\end{definition}
For every field $k/\Q$ and $b\in B(k)$, the group scheme $Z_b$ (respectively $A_b$) is the centraliser $Z_G(\kappa_b) \subset G$ of $\kappa_b$ in $G$ (respectively $Z_H(\kappa_b) \subset H$). 
We have $Z = A \cap (G\times B)$.
Since $\kappa$ lands in the regular locus of $V$, $A$ and $Z$ are commutative group schemes.
To state the next lemma, recall that $V^{\reg}\subset V$ denotes the open subscheme of regular elements and that $\pi\colon V\rightarrow B$ and $p \colon \lieh \rightarrow B$ denote the morphisms of taking invariants.

\begin{lemma}\label{lemma: Z and A are descents of regular centralisers}
Let $v\colon S\rightarrow V^{\reg}$ be a morphism with $b = \pi(v) \in B(S)$. 
Then there is a canonical isomorphism $Z_G(v) \simeq Z_b$.
Similarly if $v\colon S\rightarrow \lieh^{\reg}$ is a morphism with invariants $b = p(v) \in B(S)$ then there is a canonical isomorphism $Z_H(v) \simeq A_b$. 
\end{lemma}
\begin{proof}
The isomorphism $Z_G(v) \simeq Z_b$ follows from \cite[Proposition 4.1]{Thorne-thesis} and a very similar proof works for $A$; we briefly sketch it.
The morphism $H\times B \rightarrow \lieh^{\reg} , (h,b) \mapsto h\cdot \kappa_b$ is smooth and surjective \cite[Lemma 3.3.1]{Riche-KostantSectionUniversalCentralizer}, so has sections \'etale locally. 
It follows that $v$ is $H$-conjugate to $\kappa_b$ \'etale locally on $S$. 
Conjugating defines isomorphisms $Z_H(v) \simeq A_b$, again \'etale locally on $S$.
Since $A_b$ is commutative, these isomorphisms do not depend on the choice of element by which we conjugate $v$ to $\kappa_b$. 
Using \'etale descent, these isomorphisms glue to give an isomorphism of group schemes $Z_H(v) \simeq A_b$.
\end{proof}

The next lemma gives a useful description of the fibres of $Z\rightarrow B$.

\begin{lemma}\label{lemma: calculation stabiliser regular element}
Let $k/\Q$ be a field and $x\in V(k)$ a regular element, with Jordan decomposition $x = x_s+x_n$.
Let $\mathfrak{c}\subset V$ be a Cartan subspace containing $x_s$ and let $C\subset H$ denote the maximal torus with Lie algebra $\mathfrak{c}$. 
Let $H_{sc}\rightarrow H$ be the simply connected cover of $H$ and $C_{sc}\rightarrow C$ its restriction to $C$.
Then there is a canonical isomorphism
\begin{align}
\Hom(Z_G(x),\F_2) \simeq \image\left( \frac{X^*(C)}{2X^*(C)+\Z\Phi_{\mathfrak{c}}(x)} \rightarrow   \frac{X^*(C_{sc})}{2X^*(C_{sc})+\Z\Phi_{\mathfrak{c}}(x)} \right).
\end{align}
\end{lemma}
\begin{proof}
A theorem of Steinberg \cite[Theorem 8.1]{Steinberg-Endomorphismsofalgebraicgroups} shows that $(H_{sc})^{\theta}$ is connected. 
Therefore $Z_G(x) =\image( Z_{(H_{sc})^{\theta}}(x) \rightarrow Z_{H^{\theta}}(x))$.
Now use \cite[Corollary 2.9]{Thorne-thesis}.
\end{proof}

Let $B^{\rs}$ denote the image of the subscheme regular semisimple elements in $V$ under $\pi\colon V\rightarrow B$.
Then $B^{\rs}$ is also the complement of the discriminant locus $(\Delta=0)$ in $B$, by Part 3 of Proposition \ref{proposition: invariant theory (G,V) stable involution}.
For a $B$-scheme $X$, we denote its restriction to $B^{\rs}$ by $X^{\rs}$. 
The group scheme $Z^{\rs} \rightarrow B^{\rs}$ is finite \'etale and $A^{\rs} \rightarrow B^{\rs}$ is a family of maximal tori. 
\begin{definition}
Let $\Lambda \rightarrow B^{\rs}$ be the character group of $A^{\rs}$.
\end{definition}
In other words, $\Lambda$ is the Cartier dual $\Hom(A^{\rs},\G_m)$ of $A^{\rs}$.
The $B^{\rs}$-scheme $\Lambda$ is an \'etale sheaf of root lattices in the sense of \S\ref{subsubsection: root lattices}.
In particular, it comes equipped with a pairing $\langle \cdot ,\cdot\rangle \colon \Lambda \times \Lambda\rightarrow \Z$. 
This pairing induces an alternating pairing $(\cdot,\cdot) \colon \Lambda/2\Lambda\times \Lambda/2\Lambda \rightarrow \F_2$ which might be degenerate.
Setting $N_{\Lambda} \coloneqq \image(\Lambda/2\Lambda\rightarrow \Lambda^{\vee}/2\Lambda^{\vee})$, we see \cite[Lemma 2.11]{Thorne-thesis} that $(\cdot,\cdot)$ descends to a nondegenerate pairing on $N_{\Lambda}$.
Lemma \ref{lemma: calculation stabiliser regular element} implies:

\begin{lemma}\label{lemma: canonical iso stabiliser and root lattice}
There exists a canonical isomorphism $Z^{\rs} \simeq N_{\Lambda}$.
\end{lemma}

We use the isomorphism of Lemma \ref{lemma: canonical iso stabiliser and root lattice} to transport the pairing from $N_{\Lambda}$ to $Z^{\rs}$: we thus obtain a nondegenerate pairing $Z^{\rs} \times Z^{\rs} \rightarrow \F_2$.

It follows from Lemma \ref{lemma: singularities of C_b using dynkin diagram} that the restriction $C^{\rs} \rightarrow B^{\rs}$ is a family of smooth projective curves; write $J^{\rs} \rightarrow B^{\rs}$ for the relative Jacobian of the family of smooth projective curves $C^{\rs} \rightarrow B^{\rs}$ \cite[\S9.3; Theorem 1]{BLR-NeronModels}.
The next result is one of the main results of Thorne's thesis and a first step towards relating the curves $C^{\rs} \rightarrow B^{\rs}$ to the representation $(G,V)$.

\begin{proposition}\label{proposition: isomorphism 2-torsion and Kostant section centraliser thorne}
There exists a canonical isomorphism $J^{\rs}[2] \simeq Z^{\rs}$ of finite \'etale group schemes that sends the Weil pairing on $J^{\rs}[2]$ to the pairing on $Z^{\rs}$ defined above.
\end{proposition}
\begin{proof}
Since both group schemes are finite \'etale and $B^{\rs}$ is normal, it suffices to prove the statement above the generic point of $B^{\rs}$ by \cite[Tag \href{https://stacks.math.columbia.edu/tag/0BQM}{0BQM}]{stacksproject}.
In that case the statement follows from \cite[Corollary 4.12]{Thorne-thesis}.
\end{proof}
Part 2 of Proposition \ref{proposition: consequences of monodromy} implies that the isomorphism $J^{\rs}[2] \simeq Z^{\rs}$ is unique.

\subsection{Monodromy of \texorpdfstring{$J^{\rs}[2]$}{Jrs[2]}}\label{subsection: monodromy}

We give some additional properties of the group scheme $J^{\rs}[2]
 \rightarrow \bigB^{\rs}$, which by Lemma \ref{lemma: canonical iso stabiliser and root lattice} and Proposition \ref{proposition: isomorphism 2-torsion and Kostant section centraliser thorne} we may identify with $N_{\Lambda}\rightarrow B^{\rs}$. 
Before we state them, we recall some definitions and set up notation. 

Recall from \S\ref{subsection: a stable Z/2Z-grading} that $T$ is a split maximal torus of $H$ with Lie algebra $\liet$ and Weyl group $W$. 
Let $L \coloneqq X^*(T)$ be its character group and $N_L \coloneqq \image(L/2L \rightarrow L^{\vee}/2L^{\vee})$.
Consider the composition 
 $\liet \rightarrow \liet \GIT W \xrightarrow{\sim} \lieh \GIT H \xrightarrow{\sim} V \GIT G = B,$
where $\liet \rightarrow \liet \GIT W$ is the natural projection, $\liet \GIT W \xrightarrow{\sim} \lieh \GIT H$ the Chevalley restriction isomorphism (Proposition \ref{proposition: classical invariant theory of H on lieh}), and $\lieh \GIT H \xrightarrow{\sim} V \GIT G$ is the isomorphism induced from the inclusion $V\subset \lieh$ (Proposition \ref{proposition: invariant theory (G,V) stable involution}).
Restricting to regular semisimple elements defines a finite \'etale cover $f\colon\liet^{\rs} \rightarrow \bigB^{\rs}$ with Galois group $W$.

\begin{proposition}\label{proposition: monodromy of J[2]}
    The finite \'etale group scheme $J^{\rs}[2]\rightarrow B^{\rs}$ becomes trivial after the base change $f\colon\liet^{\rs} \rightarrow \bigB^{\rs}$, where it becomes isomorphic to the constant group scheme $N_L$. The monodromy action is induced by the natural action of $W$ on $L$.
\end{proposition}
\begin{proof}
    Since $J^{\rs}[2]$ is isomorphic to $N_{\Lambda}$, it suffices to prove that the torus $A \rightarrow \bigB^{\rs}$ is isomorphic to the constant torus $T\times \liet^{\rs} \rightarrow \liet^{\rs}$ after pulling back along $f$, with monodromy given by the action of $W$ on $T$.  
    
    To prove this, note that by Lemma \ref{lemma: Z and A are descents of regular centralisers}, if $x\colon S \rightarrow \lieh^{\rs}$ is an $S$-point with invariants $b=p(x)\in B^{\rs}(S)$, then $Z_{\bigH}(x) \simeq A_b$ as group schemes over $S$. (Here $\lieh^{\rs}\subset \lieh$ denotes the subset of regular semisimple elements.)
    In particular, we can apply this to the $\liet^{\rs}$-point $i\colon \liet^{\rs} \rightarrow \lieh^{\rs}$ (where $i$ is the inclusion map), giving an isomorphism $T\times{\liet^{\rs}} \simeq A_{\liet^{\rs}} $.
    Since this isomorphism is induced by \'etale locally conjugating $i$ to $\kappa$ by elements of $\bigH$, the monodromy action is indeed given by the natural action of $W$ on $\bigT$. 
    \end{proof}
    
Proposition \ref{proposition: monodromy of J[2]} shows that it suffices to understand the $W$-action on $N_L$ if we wish to understand the group scheme $J^{\rs}[2]$.
To this end, we perform some root system calculations in the following proposition.

\begin{proposition}\label{proposition: consequences of monodromy}
Suppose that $H$ is not of type $A_1$.
\begin{enumerate}
\item $L^{\vee}/2L^{\vee}$ has no nonzero $W$-invariant elements.
\item The $\F_2[W]$-representation $N_L$ is absolutely irreducible, so every $W$-equivariant automorphism $N_L\rightarrow N_L$ is the identity.
\item There exists an element $w \in W$ that has no nonzero fixed points on $N_L$.
\end{enumerate}
\end{proposition}
\begin{proof}
\begin{enumerate}
\item Note that $L^{\vee}/2L^{\vee} = \Hom(L/2L,\F_2)$, so let $f \colon L/2L\rightarrow \F_2$ be a nonzero $W$-invariant functional.
If $f$ vanishes on a root of $L$, then $f$ vanishes on all of them since they form a single $W$-orbit.
Since the roots of $L$ generate $L/2L$, it follows that $f(\alpha) =1$ for every root. 
Since we have assumed that $L$ is not of type $A_1$, there exists roots $\alpha,\beta$ such that $\alpha+\beta$ is also a root. 
But then we would have $1 = f(\alpha+ \beta) = f(\alpha) + f(\beta) = 1+1 = 0$. This is a contradiction hence no such nonzero $f$ exists.
%
%
%
\item Let $S\subset N_{L}\otimes \bar{\F}_2$ be a $W$-stable subspace, and assume that $v\in S$ is a nonzero element. Since the pairing $(\cdot,\cdot)$ on $N_L$ is nondegenerate, there exists a root $\alpha \in L$ such that $(v,\alpha) \neq 0$ in $\bar{\F}_2$.
If $w_{\alpha}\in W$ denotes the reflection associated with $\alpha$, then $w_a(v) = v-(v,\alpha)\alpha$ also lies in $S$.
It follows that $w_a(v)-v = (v,\alpha)\alpha$ lies in $S$ hence $\alpha$ also lies in $S$. Since $W$ acts transitively on the roots, every root is contained in $S$. Since the roots generate $L$, it follows that $S = N_L\otimes \bar{\F}_2$, as claimed. 
The second claim follows from Schur's lemma.
(We thank Beth Romano for helping us with the proof of this fact.)
\item We first consider the case that the pairing on $L/2L$ is nondegenerate, which is equivalent to the projection map $L/2L \rightarrow N_L$ being injective.
Since $L/2L$ and $L^{\vee}/2L^{\vee}$ have the same order, the latter statement is also equivalent to the fact that $L/2L\rightarrow L^{\vee}/2L^{\vee}$ is an isomorphism.
We show that in this case it suffices to take a Coxeter element $w_{cox}$ of $W$.
Indeed, let $H_{sc} \rightarrow H$ be the simply connected cover of $H$, let $\pi_1$ be the centre of $H_{sc}$ and let $T_{sc}$ be the preimage of $T$ in $H_{sc}$.
It is a classical fact that the inclusion $\pi_1\subset T_{sc}$ restricts to an equality $\pi_1 = T_{sc}^{w_{cox}}$, see \cite[Theorem 1.6]{DwyerWilkerson-centerscoxeterelements}.\footnote{The reference assumes that $T$ is a maximal torus in a compact Lie group, but this implies the corresponding result for a maximal torus in a semisimple group over a field of characteristic zero.}
Taking $2$-torsion implies that $\pi_1[2]  = T_{sc}[2]^{w_{cox}}$.
Since the map $L/2L\rightarrow L^{\vee}/2L^{\vee}$ is an isomorphism, the same is true for the map $T_{sc}[2]\rightarrow T[2]$ which has kernel $\pi_1[2]$, hence $T_{sc}[2]^{w_{cox}} = \pi_1[2] = \{1\}$.
Since $T_{sc}[2] \simeq L/2L$, we have shown that $(L/2L)^{w_{cox}} = N_L^{w_{cox}} = 0$, as claimed.

We now consider the general case. Let $S$ be a root basis of $L$, which is an $\F_2$-basis of the vector space $L/2L$.
Since $L/2L \rightarrow N_L$ is surjective, there exists a subset $S_M \subset S$ projecting onto a basis of $N_L$. 
Let $M$ be the $\F_2$-span of $S_M$. Then $M$ is a (possibly reducible) root lattice associated with the sub-root system generated by $S_M$, and the composition $M/2M \hookrightarrow L/2L \twoheadrightarrow N_L$ is an isomorphism.
The pairing $(\cdot,\cdot)$ on $L/2L$ restricts to a pairing on $M/2M$, and the previous sentence shows that this pairing on $M/2M$ is nondegenerate.
It follows that $M$ is a direct sum of irreducible root lattices of the form considered in the first case of this proof.
Let $w$ be a Coxeter element with respect to $S_M$, i.e. a product of the simple reflections in $S_M$.
Then $(M/2M)^{w} = 0$ by the first case of the proof, so $N_L^w=0$ too.
\end{enumerate}
\end{proof}


\section{The mildly singular locus}\label{section: the mildly singular locus}

We keep the notations from Section \ref{section: recollections of Thorne's thesis}.
Recall from \S\ref{subsection: The universal centraliser} that $B^{\rs} \subset B$ denotes the locus where the discriminant polynomial $\Delta$ is nonzero, and that the family of curves $C\rightarrow B$ is smooth exactly above $B^{\rs}$ (Lemma \ref{lemma: singularities of C_b using dynkin diagram}).
In this chapter we introduce an open subset $B^1 \subset B$ strictly containing $B^{\rs}$ where we allow the fibres of $C\rightarrow B$ to have one nodal singular point. 
We therefore call $B^1\setminus B^{\rs}$ the `mildly singular locus' of $B$.
We then extend some results concerning the representation $V$ and the family of curves from $B^{\rs}$ to $B^1$ in \S\ref{subsection: representation theory over B1} and \S\ref{subsection: the family of curves over B 1}, and generalise Proposition \ref{proposition: isomorphism 2-torsion and Kostant section centraliser thorne} in \S\ref{subsection: summary}.
This will be useful for the construction of orbits in \S\ref{section: constructing orbits} and for the analysis of integral orbits of square-free discriminant in \S\ref{subsection: square-free discriminant case}.
To avoid stating the same assumption repeatedly, we will make the following assumption throughout the rest of this paper:

\begin{convention}\label{convention: H not of type A1}
The group $H$ is not of type $A_1$.
\end{convention}

\subsection{The discriminant locus}\label{subsection: the discriminant locus}

Recall that we have fixed a maximal torus $T\subset H$ in \S\ref{subsection: a stable Z/2Z-grading}.
Recall from \S\ref{subsection: the chevalley morphism} that the discriminant polynomial $\Delta \in \Q[\lieh]^H$ is the image of $\prod_{\alpha \in \Phi_{\liet}} \alpha\in \Q[\liet]^W$ under the isomorphism $\Q[\liet]^W \xrightarrow{\sim} \Q[\lieh]^{H}$ of the Chevalley restriction theorem.
Using the isomorphism $\Q[\lieh]^H\xrightarrow{\sim} \Q[V]^G=\Q[B]$ from Proposition \ref{proposition: invariant theory (G,V) stable involution}, we view $\Delta$ as an element of $\Q[B]$.

\begin{lemma}\label{lemma: discriminant irreducible}
	For every field $k/\Q$, $\Delta$ is irreducible in $k[B]$.
\end{lemma}
\begin{proof}
	It suffices to prove that we cannot partition $\Phi_{\liet}$ into two nonempty $W$-invariant subsets. 
	Equivalently, we need to prove that $W$ acts transitively on $\Phi_{\liet}$.
	This is true since $\Phi_{\liet}$ is irreducible and simply laced.
\end{proof}

We write $D$ for the subscheme of $B$ cut out by $\Delta$. Lemma \ref{lemma: discriminant irreducible} implies:

\begin{corollary}
The scheme $D$ is geometrically integral.
\end{corollary}

Write $D_{sing}$ for the singular locus of $D$, a closed subscheme of $D$. 

\begin{definition}
We define $B^1$ as the complement of $D_{sing} $ in $B$.
We define $D^1$ as the complement of $D_{sing}$ in $D$; we call $D^1$ the \define{mildly singular locus}.
\end{definition}
The subscheme $B^1\subset B$ is open and we have inclusions $B^{\rs} \subsetneq B^1 \subsetneq B$.
Since $D$ is geometrically integral, the complement of $B^1$ in $B$ has codimension $\geq 2$.
(In fact, Lemma \ref{lemma: points of D1 are those points with derived subgroup A1} shows that it has codimension exactly $2$.)
As a general piece of notation, if $X$ is a $B$-scheme we write $X^1$ for its restriction to $B^1$.

\subsection{Representation theory over \texorpdfstring{$B^1$}{B1}}\label{subsection: representation theory over B1}


If $b$ is a point of $B$ we write $\lieh_b$ for the fibre of the adjoint quotient $\lieh \xrightarrow{p} \lieh \GIT  H = B$ along this point.

\begin{lemma}\label{lemma: points of D1 are those points with derived subgroup A1}
Suppose that $k/\Q$ is an algebraically closed field and $b\in B(k)$. Then $b \in D^1(k)$ if and only if some (equivalently, every) semisimple $x\in \lieh_b(k)$ has the property that the derived subgroup of $Z_H(x)$ is of type $A_1$.
\end{lemma}
\begin{proof}
Since every two semisimple elements in $\lieh_b(k)$ are $H(k)$-conjugate (Proposition \ref{proposition: classical invariant theory of H on lieh}), requiring the last claim for some semisimple element of $\lieh_b(k)$ is equivalent to requiring it for all of them.
Let $T\subset H$ be the fixed maximal torus of \S\ref{subsection: a stable Z/2Z-grading}, with root system $\Phi_{\liet}$ and Weyl group $W$.
Let $x$ be an element of $\liet$ with invariants $b$.
By Lemma \ref{lemma: centraliser semisimple element}, $Z_H(x)$ is a reductive group with root system $\Phi_{\liet}(x) = \{ \alpha \in \Phi_{\liet} \mid \alpha(x) = 0\}$ and its Weyl group $W_x$ is the subgroup of $W$ generated by the reflections through $\Phi_{\liet}(x)$.

To prove the lemma, we need to prove that $b\in D^1(k)$ if and only if $\Phi_{\liet}(x)$ is of type $A_1$. 
Let $b_x $ be the image of $x$ in $\liet\GIT W_x$ and let $D_x\subset \liet \GIT W_x$ be the discriminant locus of $Z_H(x)$, with smooth locus $D_x^1$.
By Lemma \ref{lemma: etale local structure GIT quotient}, $D_x\rightarrow D$ is etale at $b_x$ hence $b \in D^1(k)$ if and only if $b_x \in D_x^1(k)$. 
So it suffices to prove that $b_x \in D_x^1(k)$ if and only if $\Phi_{\liet}(x)$ is of type $A_1$.

Firstly, suppose that $\Phi_{\liet}(x) = \{\alpha,-\alpha\}$ is of type $A_1$, so $W_x = \{1,w_{\alpha}\}$ is generated by the reflection through $\alpha$. 
Then one can compute $D_x$ explicitly: $\liet \GIT W_x$ is given, up to taking a product with an affine space, by the quotient of $\Spec k[X]$ by the $\Z/2\Z$-action $X\mapsto -X$.
This quotient is $\Spec k[X^2]$ and $D_x$ is the vanishing locus of $X^2$, hence $D_x$ is smooth.
Therefore $b_x \in D_x^1(k)$, which proves one direction.

Conversely, suppose that $\Phi_{\liet}(x)$ is not of type $A_1$. If $\Phi_{\liet}(x)$ were empty then $b\in B^{\rs}(k)$, so $\Phi_{\liet}(x)$ is nonempty and of rank $\geq 2$.
We need to prove that $D$ is singular at $b$. 
Since the singular locus of $D$ is closed and $x$ is the specialisation of a point $y$ for which the rank of $\Phi_{\liet}(y)$ is exactly $2$, we may assume that $\Phi_{\liet}(x)$ is either of type $A_2$ or $A_1\times A_1$.
In both cases, one can compute explicitly that $D_x$ is not smooth at $b_x$, as required. 
\end{proof}

Recall from Proposition \ref{proposition: isomorphism 2-torsion and Kostant section centraliser thorne} that if $k/\Q$ is a field and $b\in B^{\rs}(k)$, we have an isomorphism $Z_b \simeq J_b[2]$ of finite \'etale $k$-groups. 
So if $g$ denotes the common arithmetic genus of the curves $C_b$, the group scheme $Z_b$ has order $2^{2g}$.

\begin{lemma}\label{lemma: centraliser has order 2^2g-1 over D1}
If $b\in D^1(k)$, the group scheme $Z_b$ has order $2^{2g-1}$. 
\end{lemma}
\begin{proof}
The group scheme $Z_b$ is the centraliser of the element $\kappa_b$, which is regular by Proposition \ref{proposition: Kostant section properties}.
By Lemma \ref{lemma: calculation stabiliser regular element}, it suffices to prove that if $L$ is a root lattice of the same type as $H$, $N_L = \image(L/2L\rightarrow L^{\vee}/2L^{\vee})$ and $\alpha \in L$ is a root, then $\alpha$ is nonzero in $N_L$.
Since $L$ is not of type $A_1$ (Convention \ref{convention: H not of type A1}), there exists a root $\beta$ with $(\alpha,\beta) = -1$.
Therefore $\alpha \not\in 2L^{\vee}$, so $\alpha$ is nonzero in $N_L$, as claimed.
\end{proof}


Before we state the last result of this section, we record a useful lemma. 
\begin{lemma}\label{lemma: centre of centraliser semisimple element connect}
    Let $k/\Q$ be a field and $x\in \lieh(k)$ a semisimple element with centraliser $L\coloneqq Z_H(x)$.
    Then the centre of $L$ is connected.
\end{lemma}
\begin{proof}
    We may assume that $k$ is algebraically closed and that $x$ lies in $\liet(k)$, the Lie algebra of the maximal torus $T\subset H$ fixed in \S\ref{subsection: a stable Z/2Z-grading}.
    It suffices to prove that the character group of the centre of $L$ is torsion-free.
    By Lemma \ref{lemma: centraliser semisimple element} this group can be identified with $X^*(T)/\Z\Phi_{\liet}(x)$.
    Since $H$ is adjoint, $X^*(T)= \Z\Phi_{\liet}$. 
    The definition of the root system $\Phi_{\liet}(x)$ shows that it is $\Q$-closed in the sense of \cite[\S3.5]{Slodowy-simplesingularitiesalggroups}.
    By \cite[Proposition 3.5]{Slodowy-simplesingularitiesalggroups}, every root basis of $\Phi_{\liet}(x)$ can be extended to a root basis of $\Phi_{\liet}$. 
    This implies that $\Z\Phi_{\liet}(x)$ is a direct summand of $\Z\Phi_{\liet}$ so the quotient $\Z\Phi_{\liet}/\Z\Phi_{\liet}(x)$ is indeed torsion-free.
\end{proof}

To state the next proposition, recall that $V^{\reg}\subset V$ denotes the open subscheme of regular elements and that we have fixed a Kostant section $\kappa$ in \S\ref{subsection: The universal centraliser}.

\begin{proposition}\label{proposition: action map is surjective over B1}
The action map $G\times B^1 \rightarrow V^{\reg}|_{B^1},\, (g,b) \mapsto g\cdot \kappa_b$ is surjective.
\end{proposition}
\begin{proof} 
Part 2 of Proposition \ref{proposition: invariant theory (G,V) stable involution} implies that this map is surjective when restricted to $B^{\rs}$.
Therefore it suffices to prove that if $k/\Q$ is algebraically closed and $b\in D^1(k)$, then every two elements $x,y \in V^{\reg}_b(k)$ are $G(k)$-conjugate.
Let $x = x_s + x_n$ and $y = y_s + y_n$ be the Jordan decompositions of $x$ and $y$.
Since the semisimple parts $x_s$ and $y_s$ are $G(k)$-conjugate by Proposition \ref{proposition: basic invariant theory of vinberg representation full generality}, we may assume that $x_s = y_s$. 
The centraliser $L \coloneqq Z_H(x_s)$ is a reductive group with derived subgroup of type $A_1$ (Lemma \ref{lemma: points of D1 are those points with derived subgroup A1}); write $\mathfrak{l}$ for its Lie algebra.
The involution $\theta$ restricts to a stable involution $\theta|_L$ on $L$ \cite[Lemma 2.5]{Thorne-thesis}.
Since $x$ and $y$ are regular, $x_n, y_n$ are regular nilpotent elements of $\mathfrak{l}^{\theta= -1}$.
Therefore to prove the lemma, it suffices to prove that $L^{\theta} \cap G$ acts transitively on the regular nilpotents in $\mathfrak{l}^{\theta = -1}$.
(Note that $L^{\theta}\subset H^{\theta}$ but we don't have $L^{\theta} \subset G$ in general.)

We first claim that $L^{\theta}$ acts transitively on the regular nilpotents in $\mathfrak{l}^{\theta = -1}$.
To this end, let $Z(L)$ denote the centre of $L$ and consider the exact sequence
\begin{align}\label{equation: centre to pgl2 exact sequence}
0 \rightarrow Z(L) \rightarrow L \rightarrow \PGL_2 \rightarrow 0.
\end{align}
The involution $\theta$ preserves $Z(L)$ (acting via inversion by \cite[Lemma 2.7(3)]{Thorne-thesis}) and by Lemma \ref{lemma: stable involutions over alg closed field are conjugate} we may choose the isomorphism $L/Z(L) \simeq \PGL_2$ so that $\theta$ corresponds to the standard stable involution $\xi = \Ad\left(\text{diag}(1,-1) \right)$ of $\PGL_2$ from \S\ref{subsection: a stable Z/2Z-grading}.
An elementary computation in $\liesl_2$ (or Proposition \ref{proposition: H^theta acts simply transitively regular nilpotents}) shows that $\PGL_2^{\xi}$ acts transitively on the regular nilpotents in $\mathfrak{l}^{\xi = -1}$. 
To prove the claim, we only need to show that $L^{\theta} \rightarrow \PGL_2^{\xi}$ is surjective.
Since $Z(L)$ is connected (Lemma \ref{lemma: centre of centraliser semisimple element connect}), $Z(L)/(1-\theta)Z(L)$ is trivial and therefore taking $\theta$-invariants of \eqref{equation: centre to pgl2 exact sequence} shows that indeed $L^{\theta} \rightarrow \PGL_2^{\xi}$ is surjective, proving the claim.

To prove that $L^{\theta}\cap G$ acts transitively on the regular nilpotents in $\mathfrak{l}^{\theta = -1}$, it suffices to prove by the previous paragraph that $L^{\theta}\cap G$ surjects onto $\PGL_2^{\xi}$. 
We first claim that there exists a semisimple element $t\in V$ with centraliser $M = Z_H(t)$ such that $L\subset M$ and such that the derived subgroup of $M$ is of type $A_2$.
Indeed, take a Cartan subspace $\mathfrak{c} \subset V$ containing $x_s$; then $\Phi_{\mathfrak{c}}(x_s) = \{\pm \alpha\}$ for some root $\alpha$. 
Since $\Phi_{\mathfrak{c}}$ is not of type $A_1$,  there exists a root $\beta$ such that $\{\pm \alpha,\pm \beta, \pm(\alpha+\beta)\} \subset \Phi_{\mathfrak{c}}$. 
Taking $t$ to be an element of $\mathfrak{c}$ that vanishes exactly on those roots satisfies the requirements.

Again $\theta$ restricts to a stable involution on $M$ and the isomorphism $M/Z(M) \simeq \PGL_3$ can be chosen so that $\theta$ agrees with the standard stable involution $\psi$ of $\PGL_3$ from \S\ref{subsection: a stable Z/2Z-grading}. 
Again $Z(M)$ is connected by Lemma \ref{lemma: centre of centraliser semisimple element connect} and taking $\theta$-invariants of the analogue of the sequence \eqref{equation: centre to pgl2 exact sequence} for $M$ gives an exact sequence
\begin{align}
1\rightarrow Z(M)^{\theta} \rightarrow M^{\theta} \rightarrow \PGL_3^{\psi}\rightarrow 1.
\end{align}
A component group calculation (Corollary \ref{corollary: explicit description component group in terms of fundamental group}) shows that $\PGL_3^{\psi}$ is connected.
Therefore the identity component $(M^{\theta})^{\circ}$ maps surjectively onto $\PGL_3^{\psi}$. 
Since $(M^{\theta})^{\circ}\subset M^{\theta}\cap G$, this implies that $M^{\theta} = (M^{\theta} \cap G) \cdot Z(M)^{\theta}$. 
It follows that $L^{\theta} = (L^{\theta} \cap G) \cdot Z(L)^{\theta}$. 
Indeed, if $l\in L^{\theta}$ there exists an element $z\in Z(M)^{\theta}$ such that $lz\in M^{\theta}\cap G$.
Since $Z(M)^{\theta} \subset Z(L)^{\theta}$, we have $lz\in L^{\theta}\cap G$. 
We have established the equality $L^{\theta} = (L^{\theta} \cap G) \cdot Z(L)$, and it implies that $L^{\theta}\cap G$ surjects onto $\PGL_2^{\xi}$. 
\end{proof}

\begin{remark}
Proposition \ref{proposition: action map is surjective over B1} is false when $H$ is of type $A_1$.
Indeed, in that case $0 \in B^1 = B$, $G = \left\{\begin{pmatrix} 0 & x \\ y & 0 \end{pmatrix} \mid xy=0\right\}\simeq \G_m$ and there exist two $G$-orbits on $V_0^{\reg} = \left\{\begin{pmatrix} 0 & x \\ y & 0 \end{pmatrix} \mid xy=0\right\}$.
The problem is that, in the notation of the proof, $L^{\theta} \cap G=G$ does not surject onto the disconnected group $\PGL_2^{\xi}$.
\end{remark}



\subsection{Geometry over \texorpdfstring{$B^1$}{B1}}\label{subsection: the family of curves over B 1}

Recall from \S\ref{subsection: A family of curves} that we have introduced a family of projective curves $C\rightarrow B$ which is smooth exactly above $B^{\rs}$.

\begin{lemma}\label{lemma: curves over B1 have nodal singularity}
Let $k/\Q$ be a field and $b\in B(k)$. Then $b \in D^1(k)$ if and only if the curve $C_b$ has a unique nodal singularity.
\end{lemma}
\begin{proof}
A node is a simple singularity of type $A_1$. 
Therefore the lemma follows from Lemmas \ref{lemma: singularities of C_b using dynkin diagram} and \ref{lemma: points of D1 are those points with derived subgroup A1}.
\end{proof}

The fibres of the morphism $C\rightarrow B$ may be reducible. 
However, this does not happen over $B^1$:

\begin{lemma}\label{lemma: fibres over B1 of curves are integral}
The fibres of $C^1 \rightarrow B^1$ are geometrically integral. 
\end{lemma}
\begin{proof}
The geometric fibres of $C\rightarrow B$ are reduced, connected, and over $B^{\rs}$ these fibres are smooth. 
Therefore it suffices to prove that $C_b$ is irreducible if $k/\Q$ is algebraically closed and $b\in D^1(k)$. This can be verified by computation in each case, or by the following uniform argument.

Let $Z\subset D^1$ be the locus above which the fibres fail to be geometrically integral.
Then $Z$ is closed by \cite[Théorème 12.2.1(x)]{EGAIV-3}.
We claim that $Z$ is also open. 
To prove this, it suffices to prove that $Z$ is closed under generalisation \cite[Tag \href{https://stacks.math.columbia.edu/tag/0903}{0903}]{stacksproject}. 
By \cite[Tag \href{https://stacks.math.columbia.edu/tag/054F}{054F}]{stacksproject} this amounts to showing that for every complete discrete valuation ring $R$ and morphism $\Spec R\rightarrow D^1$, the generic point of $\Spec R$ lands in $Z$ if the closed point of $\Spec R$ does.
If $b \in D^1(k)$ then $C_b$ has a unique nodal singularity by Lemma \ref{lemma: curves over B1 have nodal singularity}, so either $C_b$ is irreducible with one node or a union of two irreducible components intersecting transversally.
Therefore if $b\in D^1(k)$, then $C_b$ is reducible (in other words, $b\in Z$) if and only if the generalised Jacobian $\Pic^0_{C_b/k}$ is an abelian variety \cite[\S9.2, Example 8]{BLR-NeronModels}. 
On the other hand, since the locus of $\Spec R$ where the relative Jacobian $\Pic^0_{C_R/R}\rightarrow \Spec R$ (which exists and is a semi-abelian scheme by \cite[\S9.3, Theorem 7]{BLR-NeronModels}) is an abelian variety is open (this follows from looking at torsion points), we conclude that the generic fibre of $\Pic^0_{C_R/R}$ is an abelian variety if its special fibre is. 
It follows that $Z$ is closed under generalisation, proving the claim.
Since $D^1$ is irreducible and $Z$ is open and closed, it follows that $Z$ is empty or equal to $D^1$; we will exclude the latter case.

To this end, it suffices to prove that $C_{\eta}$ is geometrically integral, where $\eta$ is the generic point of $D$. 
Assume by contradiction that this is not the case. 
As observed in the previous paragraph, this implies that $\Pic^0_{C_{\eta}/\eta}$ is an abelian variety.
Therefore the finite \'etale group scheme $J^{\rs}[2] \rightarrow B^{\rs}$ is unramified along $D$.
By Zariski--Nagata purity for finite \'etale covers \cite[Tag \href{https://stacks.math.columbia.edu/tag/0BJE}{0BJE}]{stacksproject}, this implies that $J^{\rs}[2]$ extends to a finite \'etale cover over $B$. Since $B$ is isomorphic to affine space over $k$, this cover must be trivial.
However, a monodromy calculation (Propositions \ref{proposition: monodromy of J[2]} and \ref{proposition: consequences of monodromy}) shows that $J^{\rs}[2]$ is nontrivial. (Recall that we have excluded that $H$ is of type $A_1$.) 
This is a contradiction, completing the proof of the lemma.
%
%
%
%
%
\end{proof}

Since $C^1 \rightarrow B^1$ has geometrically integral fibres by Lemma \ref{lemma: fibres over B1 of curves are integral}, the group scheme $\Pic^0_{C^1/B^1}$ is well-defined and we denote it by $J^1\rightarrow B^1$ \cite[\S9.3, Theorem 1]{BLR-NeronModels}. It is a semi-abelian scheme.
The $2$-torsion subgroup $J^1[2]\rightarrow B^1$ is a quasi-finite \'etale group scheme; we may therefore view it as a sheaf on the \'etale site of $B^1$.

\begin{lemma}\label{lemma: pushforward of Jrs[2] is J1[2]}
Let $j\colon B^{\rs} \hookrightarrow B^1$ be the open inclusion. Then $j_*J^{\rs}[2]  = J^1[2]$ as \'etale sheaves on $B^1$. 
\end{lemma}
\begin{proof}
Consider the natural morphism $\phi\colon J^1[2] \rightarrow j_*j^*J^1[2] = j_*J^{\rs}[2]$ obtained by adjunction. 
Since $J^1 \rightarrow B^1$ is separated, $J^1[2] \rightarrow B^1$ is separated as well so $\phi$ is injective. 
To prove that $\phi$ is an isomorphism, it suffices to check this at geometric points of $B^1$. 
Combining the last two sentences, it suffices to prove that $(J^1[2])_{\bar{b}}$ and $(j_*J^{\rs}[2])_{\bar{b}}$ have the same cardinality for all geometric points $\bar{b}$ of $B^1$, or even that the cardinality of the latter is bounded above by the cardinality of the first.
This is obvious if $\bar{b}$ lands in $B^{\rs}$, so assume that $\bar{b}$ lands in $D^1$. 

By Lemma \ref{lemma: fibres over B1 of curves are integral}, $C_{\bar{b}}$ is integral and has a unique singularity, which is a node. It follows that $J^1_{\bar{b}}$ has order $2^{2g-1}$, where $g$ is the arithmetic genus of $C_b$. 
On the other hand, the order of $(j_*J^{\rs}[2])_x$ for $x\in D^1$ can only go down under specialisation.
It therefore suffices to prove that if $\eta$ denotes the generic point of $D$, then $(j_*J^{\rs}[2])_{\eta}$ has order $2^{2g-1}$.
In fact, we claim that $\phi_{\eta}$ is an isomorphism.

To this end, let $K$ be the fraction field of the discrete valuation ring $\O_{B,\eta}$ and let $j_{\eta}$ be the inclusion $\Spec K \hookrightarrow \Spec \O_{B,\eta}$. 
Then the pullback of $j_*J^{\rs}[2]$ along $\Spec\O_{B,\eta} \rightarrow B$ equals $(j_{\eta})_*J_K^{\rs}[2]$, where $J_K^{\rs}$ denotes the pullback $J^{\rs}$ along the generic point of $B$.
The curve $C_{\O_{B,{\eta}}}$ is regular, since $\Spec \O_{B,\eta} \rightarrow B$ hence $C_{\O_{B,\eta}}\rightarrow C$ is formally smooth and the total space $C$ is smooth.
Therefore a result of Raynaud \cite[\S9.5, Theorem 1]{BLR-NeronModels} shows that the identity component of the Picard functor of $C_{\O_{B,\eta}}$, which equals $J_{\O_{B,\eta}}^1$ by definition, is isomorphic to the N\'eron model of $J_{K}$. 
By the N\'eron mapping property, this shows that  $J_{\eta}^1[2] = (j_*J^{\rs}[2])_{\eta}$. 
This completes the proof of the claim hence that of the lemma.
\end{proof}

\subsection{Summary of properties of \texorpdfstring{$B^1$}{B1}}\label{subsection: summary}

We summarise the properties of $D^1$ in the next theorem.

\begin{theorem}\label{theorem: summary properties mildly singular locus}
Let $k/\Q$ be an algebraically closed field and $b\in B(k)$. 
Then the following are equivalent:
\begin{enumerate}
\item $b\in D^1(k)$;
\item for every semisimple $v\in V_b(k)$, the derived subgroup of $Z_H(v)$ is of type $A_1$;
\item $C_b$ is irreducible and has a unique singular point, which is a node.
\end{enumerate}
\end{theorem}
\begin{proof}
Combine Lemmas \ref{lemma: points of D1 are those points with derived subgroup A1} and \ref{lemma: curves over B1 have nodal singularity}.
\end{proof}

\begin{theorem}\label{theorem: isomorphism centraliser 2-torsion over B1}
The isomorphism $J^{\rs}[2]\simeq Z^{\rs}$ from Proposition \ref{proposition: isomorphism 2-torsion and Kostant section centraliser thorne} uniquely extends to an isomorphism $J^1[2] \simeq Z^1$ of separated \'etale group schemes over $B^1$. 
\end{theorem}
\begin{proof}
Since $J^{\rs}[2]$ and $Z^{\rs}$ are dense in $J^1[2]$ and $Z^1$ respectively, uniqueness is clear. 
For the existence, denote the open immersion $B^{\rs} \hookrightarrow B^1$ by $j$.
Consider the composition $\psi\colon Z^1 \rightarrow j_* Z^{\rs} \xrightarrow{\sim} j_*J^{\rs}[2] \xrightarrow{\sim} J^1[2]$ of the adjunction morphism $Z^1\rightarrow j_*Z^{\rs}$, the pushforward of the isomorphism $Z^{\rs} \xrightarrow{\sim} J^{\rs}[2]$ along $j$ and the isomorphism $j_*J^{\rs}[2] \xrightarrow{\sim} J^1[2]$ of Lemma \ref{lemma: pushforward of Jrs[2] is J1[2]}.
Since $Z^1\rightarrow B^1$ is separated, $\psi$ is injective. 
It therefore suffices to prove that $Z^1_{\bar{b}}$ and $J^1_{\bar{b}}$ have the same cardinality for every geometric point $\bar{b}$ of $D^1$. 
If $g$ denotes the common arithmetic genus of the fibres of $C\rightarrow B$, then $Z^1_{\bar{b}}$ has cardinality $2^{2g-1}$ by Lemma \ref{lemma: centraliser has order 2^2g-1 over D1}.
On the other hand, $J^1_{\bar{b}}$ also has cardinality $2^{2g-1}$ by Lemmas \ref{lemma: curves over B1 have nodal singularity} and \ref{lemma: fibres over B1 of curves are integral}.
\end{proof}

Although it is not necessary for this paper, we expect that an isomorphism similar to the one of Theorem \ref{theorem: isomorphism centraliser 2-torsion over B1} holds over the whole of $B$, but we have not been able to prove it.

\section{The compactified Jacobian}\label{section: the compactified jacobian}

We keep the notations from Section \ref{section: recollections of Thorne's thesis}.
Recall from \S\ref{subsection: The universal centraliser} that the family of smooth projective curves $C^{\rs} \rightarrow B^{\rs}$ has Jacobian variety $J^{\rs} \rightarrow B^{\rs}$ which is itself a smooth and projective morphism.
The goal of this chapter is to extend the $B^{\rs}$-scheme $J^{\rs}$ to a proper $B$-scheme $\bar{J}$ with good geometric properties. 
We achieve this using the theory of the compactified Jacobian of Altman--Kleiman \cite{AltmanKleiman-CompactifyingThePicardScheme}, extended by Esteves \cite{Esteves-compactifyingrelativejacobian} to incorporate reducible curves.
Its construction is given in \S\ref{subsection: the definition} and its basic properties are summarised in Theorem \ref{theorem: summary compactified jacobian}.
Note that the occurrence of reducible fibres is the reason why the definition of $\bar{J}$ is more involved here than in our previous work \cite[\S4.3]{Laga-E6paper}, which treats the $E_6$ case and where only irreducible fibres are present.

The results of this chapter will be useful for the construction of orbits in \S\ref{section: constructing orbits} (specifically \S\ref{subsection: a universal torsor}) and the construction of integral representatives in \S\ref{subsection: proof of integral representatives general case}.

\subsection{Generalities on sheaves}\label{subsection: generalities on sheaves}

The following material is largely taken from \cite{Esteves-compactifyingrelativejacobian, MeloRapgnettaViviani-FinecompactifiedJacobians}.
Let $k$ be an algebraically closed field. By a \define{curve} we mean a reduced projective scheme of pure dimension $1$ over $k$.

\begin{definition}
A coherent sheaf $I$ on a connected curve $X$ is said to be
\begin{enumerate}
\item \define{rank-1} if $I_{\eta} \simeq \O_{X,\eta}$ as $\O_{X,\eta}$-modules for every generic point $\eta\in X$;
\item \define{torsion-free} if the associated points of $I$ are precisely the generic points of $X$;
\item \define{simple} if $\End_k(I) = k$.
\end{enumerate}
\end{definition}

We remark that the first two conditions imply the third if $X$ is irreducible and that every torsion-free rank-$1$ sheaf on a smooth curve is invertible.

To obtain a well-behaved moduli problem of torsion-free rank $1$ sheaves on a reducible curve, we use stability conditions introduced by Esteves \cite{Esteves-compactifyingrelativejacobian}.
A \define{subcurve} $Z$ of a curve $X$ is a closed $k$-subscheme that is reduced and of pure dimension $1$.
If $I$ is a torsion-free sheaf on $X$, its restriction to a subcurve $I|_Z$ is not necessarily torsion-free; it contains a biggest torsion subsheaf and the quotient of $I|_Z$ by this subsheaf is denoted by $I_Z$. 
The sheaf $I_Z$ is the unique torsion-free quotient of $I$ whose support is equal to $Z$. 

\begin{definition}
Let $E$ be a vector bundle on a connected curve $X$ of rank $r\geq 1$ and degree $-rd$.
Let $I$ be a torsion-free rank-$1$ sheaf on $X$ with Euler characteristic $\chi(I)=d$.
We say that $I$ is \define{$E$-semistable} if for every nonempty proper subcurve $Y\subsetneq X$ we have that 
\begin{align}\label{equation: definition E semistable}
\chi(I_Y) \geq -\frac{\deg(E|_Y)}{r}.
\end{align}
We say that $I$ is \define{$E$-stable} if for every nonempty proper subcurve the inequality \eqref{equation: definition E semistable} is strict.
\end{definition}

Given a vector bundle $E$ on $X$, we may define its \define{multislope} $\underline{q}^E = \{q^E_{C_i}\}$ as follows. It is a tuple of rational numbers, one for each irreducible component $C_i$ of $X$, defined by setting 
\begin{align*}
q^E_{C_i} \coloneqq -\frac{\deg(E|_{C_i})}{\rank E }.
\end{align*}
If $Y\subset X$ is a subcurve, write $q_Y^E \coloneqq \sum_{C_i \subset Y} q_{C_i}^E$, where the sum is taken over those irreducible components $C_i$ that are contained in $Y$.
If $E$ is of rank $r$ and degree $-rd$ then $q_X^E = d$. 
When the vector bundle $E$ is clear from the context we omit the superscript from the notation $\underline{q}^E$.

\begin{definition}\label{definition: general polarisation}
Let $X$ be a curve and $E$ a vector bundle on $X$ of rank $r$ and degree $-rd$ with multislope $\underline{q}$.
We say that $E$ is \define{general} if $q_Y \not\in \Z$ for any nonempty proper subcurve $Y\subsetneq X$.
\end{definition}

If $I$ is torsion-free rank-$1$ on $X$, then $I$ is $E$-semistable if and only if $\chi(I_Y) \geq q_Y$ for every nonempty proper subcurve $Y \subset X$, and $E$-stable if every such inequality is strict.
Therefore if $E$ is general, a torsion-free rank-$1$ sheaf on $X$ is $E$-semistable if and only if it is $E$-stable.

The next lemma shows that a family of simple torsion-free rank-$1$ sheaves has no unexpected endomorphisms.
For a quasi-coherent sheaf $\sh{F}$ on a scheme $X$, we write $\mathcal{E}nd(\sh{F})$ for the sheaf of $\O_X$-module endomorphisms of $\sh{F}$, which is again a quasi-coherent sheaf on $X$.

\begin{lemma}\label{lemma: relative family of simple sheaves is globally simple}
Let $p\colon \mathcal{X} \rightarrow T$ be a flat family of projective curves whose geometric fibres are reduced and connected.
Let $I$ be a locally finitely presented $\O_{\mathcal{X}}$-module, flat over $T$, whose geometric fibres above $T$ are simple torsion-free rank-$1$.
Then $p_*\mathcal{E}nd(I) = \O_T$.
\end{lemma}
\begin{proof}
Use \cite[Corollary (5.3)]{AltmanKleiman-CompactifyingThePicardScheme} and the assumption that $I$ is simple in each geometric fibre. 
\end{proof}

\subsection{The definition}\label{subsection: the definition}

Recall from \S\ref{subsection: A family of curves} that $C\rightarrow B$ is a flat projective morphism whose geometric fibres are reduced connected curves, and that this morphism has sections $\infty_1,\dots,\infty_m\colon B\rightarrow C$ landing in the smooth locus.

\begin{lemma}\label{lemma: irreducible components of Cb geom irreducible and contain section}
For every field $k/\Q$ and $b\in B(k)$, the irreducible components of $C_b$ are geometrically irreducible. 
Moreover, every such irreducible component contains $\infty_{i,b}$ in its smooth locus for some $i$. 
\end{lemma}
\begin{proof}
The first claim follows from the second one. 
For the second one, we may assume that $k$ is algebraically closed.
Consider the line bundle $\sh{L} = \O_C(\infty_1+ \cdots + \infty_m)$ on $C$ associated with the divisors $\infty_i$ of $C$.
For every $b\in B(k)$, $\sh{L}_b$ is ample if and only if every irreducible component of $C_b$ contains $\infty_{i,b}$ for some $i$. 
Moreover, the locus of elements $b\in B$ for which $\sh{L}_b$ is ample is open \cite[Corollaire (9.6.4)]{EGAIV-3}, $\G_m$-invariant (with respect to the $\G_m$-action on $C\rightarrow B$ introduced in \S\ref{subsection: A family of curves}) and contains the central point by Lemma \ref{lemma: bijection between marked points, irreducible components and certain nilpotents}.
These three facts imply that it must be the whole of $B$.
%
\end{proof}

In order to define a compactified Jacobian of $C\rightarrow B$, we first construct a vector bundle $E$ on $C$ using properties of the central fibre $C_0$. 
Recall from Lemma \ref{lemma: bijection between marked points, irreducible components and certain nilpotents} that each of the $m$ irreducible components of $C_0$ contains a unique marked point $\infty_{i,0}$. 
Let $\underline{q} = \{q_1,\dots,q_m\}$ be a tuple of rational numbers such that $\sum_{i =1}^{m} q_i = \chi(\O_{C_0}) =1-p_a(C_0)$ and $\sum_{ i \in I } q_i \not\in \Z$ for every nonempty proper subset $I\subset \{1,\dots,m\}$; it is easy to see that such a tuple exists.
Write $q_i = e_i/r$ for some $e_i \in \Z$ and $r\in \Z_{\geq 1}$.
By further multiplying $e_i$ and $r$, we may assume that $r\geq m$. 
Let $E$ be the following vector bundle on $C$:
\begin{align*}
E = \O_C(-e_1 \cdot \infty_1) \oplus \cdots \oplus \O_C(-e_m\cdot  \infty_m) \oplus \O_C^{\oplus r-m}.
\end{align*}
Since the image of $\infty_i \colon B\rightarrow C$ is a divisor of $C$, the line bundles $\O_C(-e_i \cdot \infty_i)$ are well-defined.
Note that the vector bundle $E|_{C_0}$ has multislope $\underline{q}$ by construction. 
For every geometric point $b$ of $B$, $E|_{C_b}$ is a vector bundle of rank $r$ and degree $-r (1-p_a(C_0))$ on the curve $C_b$.

\begin{lemma}\label{lemma: E is general vector bundle}
For every geometric point $b$ of $B$, the vector bundle $E|_{C_b}$ is general in the sense of Definition \ref{definition: general polarisation}.
\end{lemma}
\begin{proof}
    Follows from Lemma \ref{lemma: irreducible components of Cb geom irreducible and contain section} and the construction of $E$.
\end{proof}

We are now ready to define the compactified Jacobian associated with $E$.
We assume we have made a choice of $\underline{q}$ and $E$ as above.
Consider the functor
\begin{align}
\bar{\mathbb{J}}_E \colon \{B\text-{Schemes} \} \rightarrow \{\text{Sets}\}
\end{align}
sending a $B$-scheme $T$ to the set of equivalence classes of pairs $(I,\phi)$, where
\begin{itemize}
\item $I$ is a locally finitely presented $\O_{C_T}$-module, flat over $T$, with the property that for every geometric point $t$ of $T$, $I_t$ is simple torsion-free rank-$1$, $\chi(I_t) = \chi(\O_{C_0})$ and $I_t$ is $E_t$-stable;
\item $\phi$ is an isomorphism $\infty_{1,T}^*I \simeq \O_T$ of $\O_T$-modules.
\end{itemize}
We say two pairs $(I,\phi)$ and $(I',\phi')$ are equivalent if there is an isomorphism $I \simeq I'$ mapping $\phi$ to $\phi'$.
We have the following basic representability result \cite[Theorem B]{Esteves-compactifyingrelativejacobian}:
\begin{proposition}[Esteves]
    The functor $\bar{\mathbb{J}}_E$ is representable by a $B$-scheme $\bar{J}_E$.
\end{proposition}
\begin{proof}
Let $F$ be the functor from $B$-schemes to sets, sending a $B$-scheme $T$ to the set of equivalence classes of locally finitely presented $\O_{C_T}$-modules $I$, flat over $T$, with the property that for every geometric point $t$ of $T$, $I_t$ is simple torsion-free rank-$1$, $\chi(I_t) = \chi(\O_{C_0})$ and $I_t$ is $E_t$-stable. (In contrast to $\bar{\mathbb{J}}_E$, we omit the rigidification $\phi$.)
Here we say $I$ and $I'$ are equivalent if there exists an invertible sheaf $\sh{L}$ on $T$ such that $I' \simeq \sh{L}_{C_T} \otimes I$.
Let $F^{et}$ denote the \'etale sheafification of $F$.
By \cite[Proposition 34]{Esteves-compactifyingrelativejacobian}, the functor $F^{et}$ is representable by an open subspace of the algebraic space parametrising simple torsion-free rank-$1$ sheaves with no Euler characteristic or stability condition. 
By Lemma \ref{lemma: irreducible components of Cb geom irreducible and contain section} and \cite[Theorem B]{Esteves-compactifyingrelativejacobian}, the latter algebraic space is in fact a scheme, so $F^{et}$ is representable by a scheme as well.

On the other hand, the forgetful morphism $\bar{\mathbb{J}}_E\rightarrow F, (T,\phi) \mapsto T$ is an isomorphism of functors, since every $I\in F(T)$ is equivalent to another element $I' \in F(T)$ admitting a rigidification.
Since elements of $\bar{\mathbb{J}}_E$ have no nontrivial automorphisms by Lemma \ref{lemma: relative family of simple sheaves is globally simple}, \'etale descent of quasi-coherent sheaves implies that $\bar{\mathbb{J}}_E$ is an \'etale sheaf, so we have natural identifications $\bar{\mathbb{J}}_E = F = F^{et}$.
Since $F^{et}$ is representable by a scheme by the previous paragraph, the same is true for $\bar{\mathbb{J}}_E$.
\end{proof}

\begin{definition}\label{definition: compactified jacobian}
We call $\bar{J}_E$ a compactified Jacobian of $C\rightarrow B$ associated with $E$.
\end{definition}
If $C\rightarrow B$ has reducible fibres, different choices of $\underline{q}$ may give rise to different compactified Jacobians.
For our purposes, these differences will be harmless and for the remainder of this paper we fix a choice of $\underline{q}$ and $E$ as above and we simply write $\bar{J} = \bar{J}_E$.

\begin{lemma}\label{lemma: compactified jacobian in case of integral curves}
Let $k$ be a field and $b\in B(k)$ such that the curve $C_b$ is integral. Then $\bar{J}_b$ parametrises torsion-free rank-$1$ sheaves on $C_b$ with degree zero, i.e. Euler characteristic $1-p_a(C_0)$.
\end{lemma}
\begin{proof}
    If $C_b$ is integral, the $E_b$-stability condition and the simplicity of the sheaves are automatic.
\end{proof}

\subsection{Basic properties of \texorpdfstring{$\bar{J}$}{Jbar}}\label{subsection: basic properties}

\begin{lemma}
The morphism $\bar{J } \rightarrow B$ is projective\footnote{There are several nonequivalent definitions of a projective morphism but in this case they all agree, see \cite[Tag \href{https://stacks.math.columbia.edu/tag/0B45}{0B45}]{stacksproject}.}.
\end{lemma}
\begin{proof}
Since the vector bundle $E$ is chosen to be general (Lemma \ref{lemma: E is general vector bundle}), the notions of $E$-stable and $E$-semistable agree. 
Therefore, a theorem of Esteves \cite[Theorem C.1 and C.4]{Esteves-compactifyingrelativejacobian} shows that $\bar{J}$ is quasi-projective. 
Moreover, \cite[Theorem A.1]{Esteves-compactifyingrelativejacobian} shows that $\bar{J} \rightarrow B$ is universally closed. 
We conclude that $\bar{J}$ is projective over $B$.
\end{proof}

Recall from \S\ref{subsection: A family of curves} that we have defined a $\G_m$-action on $C$ such that $C\rightarrow B$ is $\G_m$-equivariant with respect to the square of the usual $\G_m$-action on $B$. 
By functoriality, this induces a $\G_m$-action on $\bar{J}$ too such that $\bar{J} \rightarrow B$ is $\G_m$-equivariant (again with respect to the square of the usual $\G_m$-action on $B$).
The following argument will be used in the next two lemmas: if $U\subset \bar{J}$ is an open $\G_m$-invariant subset containing the central fibre $\bar{J}_0$, then $U = \bar{J}$. 
Indeed, by the properness of $\bar{J}\rightarrow B$ the complement of $U$ in $\bar{J}$ projects to a closed $\G_m$-invariant subset of $B$ that does not contain the central point $0\in B$, so must be empty.

\begin{lemma}
The variety $\bar{J}$ is smooth.
\end{lemma}
\begin{proof}
The family $C\rightarrow B$ is a semi-universal deformation of the plane curve singularity $C_0$ (Proposition \ref{proposition: first properties of the family of curves}). 
Therefore \cite[Fact 4.2(ii)]{MeloRapgnettaViviani-FinecompactifiedJacobians} implies that $\bar{J}$ is smooth in a neighbourhood of $\bar{J}_0$. 
Since the smooth locus of $\bar{J}$ is open, $\G_m$-invariant and contains $\bar{J}_0$, it must be the whole of $\bar{J}$. 
\end{proof}
We emphasise that the fibres of $\bar{J}\rightarrow B$ might be singular above points that do no lie in $B^{\rs}$.
In fact, we have a precise description of the smooth locus:

\begin{lemma}\label{lemma: barJ is B-flat and smooth locus is line bundles}
The morphism $\bar{J} \rightarrow B$ is flat of relative dimension $p_a(C_0)$. The smooth locus of $\bar{J} \rightarrow B$ coincides with the locus of invertible sheaves.
\end{lemma}
\begin{proof}
    By  \cite[Theorem 5.5(ii)]{MeloRapgnettaViviani-FinecompactifiedJacobians}, the morphism $\bar{J} \rightarrow B$ is flat in a neighbourhood of $\bar{J}_0$. 
    Since the flat locus is open and $\G_m$-invariant, it follows that it must equal the whole of $\bar{J}$. 
    The claim about the smooth locus is \cite[Theorem 5.5(iii)]{MeloRapgnettaViviani-FinecompactifiedJacobians}.
 \end{proof}

\begin{lemma}\label{lemma: barJ is integral and has connected fibres and integral fibres when curve is integral}
The geometric fibres of $\bar{J} \rightarrow B$ are reduced and connected. Consequently, $\bar{J}$ is geometrically integral.
Moreover, if $k/\Q$ is an algebraically closed field and $b\in B(k)$ is such that $C_b$ is integral, then $\bar{J}_b$ is integral.
\end{lemma}
\begin{proof}
Since all the fibres of $C\rightarrow B$ have planar singularities, \cite[Theorem A(i)-(iii)]{MeloRapgnettaViviani-FinecompactifiedJacobians} shows that $\bar{J} \rightarrow B$ has geometrically reduced and connected fibres.
Since $\bar{J}$ is $B$-flat, this implies that $\bar{J}$ is geometrically connected. Since $\bar{J}$ is smooth, it follows that it is geometrically irreducible.
To establish the last claim, \cite[Corollary 5.14]{MeloRapgnettaViviani-FinecompactifiedJacobians} shows that the number of irreducible components of $\bar{J}_b$ can be calculated in terms of the intersections between the irreducible components of $C_b$.
This number is always $1$ when $C_b$ is irreducible, as can be seen from \cite[Definition 5.12]{MeloRapgnettaViviani-FinecompactifiedJacobians}.
\end{proof}

For future reference, we summarise the above properties in the following theorem.
Write $\bar{J}^1$ for the restriction of $\bar{J}$ to $B^1$.
Recall from \S\ref{subsection: the family of curves over B 1} that $J^1\rightarrow B^1$ is the relative generalised Jacobian of the family of integral curves $C^1\rightarrow B^1$. 
Note that by Lemma \ref{lemma: compactified jacobian in case of integral curves} and the definition of $J^1$ we have an open embedding $J^1\rightarrow \bar{J}^1$.

\begin{theorem}\label{theorem: summary compactified jacobian}
Let $\bar{J}\rightarrow B$ be a compactified Jacobian associated with some choice of $E$ as in Definition \ref{definition: compactified jacobian}.
Then the morphism $\bar{J} \rightarrow B$ is flat, projective and restricts to $J^{\rs}$ over $B^{\rs}$. Its geometric fibres are reduced and connected. The scheme $\bar{J}$ is geometrically integral and smooth over $\Q$. 
The smooth locus of $\bar{J}^1\rightarrow B^1$ is isomorphic to $J^1\rightarrow B^1$.
The complement of $J^1$ in $\bar{J}$ has codimension $\geq 2$.
\end{theorem}
\begin{proof}
Only the last two sentences remain to be established. 
The claim about the smooth locus of $\bar{J}^1$ follows from Lemmas \ref{lemma: compactified jacobian in case of integral curves} and \ref{lemma: barJ is B-flat and smooth locus is line bundles} and the definition of $J^1$.
For the claim about the codimension, let $Z$ be the complement of $J^1$ in $\bar{J}$. 
Then $Z$ is supported above the discriminant locus $D$ of $B$.
Moreover the fibres of the map $\bar{J}|_{D^1}\rightarrow D^1$ are geometrically integral by Lemmas \ref{lemma: fibres over B1 of curves are integral} and \ref{lemma: barJ is integral and has connected fibres and integral fibres when curve is integral}, so the fibres of the map $Z|_{D^1} \rightarrow D^1$ have dimension strictly less than those of $\bar{J}|_{D^1} \rightarrow D^1$. 
Combining the last two sentences proves the claim.
\end{proof}

\subsection{The Białynicki-Birula decomposition of \texorpdfstring{$\bar{J}$}{Jbar}}\label{subsection: the BB decomposition}

We recall the Białynicki-Birula decomposition \cite{BialynickBirula-sometheoremsactionsalgebraicgroups} from geometric representation theory.
If $k$ is a field and $X$ is a scheme of finite type of $k$, we define a \define{decomposition} of $X$ to be a collection of locally closed subschemes $X_1,\dots,X_n$ of $X$ such that the underlying topological space of $X$ is a disjoint union of the underlying topological spaces of the $X_i$.
If in addition $X$ is separated over $k$ and endowed with a $\G_m$-action and if $x\in X$ we say that $\lim_{\lambda\rightarrow 0}\lambda\cdot x$ exists if the action map $\G_m\rightarrow X, \lambda\mapsto \lambda\cdot x$ extends (necessarily uniquely) to a morphism $\mathbb{A}^1_k\rightarrow X$.

\begin{proposition}\label{proposition: BB decomposition in general}
Suppose that $X$ is a smooth and separated scheme of finite type over a field $k$, endowed with a $\G_m$-action.
Then the closed subscheme of fixed points $X^{\G_m}$ is smooth; let $F_1,\dots,F_n$ denote its connected components.
Suppose in addition that $\lim_{\lambda \rightarrow 0} \lambda\cdot x$ exists for every $x\in X$.
Then there exists a decomposition of $X$ into locally closed subschemes $X_i$ and morphisms $X_i \rightarrow F_i$ which are affine space fibrations in the Zariski topology.
\end{proposition}
\begin{proof}
See \cite[Theorem 1.5]{JelisiejewSienkiewicz-BBdecompositionreductivegroups} for a modern proof, which treats the generality in which we have stated it.
We may informally describe $X_i$ as those points $x\in X$ whose limit $\lim_{\lambda \rightarrow 0} \lambda\cdot x$ lies in $F_i$, and the map $X_i \rightarrow F_i$ as taking the limit $x\mapsto \lim_{\lambda\rightarrow 0} \lambda\cdot x$.
\end{proof}

\begin{corollary}\label{corollary: BB decomposition implies open iso to affine space}
In the setting of Proposition \ref{proposition: BB decomposition in general}, assume furthermore that $X$ is geometrically integral and $X^{\G_m}$ is finite. Then there exists an open subset of $X$ isomorphic to affine space $\A_k^{\dim X}$.
\end{corollary}
\begin{proof}
Let $F_1,\dots,F_n$ denote the connected components of $X^{\G_m}$; since $X^{\G_m}$ is smooth and finite each $F_i$ is the spectrum of a separable field extension $k_i$ of finite degree over $k$. 
Let $X_1,\dots,X_n$ be the decomposition of $X$ of Proposition \ref{proposition: BB decomposition in general}. 
Then each $X_i$ is isomorphic to $\A^{n_i}_{k_i}$ for some integer $n_i\geq 0$. 
There exists an $X_i$, say $X_1$, which is of maximal dimension $\dim X$.
Since $X_1$ is locally closed, it is an open subset of its closure $\bar{X}_1 = X$, so $X_1$ is an open subset of $X$.
Since $X$ is geometrically irreducible, the same is true for $X_1$. 
This implies that $X_1 \times_k k_1$ is irreducible, so $k_1 = k$. 
Therefore $X_1$ is isomorphic to $\A^{\dim X}_k$. 
\end{proof}

\begin{remark}
The proof shows that under the assumptions of Corollary \ref{corollary: BB decomposition implies open iso to affine space}, $X$ is even decomposed into affine cells. 
\end{remark}

We will apply Corollary \ref{corollary: BB decomposition implies open iso to affine space} to the compactified Jacobian $\bar{J}\rightarrow B$ constructed in \S\ref{subsection: the definition}.
Recall from \S\ref{subsection: basic properties} that $\bar{J}$ inherits a $\G_m$-action from $C$. 
We denote the central fibre of $\bar{J}$ by $\bar{J}_0$.

\begin{lemma}\label{lemma: Gm-fixed points compactified jacobian are finite}
The set of $\G_m$-fixed points $\bar{J}_0^{\G_m}$ is finite.
\end{lemma}
\begin{proof}
This follows from calculations of Beauville \cite[\S4.1]{Beauville-rationalcurvesK3} if $C_0$ is integral.
It seems likely that one can extend his analysis to reducible curves, but we will proceed differently.
We will assume that all schemes are base changed to a fixed algebraic closure $k$ of $\Q$. 
Let $J(C_0)$ be the generalised Jacobian of $C_0$ parametrising line bundles having multidegree zero, i.e. degree zero on each irreducible component of $C_0$.
Then $J(C_0)$ is an algebraic group acting on $\bar{J}_0$, compatibly with the $\G_m$-actions on $J(C_0)$ and $\bar{J}_0$.

First we claim that the closure of every $\G_m$-orbit of a point in $J(C_0)$ contains the identity.
Indeed, every point of $[L]\in J(C_0)(k)$ is represented by a Cartier divisor $D_1-\deg(D_1)\infty_{1,0} + \dots+ D_m - \deg(D_m)\infty_{m,0}$, where $D_i$ is a Cartier divisor supported on the smooth affine part of the irreducible component $C_{0,i}$ of $C_0$ containing $\infty_{i,0}$. 
Since every smooth point $P$ of $C_{0,i}$ satisfies $\lim_{\lambda\rightarrow \infty} \lambda\cdot P = \infty_{i,0}$ (as can be seen from the definition of the $\G_m$-action on $C_0$ in \S\ref{subsection: A family of curves} and Table \ref{table: introduction different cases}), we see that $\lambda\cdot [L] \rightarrow 0$ as $\lambda\rightarrow \infty$, proving the claim.

Secondly, we claim that the action of $J(C_0)$ on $\bar{J}_0$ has finitely many orbits. 
Indeed, let $p\in C_0$ be the unique singular point.
Since $p$ is an ADE-singularity, there are only finitely many isomorphism classes of torsion-free rank-$1$ modules over the completed local ring $\widehat{\O}_{C_0,p}$ (see \cite{Knorrer-CohenMacaulaymodulesonhypersurfacesingularities}, in fact this property can be used to characterise ADE-singularities amongst Gorenstein singularities).
It therefore suffices to prove that if $[\sh{F}],[\sh{G}] \in \bar{J}_0(k)$ are two sheaves whose completed stalks at $p$ are isomorphic, then $\sh{F} \simeq \sh{G} \otimes \sh{L}$, where $\sh{L}$ is a line bundle on $C_0$ whose multidegree can only take finitely many values (independently of $\sh{F}$ and $\sh{G}$).
To prove this, consider the Hom-sheaf $\sh{H} = \mathcal{H}om(\sh{F},\sh{G})$ and the endomorphism sheaf $\sh{E} = \mathcal{E}nd(\sh{F})$.
Since $\sh{F}|_{C_0 \setminus \{p\}}$ is a line bundle, $\sh{E}$ is a coherent commutative $\O_{C_0}$-algebra which is generically isomorphic to $\O_{C_0}$, and $\sh{H}$ is a coherent $\sh{E}$-module.
Since the formation of $\sh{H}$ and $\sh{E}$ commutes with flat base change \cite[Tag \href{https://stacks.math.columbia.edu/tag/0C6I}{0C6I}]{stacksproject}, the completed stalk $\sh{H}_p \otimes \widehat{\O}_{C_0,p}$ is free of rank $1$ over $\sh{E}_p \otimes \widehat{\O}_{C_0,p}$.
It follows that $\sh{H}_p$ is free of rank $1$ over $\sh{E}_p$, so the stalks $\sh{F}_p, \sh{G}_p$ are isomorphic $\O_{C_0,p}$-modules. (An isomorphism is given by choosing an $\sh{E}_p$-generator of $\sh{H}_p$.)
By spreading out such an isomorphism, we may find an open subset $U \subset C_0$ containing $p$ and an isomorphism $\phi_U \colon \sh{F}|_U \xrightarrow{\sim} \sh{G}|_U$. 
The restrictions of $\sh{F}, \sh{G}$ to $C_0 \setminus \{p\}$ are line bundles.
Since $C_0$ is connected and $p$ is its unique singular point, the complement $C_0 \setminus U$ is a union of finitely many points.
We may therefore find an open subset $V \subset C_0 \setminus \{p \}$ containing those points and an isomorphism $\phi_V \colon \sh{F}|_V \xrightarrow{\sim} \sh{G}|_V$.
The transition map $(\phi_V)|_{U\cap V}^{-1} \circ (\phi_U)|_{U\cap V}\colon \sh{F}|_{U\cap V} \xrightarrow{\sim } \sh{F}|_{U\cap V}$ defines an element $f\in \HH^0(U\cap V,\O_{C_0}^{\times})$.
Let $\sh{L}$ be the line bundle on $C_0$ obtained by glueing $\O_U$ and $\O_V$ along the automorphism $f$.
One can then explicitly check that the maps 
\begin{align*}
\sh{F}_U \otimes \O_U \rightarrow \sh{G}_U:  s\otimes 1\mapsto \phi_U(s), \\
\sh{F}_V \otimes \O_V \rightarrow \sh{G}_V: s\otimes 1\mapsto \phi_V(s),
\end{align*}
glue to an isomorphism $\sh{F} \otimes \sh{L} \simeq \sh{G}$. 
The multidegree of $\sh{L}$ can only take on finitely many values (when we vary $\sh{F}$ and $\sh{G}$ in $\bar{J}_0$) because of the $E$-stability condition imposed on sheaves in $\bar{J}_0$.
This completes the proof of the claim.
(We thank Jesse Leo Kass for his help with the proof of this claim.)

We now use the last two paragraphs to show that $\bar{J}_0^{\G_m}$ is finite. 
Indeed, by the second claim, it suffices to prove that every $J(C_0)$-orbit contains at most one $\G_m$-fixed point. 
If $x\in \bar{J}_0^{\G_m}$ and $g\in J(C_0)$ are such that $g\cdot x \in \bar{J}_0^{\G_m}$, then $g^{-1} (\lambda\cdot g )$ lies in the stabiliser of $x$ in $J(C_0)$ for all $\lambda\in \G_m$. 
Since this stabiliser is closed, the first claim implies that it contains $\lim_{\lambda\rightarrow \infty} g^{-1} \lambda\cdot g = g^{-1}$. 
Therefore $g$ lies in the stabiliser of $x$, that is $g\cdot x = x$. 
We conclude that $x$ is the only $\G_m$-fixed point in the $J(C_0)$-orbit of $x$. 
\end{proof}

\begin{theorem}\label{theorem: compactified jacobian has a dense open isomorphic to affine space}
The variety $\bar{J}$ has a dense open subset isomorphic to affine space $\A_{\Q}^d$ for some $d\geq1$.
\end{theorem}
\begin{proof}
The compactified Jacobian is a smooth, geometrically integral and quasi-projective scheme over $\Q$ (Theorem \ref{theorem: summary compactified jacobian}).
Since $\bar{J} \rightarrow B$ is proper and $\lim_{\lambda\rightarrow 0} \lambda\cdot b$ exists for every $b\in B$, $\lim_{\lambda \rightarrow 0} \lambda\cdot x$ exists for every $x\in \bar{J}$.
We wish to apply Corollary \ref{corollary: BB decomposition implies open iso to affine space}, so it suffices to prove that the fixed point locus $\bar{J}^{\G_m}$ is finite. 
Since $B^{\G_m}$ consists of the central point $0$, it suffices to prove that $\bar{J}_0^{\G_m}$ is finite, which is exactly Lemma \ref{lemma: Gm-fixed points compactified jacobian are finite}.
\end{proof}

\begin{remark}\label{remark: BB decomposition mumford rep}
The Białynicki-Birula decomposition gives a canonical decomposition of $\bar{J}$ into locally closed subschemes isomorphic to affine space. 
This can be made very explicit when $H$ is of type $A_{2}$.
In that case $C\rightarrow B$ is the universal family of elliptic curves with short Weierstrass form $y^2 = x^3 +p_2x+p_3$ (see Table \ref{table: introduction different cases}), and it turns out \cite[Proposition 7.3]{MeloRapgnettaViviani-FinecompactifiedJacobians} that there is an isomorphism of $B$-schemes $\bar{J} \simeq C$ respecting the $\G_m$-actions. 
The BB decomposition of $C$ has the form
\begin{align}
    C = C^{\circ} \sqcup \infty(B)
\end{align}
where $C^{\circ}$ is the affine part and $\infty(B)\simeq B$ is the image of the marked section at infinity $\infty\colon B\rightarrow C$.
The fact that these strata are isomorphic to affine space can be seen directly: $\infty(B) \simeq B=\A^2_{\Q}$ and $C^{\circ}$ is isomorphic to the closed subscheme $\{(x,y,p_2,p_3) \mid y^2= x^3+p_2x+p_6\}\subset \A^4_{\Q}$, which via projection to the first three coordinates is isomorphic to $\A^3_{\Q}$.

In case $A_{2g}$ it is true (but we do not show) that the BB decomposition stratifies the smooth locus of $\bar{J}\rightarrow B$ according to the Mumford representation of a divisor on a hyperelliptic curve \cite[IIIa]{Mumford-tatalecturesiii}.
It would be interesting to obtain a similarly concrete interpretation of this decomposition for other families of curves studied here.
\end{remark}

\section{Constructing orbits}\label{section: constructing orbits}

We keep the notations from Section \ref{section: recollections of Thorne's thesis}.
The goal of this chapter is to construct for every $b\in B^{\rs}(\Q)$ and every element of $\Sel_2 J_b$ a $G(\Q)$-orbit of $V_b(\Q)$, see Corollary \ref{corollary: inject 2-Selmer orbits}.
The technical input is the Zariski triviality of a certain universal torsor on $J^{\rs}$ in \S\ref{subsection: a universal torsor}, see Theorem \ref{theorem: universal torsor is Zariski trivial}.
This will be achieved using generalities concerning torsors on open subsets of affine spaces developed in \S\ref{subsection: torsors on open subsets affine space}.

\subsection{Torsors on open subsets of affine space}\label{subsection: torsors on open subsets affine space}

The purpose of this subsection is to prove the following theorem, which will be useful in the proof of Theorem \ref{theorem: universal torsor is Zariski trivial}.

\begin{theorem}\label{theorem: torsors open subsets affine space codim 2}
Let $k$ be a field of characteristic zero and $X$ an open subset of $\A_k^n$ whose complement has codimension $\geq 2$. 
Let $G$ be a reductive group over $k$ and let $T\rightarrow X$ be a $G$-torsor.
Suppose that $X$ contains a $k$-rational point over which $T$ is trivial. Then $T$ is Zariski locally trivial.
\end{theorem}

\begin{example}\label{example: illustration of theorem on torsors affine spaces PGL2}
We illustrate Theorem \ref{theorem: torsors open subsets affine space codim 2} in the concrete case $G = \PGL_2$.
If $k$ and $X$ are as in the theorem, then a $\PGL_2$-torsor can alternatively be viewed as a Severi--Brauer curve $\mathcal{C} \rightarrow X$.
In other words, $\mathcal{C} \rightarrow X$ is a smooth projective family of genus zero curves (i.e. conics).
Suppose that $X$ contains a point $x\in X(k)$ such that the conic $\mathcal{C} _x$ has a $k$-rational point.
Then Theorem \ref{theorem: torsors open subsets affine space codim 2} says that $\mathcal{C} $ is a projective bundle over $X$.
This implies that all the other fibres of $\mathcal{C} \rightarrow X$ also contain a $k$-rational point.
\end{example}

Theorem \ref{theorem: torsors open subsets affine space codim 2} might be unsurprising to experts, although we have not been able to locate it explicitly in the literature.
It follows from a slight variant of the formalism developed in \cite[\S 1]{ColliotThelene-Formesquadratiquesdeuxcomplements} using pointed sets rather than abelian groups.
We take a different route and give a short proof which was suggested to us by Colliot-Th\'el\`ene. 
We thank him for letting us to include the argument here.

The crucial input is the following result of Raghunathan--Ramanathan, see \cite[Theorem 1.1]{Raghunathan-Ramanathan} and the remark immediately thereafter.
\begin{proposition}\label{proposition: raghunatan-ramanathan}
Let $k$ be a perfect field and $G$ a reductive group over $k$. 
Then every $G$-torsor on $\A^1_k$ is isomorphic to the pullback of a $G$-torsor on $\Spec k$ along the map $\A^1_k\rightarrow \Spec k$.
\end{proposition}

\begin{proof}[Proof of Theorem \ref{theorem: torsors open subsets affine space codim 2}]
Let $p\colon \A^n_k \rightarrow \A_k^{n-1}$ be the projection onto the first $n-1$ coordinates and let $p_X\colon X\rightarrow \A^{n-1}_k$ be its restriction to $X$.
If $K$ denotes the function field of $\A^{n-1}_k$, then the assumptions on $X$ imply that the generic fibre of $p_X$ is isomorphic to $\A^1_K$.
By Proposition \ref{proposition: raghunatan-ramanathan}, the restriction of $T$ to this generic fibre is induced from a torsor $T_0$ on $\Spec K$ along the map $\A^1_K\rightarrow \Spec K$. 

Since $k$ is infinite, we may choose a section $s\colon \A^{n-1}_k\rightarrow \A^n_k$ of $p$ such that its image does not contain any generic point of the complement of $X$ in $\A^n_k$.
Let $X_1 = s^{-1}(X\cap s(\A^{n-1}_k))$ and let $T_1$ be the pullback of $T$ along $s|_{X_1}\colon X_1 \rightarrow X$.
By our choice of $s$, $X_1$ is an open subset of $\A^{n-1}_k$ whose complement has codimension $\geq 2$.
The generic fibre of $T_1$ is a torsor on $\Spec K$ isomorphic to $T_0$.

Replacing $T$ by $T_1$ and iterating this process, we obtain a torsor $T_n$ on $\Spec k$ whose pullback along $X\rightarrow \Spec k$ is generically isomorphic to $T$. 
By known cases of the Grothendieck--Serre conjecture (Corollary \ref{corollary: known cases Grothendieck Serre conjecture spreading out version}), $T$ and $T_n\times_k X$ are Zariski locally isomorphic.
Since $X$ has a $k$-point above which $T$ is trivial, the same is true for $T_n\times_k X$ hence $T_n$ is trivial itself. 
We conclude that $T$ is Zariski locally trivial, as desired.
\end{proof}

\subsection{A universal torsor}\label{subsection: a universal torsor}

Recall from \S\ref{subsection: The universal centraliser} that $J^{\rs} \rightarrow B^{\rs}$ denotes the relative Jacobian of the family of smooth curves $C^{\rs} \rightarrow B^{\rs}$,
that $Z\rightarrow B$ denotes the universal stabiliser of the Kostant section $\kappa$, and that there is an isomorphism of finite \'etale group schemes $J^{\rs}[2] \simeq Z^{\rs}$ over $B^{\rs}$.

Since $J^{\rs}\rightarrow B^{\rs}$ is an abelian scheme, the multiplication-by-$2$ map $J^{\rs} \xrightarrow{\times 2} J^{\rs}$ is a $J^{\rs}[2]$-torsor. 
Pushing out this torsor along the maps $J^{\rs}[2] \xrightarrow{\sim} Z^{\rs} \hookrightarrow G$ defines a $G$-torsor $T^{\rs} \rightarrow J^{\rs}$. (This procedure is also called `changing the structure group'.)
The following theorem is one of the main technical results of this paper, and is the essential input for constructing orbits associated with elements of $J_b(\Q)$ (Theorem \ref{theorem: inject 2-descent orbits}).

\begin{theorem}\label{theorem: universal torsor is Zariski trivial}
The torsor $T^{\rs}$ is Zariski locally trivial. That is, for every $x\in J^{\rs}$ there exists an open subset $U \subset J^{\rs}$ containing $x$ such that $T^{\rs}|_U$ is trivial. 
\end{theorem}

To briefly explain why this is relevant for constructing orbits, note that if $b\in B^{\rs}(\Q)$ and $P\in J_b(\Q)$, the image of $P$ under the composition 
$J_b(\Q)/2J_b(\Q) \rightarrow \HH^1(\Q,J_b[2]) \xrightarrow{\sim} \HH^1(\Q,Z_b) \rightarrow \HH^1(\Q,G)$
coincides with the isomorphism class of the pullback of $T^{\rs}$ along $P \colon \Spec \Q \rightarrow J^{\rs}$. 
Theorem \ref{theorem: universal torsor is Zariski trivial} implies that this pullback defines the trivial class in $\HH^1(\Q,G)$, which implies that it corresponds to a $G(\Q)$-orbit of $V_b(\Q)$; see Theorem \ref{theorem: inject 2-descent orbits} for full details.

\begin{proof}[Proof of Theorem \ref{theorem: universal torsor is Zariski trivial}]
Recall from \S\ref{section: the mildly singular locus} that $B^1\subset B$ is an open subset containing $B^{\rs}$ and that the family of curves $J^1\rightarrow B^1$ is the relative (generalised) Jacobian of the family of curves $C^1\rightarrow B^1$. 
By Theorem \ref{theorem: isomorphism centraliser 2-torsion over B1} the isomorphism $J^{\rs}[2] \simeq Z^{\rs}$ extends to an isomorphism $J^1[2] \simeq Z^1$ of quasi-finite \'etale group schemes over $B^1$.
The multiplication-by-two map $J^1\xrightarrow{\times 2} J^1$ is a $J^1[2]$-torsor, and pushing out this torsor along the composition $J^1[2] \xrightarrow{\sim} Z^1 \rightarrow G$ defines a $G$-torsor $T^1\rightarrow J^1$. By construction, the restriction of $T^1$ to $J^{\rs}$ is isomorphic to $T^{\rs}$.

To prove the theorem, it suffices to prove that $T^1$ is Zariski locally trivial. 
Using known cases of the Grothendieck--Serre conjecture (Corollary \ref{corollary: known cases Grothendieck Serre conjecture spreading out version}) it even suffices to prove that $T^1$ is Zariski locally trivial when restricted to a nonempty open subset of $J^1$. 

Recall from \S\ref{section: the compactified jacobian} that we have constructed a scheme $\bar{J} \rightarrow B$ containing $J^1$ as an open subscheme. 
By Theorem \ref{theorem: summary compactified jacobian}, the complement of $J^1$ in $\bar{J}$ has codimension $\geq 2$; by Theorem \ref{theorem: compactified jacobian has a dense open isomorphic to affine space}, $\bar{J}$ contains an open dense subscheme $U$ isomorphic to affine $\Q$-space.
This implies that the complement of $U^1 \coloneqq U\cap J^1$ in $U$ has codimension $\geq 2$. 

We claim that $T^1|_{U^1}$ is Zariski locally trivial. By Theorem \ref{theorem: torsors open subsets affine space codim 2}, it suffices to prove that $T^1_x$ is trivial for some $x\in U^1(\Q)$. 
In fact, we will show the stronger statement that $\{ x\in J^{1}(\Q) \mid T^1_x \text{ is trivial} \}$ is Zariski dense in $J^{1}$.
Indeed, since $J^{1}$ is a rational variety (it contains $U^1$ as a dense open subscheme), the set $J^1(\Q)$ is dense in $J^1$.
Since the multiplication-by-two map $J^1\xrightarrow{\times 2} J^1$ is dominant, the subset $2J^1(\Q)\subset J^1(\Q)$ is still dense in $J^1$.
By construction of $T^1$, the pullback of $T^1$ along a point $x\in 2J^1(\Q)$ is trivial. This completes the proof of the claim, hence the proof of the theorem.
\end{proof}

\subsection{Constructing orbits for 2-descent elements}\label{subsection: constructing orbits for 2-descent elements}

We start by applying a well known lemma from arithmetic invariant theory recalled in \S\ref{subsection: arithmetic invariant theory} to give a cohomological description of the $G$-orbits of $V$.
For every $\Q$-algebra $R$ and $b\in B^{\rs}(R)$ we write $J_b\coloneqq J^{\rs}_b$ for the Jacobian of $C_b$.

\begin{corollary}\label{corollary: AIT for orbits of V}
Let $R$ be a $\Q$-algebra and $b\in B^{\rs}(R)$. 
Then the association $v\mapsto \{g\in G \mid g\cdot v = \kappa_b\}$ induces an injection $$\gamma_b \colon G(R) \backslash V_b(R)\hookrightarrow \HH^1(R,J_b[2]).$$ 
Its image coincides with the pointed kernel of the map $\HH^1(R,J_b[2]) \xrightarrow{\sim} \HH^1(R,Z_b) \rightarrow \HH^1(R,G)$. 
\end{corollary}
\begin{proof}
	We apply Lemma \ref{lemma: AIT full generality} to the action of $G_{B^{\rs}}$ on $V^{\rs}$. 
	Indeed, the action map $G \times B^{\rs} \rightarrow V^{\rs}, (g,b) \mapsto g\cdot \kappa_b$ is \'etale (Proposition \ref{proposition: Kostant section properties}) and it is surjective by Proposition \ref{proposition: invariant theory (G,V) stable involution}. 
	Pulling back along $b\colon \Spec R\rightarrow B^{\rs}$ and using the isomorphism $J_b[2] \simeq Z_b$ from Proposition \ref{proposition: isomorphism 2-torsion and Kostant section centraliser thorne} gives the desired bijection. 
\end{proof}

We now piece all the ingredients obtained so far together to deduce our first main theorem.

\begin{theorem}\label{theorem: inject 2-descent orbits}
Let $R$ be a local $\Q$-algebra (for example, a field of characteristic zero) and $b\in B^{\rs}(R)$.
Then the image of the $2$-descent map $J_b(R)/2J_b(R)\rightarrow \HH^1(R,J_b[2])$ lies in the image of $\gamma_b$ of Corollary \ref{corollary: AIT for orbits of V}. 
Consequently, there is a canonical injection $$\eta_b\colon J_b(R)/2J_b(R) \hookrightarrow G(R)\backslash V_b(R)$$ compatible with base change. 
\end{theorem}
\begin{proof}
By Corollary \ref{corollary: AIT for orbits of V}, it suffices to prove that the composition
$J_b(R)/2J_b(R) \rightarrow \HH^1(R,J_b[2]) \simeq \HH^1(R,Z_b) \rightarrow \HH^1(R,G)$ is trivial.
Recall that in \S\ref{subsection: a universal torsor} we have constructed a $G$-torsor $T^{\rs}\rightarrow J^{\rs}$ such that its pullback along a point $P\colon \Spec R \rightarrow J^{\rs}$ defines a $G$-torsor $T^{\rs}_P \rightarrow \Spec R$ whose isomorphism class equals the image of $P$ under the above composite map. 
Since $T^{\rs}$ is Zariski locally trivial by Theorem \ref{theorem: universal torsor is Zariski trivial}, $T^{\rs}_P$ is Zariski locally trivial. Since $R$ is a local ring, it follows that $T^{\rs}_P$ is trivial.
This completes the proof.
\end{proof}

\begin{remark}
In the proof of Theorem \ref{theorem: universal torsor is Zariski trivial} we have shown the stronger statement that the torsor $T^1\rightarrow J^1$ (a natural extension of $T^{\rs}$ to $J^1$) is Zariski locally trivial.
A straightforward adaption of the proof of Theorem \ref{theorem: inject 2-descent orbits} then shows that if $R$ is a local ring and $b\in B^1(R)$ (instead of $b\in B^{\rs}(R)$), then there exists an injection $J_b^1(R)/2J_b^1(R) \hookrightarrow G(R) \backslash V_b^{\reg}(R)$. 
We do not know if this observation is useful.
\end{remark}

\subsection{Constructing orbits for 2-Selmer elements}\label{subsection: constructing orbits for Selmer elements}

The next proposition might be well known to experts --- see for example \cite[Remark after Theorem 6.22]{PlatonovRapinchuk-Alggroupsandnumbertheory} --- but we believe it deserves to be stated explicitly.
We slightly deviate from our standing notation and allow $G$ to be an arbitrary split semisimple group in this proposition.

\begin{proposition}\label{proposition: weak hasse principle for H1(G)}
Let $G$ be a split semisimple group over a number field $k$. 
Then the kernel of $\HH^1(k,G) \rightarrow \prod_v \HH^1(k_v,G)$ (where $v$ runs over all places) is trivial. 
\end{proposition}
We emphasise that $\HH^1(k,G) \rightarrow \prod_v \HH^1(k_v,G)$ is merely a map of pointed sets, and that it need not be injective. (Take $G$ to be a special orthogonal group.)
\begin{proof}
We have an exact sequence 
$$
1 \rightarrow \mu \rightarrow G_{sc} \rightarrow G \rightarrow 1
$$
where $G_{sc}$ is simply connected and $\mu$ is a finite subgroup of a split torus (i.e. a product of $\mu_n$'s).
This sequence induces a long exact sequence in nonabelian cohomology.
Let $\alpha \in \HH^1(k,G)$ be a class with $\alpha_v = 1$ for all places $v$ of $k$. 
Since $\HH^2(k,\mu) \rightarrow \prod_v \HH^2(k_v,\mu)$ is injective by the Hasse principle for the Brauer group, we see that $\alpha$ lifts to a class $\beta \in \HH^1(k,G_{sc})$. 
Since $\mu$ is a central subgroup of $G_{sc}$, any other lift of $\alpha$ is given by $\lambda \beta$, where $\lambda \in \HH^1(k,\mu)$ is a cocycle. 
We will show that we can choose $\lambda$ so that $\lambda \beta$ is trivial.
By the Hasse principle for simply connected groups \cite[Theorem 6.6]{PlatonovRapinchuk-Alggroupsandnumbertheory}, the map $\HH^1(k,G_{sc}) \rightarrow \prod_v \HH^1(k_v,G_{sc})$ is injective. (This map is even bijective.)
If $v$ is a finite or complex place, then $\HH^1(k_v,G_{sc})$ is trivial \cite[Theorem 5.12.24(b)]{Poonen-rationalpointsonvarieties}.
If $v$ is real then $\beta_{v} \in \HH^1(k_v,G_{sc})$ has trivial image in $\HH^1(k_v,G)$ so comes from an element of $\HH^1(k_v, \mu)$.
Since $\HH^1(k,\mu) \rightarrow \prod_{v\text{ real}} \HH^1(k_v,\mu)$ is surjective (this follows from the case $\mu = \mu_n$), we may choose $\lambda \in \HH^1(k,\mu)$ such that $\lambda_v \beta_v =1$ for every real place $v$.
This implies that $\lambda \beta $ is trivial, as required.
\end{proof}

\begin{corollary}\label{corollary: inject 2-Selmer orbits}
Let $k$ be a number field and $b\in B^{\rs}(k)$. Let $\Sel_2 J_b$ be the $2$-Selmer group of the abelian variety $J_b/k$. 
Then $\Sel_2 J_b \subset \HH^1(k,J_b[2])$ is contained in the image of $\gamma_b$.
Consequently, the injection $\eta_b$ from Theorem \ref{theorem: inject 2-descent orbits} extends to an injection $$\Sel_2 J_b \hookrightarrow G(k) \setminus V_b(k).$$
\end{corollary}
\begin{proof}
We have a commutative diagram for every place $v$: 
	\begin{center}
		\begin{tikzcd}
			{J_b(k)/2\Jac_b(k)} \arrow[d] \arrow[r, "\delta"] & {\HH^1(k,\Jac_b[2])} \arrow[r] \arrow[d] & {\HH^1(k,\bigG)} \arrow[d] \\
			{\Jac_b(k_v)/2\Jac_b(k_v)} \arrow[r, "\delta_v"]     & {\HH^1(k_v,\Jac_b[2])} \arrow[r]         & {\HH^1(k_v,\bigG)}        
		\end{tikzcd}
	\end{center}
	By Corollary \ref{corollary: AIT for orbits of V} it suffices to prove that $2$-Selmer elements in $\HH^1(k,\Jac_b[2])$ are killed under the composition $\HH^1(k,\Jac_b[2]) \xrightarrow{\sim} \HH^1(k,Z_G(\kappa_b)) \rightarrow \HH^1(k,G)$. 
	By definition, an element of $\Sel_2 \Jac_b$ consists of a class in $\HH^1(k,\Jac_b[2])$ whose restriction to $\HH^1(k_v,\Jac_b[2])$ lies in the image of $\delta_v$ for every place $v$. 
	So by Theorem \ref{theorem: inject 2-descent orbits} the image of such an element in $\HH^1(k_v,\bigG)$ is trivial for every $v$. Proposition \ref{proposition: weak hasse principle for H1(G)} completes the proof.
\end{proof}

\subsection{Reducible orbits and marked points}\label{subsection: trivial orbits}

Recall from Definition \ref{definition: k-reducible orbits} that an element of $V^{\rs}(k)$ is $k$-reducible if it is $G(k)$-conjugate to a Kostant section.
Recall from \S\ref{subsection: A family of curves} that $\infty_1,\dots,\infty_m$ denote the set of marked points of $C\rightarrow B$.

\begin{proposition}\label{proposition: orbits corresponding to marked points are reducible}
Let $k/\Q$ be a field and $b\in B^{\rs}(k)$. 
Then the image under $\eta_b\colon J_b(k)/2J_b(k) \hookrightarrow G(k) \backslash V_b(k)$ of the subgroup of $J_b(k)/2J_b(k)$ generated by $\{\infty_2-\infty_1,\dots,\infty_m-\infty_1\}$ coincides with the set of $k$-reducible $G(k)$-orbits of $V_b(k)$.
Moreover, the set of $k$-reducible $G(k)$-orbits has the maximal size $2^{m-1}$ if and only if the inclusion $Z_G(\kappa_b) \subset Z_{H^{\theta}}(\kappa_b)$ is surjective on $k$-points.
\end{proposition}
\begin{proof}
The proof is very similar to \cite[Lemma 2.11]{Romano-Thorne-ArithmeticofsingularitiestypeE}.
For a scheme $X/k$ we write $\HH_1(X,\F_2) \coloneqq \Hom(\HH^1_{et}(X_{k^s},\F_2),\F_2)$, where $\HH^1_{et}$ denotes \'etale cohomology. We have an exact sequence of \'etale homology groups
\begin{align}\label{equation: sequence etale homology groups}
1\rightarrow \mu_2^m/\Delta(\mu_2) \rightarrow \HH_1(C_b^{\circ},\F_2) \rightarrow \HH_1(C_b,\F_2)\rightarrow 1.
\end{align}
Let $H_{sc}\rightarrow H$ be the simply connected cover of $H$ and let $C_{H_{sc}}$ be the centre of $H_{sc}$. 
By \cite[Theorem 4.10]{Thorne-thesis}, the sequence \eqref{equation: sequence etale homology groups} is isomorphic to 
\begin{align}\label{sequence: sequence centraliser sc kostant section}
1\rightarrow C_{H_{sc}}[2] \rightarrow Z_{H_{sc}^{\theta}}(\kappa_b) \rightarrow Z_G(\kappa_b) \rightarrow 1.
\end{align}
It follows that the duals of these sequences are also isomorphic. We will calculate these duals and their connecting maps in Galois cohomology.

The dual of \eqref{equation: sequence etale homology groups} is isomorphic to 
\begin{align}\label{equation: sequence 2-torsion open curve}
1\rightarrow J_b[2] \rightarrow \HH^1_{et}(C_{b,k^s}^{\circ},\F_2)\rightarrow (\mu_2^m)_{\Sigma = 0} \rightarrow 1.
\end{align}
Here we use the identification $J_b[2] = \HH^1_{et}(C_{b,k^s},\F_2)$, and $(\mu_2^m)_{\Sigma = 0}$ denotes the subset of $\mu_2^m$ of elements summing to zero.
An explicit calculation shows that the image of the connecting map $(\mu_2^m)_{\Sigma = 0}(k) \rightarrow \HH^1(k,J_b[2])$ coincides with the image of the subgroup of $J_b(k)/2J_b(k)$ generated by $\{\infty_2-\infty_1,\dots,\infty_m-\infty_1\}$ under the $2$-descent map $J_b(k)/2J_b(k) \hookrightarrow \HH^1(k,J_b[2])$.

On the other hand, we claim that the dual of \eqref{sequence: sequence centraliser sc kostant section} is isomorphic to
\begin{align}\label{equation: sequence centralisers of kostant section in H^theta}
1\rightarrow Z_G(\kappa_b) \rightarrow Z_{H^{\theta}}(\kappa_b) \rightarrow \pi_0(H^{\theta})\rightarrow 1.
\end{align}
Indeed, the identification of the first two terms follows from \cite[Corollary 2.9]{Thorne-thesis} and the existence of a nondegenerate pairing on $Z_G(\kappa_b)$ \cite[Corollary 2.12]{Thorne-thesis}.
It follows from Lemma \ref{lemma: component group meets Cartan maximal torus} that we may identify the last term with $\pi_0(H^{\theta})$. 
Next, we claim that the image of the connecting map $\pi_0(H^{\theta}) \rightarrow \HH^1(k,Z_G(\kappa_b))$ coincides with the image of the $k$-reducible orbits in $V_b(k)$ under the map $G(k) \backslash V_b(k) \hookrightarrow \HH^1(k,Z_G(\kappa_b))$ from Lemma \ref{lemma: AIT full generality}.
Indeed, consider the commutative diagram
\begin{center}
\begin{tikzcd}
G(k) \backslash V_b(k) \arrow[r] \arrow[hookrightarrow]{d} & H^{\theta}(k) \backslash V_b(k) \arrow[hookrightarrow]{d} \\ 
\HH^1(k,Z_G(\kappa_b)) \arrow[r] & \HH^1(k,Z_{H^{\theta}}(\kappa_b))
\end{tikzcd}
\end{center}
where the horizontal maps are induced by the inclusions $G\subset H^{\theta}$ and $Z_G(\kappa_b)\subset Z_{H^{\theta}}(\kappa_b)$, and the vertical maps arise from Lemma \ref{lemma: AIT full generality}.
It follows from Corollary \ref{corollary: pi0 is constant} that the map $\HH^1(k,G) \rightarrow \HH^1(k,H^{\theta})$ has trivial pointed kernel. 
Moreover all $k$-reducible elements in $V_b(k)$ are $H^{\theta}(k)$-conjugate by Proposition \ref{proposition: H^theta acts simply transitively regular nilpotents}.
Therefore the set of $k$-reducible $G(k)$-orbits corresponds to the kernel of $\HH^1(k,Z_G(\kappa_b)) \rightarrow \HH^1(k,Z_{H^{\theta}}(\kappa_b))$ which, using \eqref{equation: sequence centralisers of kostant section in H^theta}, coincides with the image of the map $\pi_0(H^{\theta})\rightarrow \HH^1(k,Z_G(\kappa_b))$.
This proves the claim and the first part of the proposition. 

To prove the remaining part, note that there are $2^{m-1}$ $k$-reducible orbits if and only if the map $(\mu_2^m)_{\Sigma = 0}(k) \rightarrow \HH^1(k,J_b[2])$ is injective.
By considering the long exact sequences associated with the isomorphic sequences \eqref{equation: sequence 2-torsion open curve} and \eqref{equation: sequence centralisers of kostant section in H^theta}, this is equivalent to the surjectivity of $\HH^0(k,Z_G(\kappa_b))\rightarrow \HH^0(k,Z_{H^{\theta}}(\kappa_b))$.
\end{proof}

\section{Integral representatives}\label{integral representatives}

We keep the notations from Section \ref{section: recollections of Thorne's thesis}.
In this chapter we introduce integral structures for $G$ and $V$ and prove that for large primes $p$, the image of the map from Theorem \ref{theorem: inject 2-descent orbits} applied to $R = \Q_p$ lands in the orbits which admit a representative in $\Z_p$.
See Theorem \ref{theorem: integral representatives exist} for a precise statement.
In \S\ref{subsection: orbits over Z}, we deduce an integrality result for orbits over $\Q$ (as opposed to orbits over $\Q_p$).

\subsection{Integral structures}\label{subsection: integral structures}

So far we have considered properties of the pair $(G,V)$ over $\Q$. In this subsection we define these objects over $\Z$.


Let $\intbigH$ (respectively $\intbigG$) be the unique (up to isomorphism) split reductive group over $\Z$ with generic fibre $H$ (respectively $G$). 
The automorphism $\bigtheta\colon \bigH \rightarrow \bigH$ extends by the same formula to an automorphism $\intbigH\rightarrow \intbigH$, still denoted by $\bigtheta$. 


\begin{lemma}
    The equality $(\bigH^{\bigtheta})^{\circ}=\bigG$ extends to an isomorphism $(\intbigH^{\bigtheta}_{\Z[1/2]})^{\circ}\simeq \intbigG_{\Z[1/2]}$, where $(\intbigH^{\bigtheta}_{\Z[1/2]})^{\circ}$ is the relative identity component of $\intbigH^{\bigtheta}_{\Z[1/2]}$.
\end{lemma}
\begin{proof}
   This follows from the fact that $(\intbigH^{\bigtheta}_{\Z[1/2]})^{\circ}$ is a reductive group scheme of the same type as $\intbigG_{\Z[1/2]}$, which follows from \cite[Remark 3.1.5]{Conrad-reductivegroupschemes}. 
\end{proof}

To obtain a $\Z$-structure for $V$, choose a Chevalley basis of $\lieg$ with respect to the maximal torus $T^{\theta}$, and choose an admissible $\Z$-form $\intbigV$ of the $G$-representation $V$ with respect to this basis \cite[Proposition 2.4]{Borel-propertieschevalley}.
Since $G$ acts faithfully on $V$, the results of \cite[\S3]{Borel-propertieschevalley} imply that $\intbigG$ is isomorphic to the Zariski closure of $G$ in $\GL(\intbigV)$.
We henceforth view $\intbigG$ as a closed subgroup of $\GL(\intbigV)$ so $\intbigV$ is a representation of $\intbigG$.

Recall from \S\ref{subsection: A family of curves} that we have fixed polynomials $p_{d_1},\dots,p_{d_r} \in \Q[\bigV]^{\bigG}$ satisfying the conclusions of Proposition \ref{proposition: first properties of the family of curves}. 
Note that those conclusions are invariant under the $\G_m$-action on $\bigB$. 
By rescaling the polynomials $p_{d_i}$ using this $\G_m$-action, we can assume they lie in $\Z[\intbigV]^{\intbigG}$. 
We may additionally assume that the discriminant $\Delta$ from \S\ref{subsection: the discriminant locus} lies in $\Z[\intbigV]^{\intbigG}$.
Define $\intbigB \coloneqq \Spec \Z[p_{d_1},\dots,p_{d_r}]$ and $\intbigB^{\rs} \coloneqq \Spec \Z[p_{d_1},\dots,p_{d_r}][\Delta^{-1}]$.
Taking invariants defines a map $\pi \colon \intbigV \rightarrow \intbigB $.

We extend the family of curves given by the equation in Table \ref{table: introduction different cases} to the family $\intbigcurve\rightarrow \intbigB$ given by that same equation.

\begin{proposition}\label{proposition: G has class number one}
    $\intbigG$ has class number $1$: $\bigG(\A^{\infty}) = \bigG(\Q)\intbigG(\widehat{\Z})$.
\end{proposition}
\begin{proof}
    The group $\intbigG$ is the Zariski closure of $\bigG$ in $\GL(\intbigV)$ and in a suitable basis of $\intbigV$, $\bigG$ contains a maximal $\Q$-split torus $T^{\theta}$ consisting of diagonal matrices in $\GL(V)$. 
	Therefore $\intbigG$ has class number $1$ by \cite[Theorem 8.11; Corollary 2]{PlatonovRapinchuk-Alggroupsandnumbertheory} and the fact that $\Q$ has class number one. 
\end{proof}

\subsection{Spreading out}\label{subsection: spreading out}

Our constructions and theorems for $(G,V)$ of the previous chapters will continue to be valid over $\Z[1/N]$ for some appropriate choice of integer $N$, in a sense we will now explain.

Let us call a positive integer $N$ \define{admissible} if the following properties are satisfied (set $S \coloneqq \Z[1/N]$): 

\begin{enumerate}
	\item\label{item: order weyl group invertible} Each prime dividing the order of the Weyl group of $H$ is a unit in $S$. (In particular, $2$ is a unit in $S$.)
	\item\label{item: spreading out discriminant} The zero locus $\underline{D}_S\rightarrow \Spec S$ of the discriminant $\Delta$ is flat and its smooth locus $\underline{D}_S^1$ coincides with the regular locus of $\underline{D}_S$. Moreover, the nonsmooth locus of $\underline{D}_S\rightarrow \Spec S$ is flat over $\Spec S$.
	\item\label{item: smooth and locus of curves} The morphism $\intbigcurve_S\rightarrow \intbigB_S$ is smooth exactly above $\intbigB_S^{\rs}$.
	\item\label{item: action map subregular slice is smooth} The affine curve $\intbigcurve_S^{\circ}$ is a closed subscheme of $\intbigV_S$ and the action map $\intbigG_S \times \intbigcurve_S^{\circ} \rightarrow \intbigV_S, (g,x)\mapsto g\cdot x$ is smooth.
	\item\label{item: characterisation nodal curves D1} For a field $k$ of characteristic not dividing $N$, $b \in \underline{D}^1(k)$ if and only if every semisimple lift $v\in \intbigV_b(k)$ has centraliser $Z_H(v)$ of semisimple rank 1, if and only if the curve $\intbigcurve_b$ has a unique nodal singularity. In that case, the curve $\intbigcurve_b$ is geometrically integral.
	\item There exist open subschemes $\intbigV^{\rs} \subset \intbigV^{\reg} \subset \intbigV_S$ such that if $k$ is a field of characteristic not dividing $N$ and $v\in \intbigV(k)$ then $v$ is regular if and only if $v\in \intbigV^{\reg}(k)$, and $v$ is regular semisimple if and only if $v\in \intbigV^{\rs}(k)$. Moreover, $\intbigV^{\rs}$ is the open subscheme defined by the nonvanishing of the discriminant polynomial $\Delta$ in $\intbigV_S$. 
	\item\label{item: smooth locus of V-> B} The morphism $\pi\colon \intbigV_S\rightarrow \intbigB_S$ is smooth exactly at $\intbigV^{\reg}$.
	\item\label{item: invariants and kostant section} $S[\intbigV]^{\intbigG} = S[p_{d_1},\dots,p_{d_r}]$. The Kostant section $\kappa$ fixed in \S\ref{subsection: A family of curves} extends to a section $\kappa\colon \intbigB_S \rightarrow \intbigV^{\reg}$ of $\pi$ satisfying the following property: for any $b\in \intbigB(\Z) \subset \intbigB_S(S)$, we have $\kappa_{N\cdot b} \in \intbigV(\Z)$. Moreover, each $G(\Q)$-orbit of Kostant sections has a representative which satisfies the same property.
	\item\label{item: image action map kostant section} Let $\intbigB^1$ be the complement of the nonregular locus of $\underline{D}$ in $\intbigB$. Then the action map $\intbigG_S \times \intbigB_S \rightarrow \intbigV^{\reg}, (g,b) \mapsto g\cdot \kappa_b$ is \'etale and its image contains $\intbigV^{\reg}|_{\intbigB_S^1}$.
	\item\label{item: isomorphism centraliser and two torsion} Let $\underline{J}^1_S\rightarrow \intbigB_S$ denote the relative generalised Jacobian of the family of integral curves $\intbigcurve_S|_{\intbigB_S^1} \rightarrow \intbigB_S^1$ \cite[\S9.3, Theorem 1]{BLR-NeronModels} and let $\underline{J}^{\rs}_S \rightarrow \intbigB_S^{\rs}$ denote its restriction to $\intbigB_S^{\rs}$. 
	Let $\underline{Z}_S\rightarrow \intbigB_S$ be the centraliser of the Kostant section $\kappa$ in $\intbigG_S$.
	Then there is an isomorphism $\underline{J}_S^{\rs}[2] \simeq \underline{Z}_S^{\rs}$ of finite \'etale group schemes over $\intbigB_S^{\rs}$ whose restriction to $\bigB^{\rs}$ is the isomorphism of Proposition \ref{proposition: isomorphism 2-torsion and Kostant section centraliser thorne}. It extends to an isomorphism $\underline{J}_S^1[2] \simeq \underline{Z}_S^1$. 
	\item\label{item: compactified jacobian spreading out} The $B$-scheme $\bar{J}$ constructed in \S\ref{section: the compactified jacobian} extends to a $\intbigB_S$-scheme $\underline{\bar{J}}_S\rightarrow \intbigB_S$ which is flat, projective, with geometrically integral fibres and whose restriction to $\intbigB_S^{\rs}$ is isomorphic to $\underline{J}_S^{\rs}$. 
	Moreover, $\underline{\bar{J}}_S\rightarrow S$ is smooth with geometrically integral fibres, and the smooth locus of the morphism $\underline{\bar{J}}_S\rightarrow \intbigB_S$ is an open subscheme of $\underline{\bar{J}}_S$ whose complement is $S$-fibrewise of codimension at least two. 
	\item\label{item: torsor for orbit stuff} The $G$-torsor $T\rightarrow J^{\rs}$ from \S\ref{subsection: a universal torsor} extends using the same definition to a $\intbigG_S$-torsor $\underline{T}_S \rightarrow \underline{J}^{\rs}_S $, and $\underline{T}_S$ is Zariski locally trivial.
\end{enumerate}

It might be possible to construct an explicit admissible integer for every pair $(G,V)$.
We will content ourselves with the following: 

\begin{proposition}
    There exists an admissible integer $N$.
\end{proposition}
\begin{proof}
    The proof is very similar to the proof of \cite[Proposition 4.1]{Laga-E6paper}.
    It follows from the results of the previous chapters and the principle of spreading out \cite[\S3.2]{Poonen-rationalpointsonvarieties}.
    We omit the details, but refer in each case to the corresponding property over $\Q$.

	Properties 1 and 2 follow from spreading out;
	Property 3 follows from Lemma \ref{lemma: singularities of C_b using dynkin diagram};
	Property 4 from the definition of $C^{\circ}$ in \S\ref{subsection: A family of curves} and Part 6 of Proposition \ref{proposition: first properties of the family of curves};
	Property 5 follows from Theorem \ref{theorem: summary properties mildly singular locus};
	Property 6 follows from an argument similar to Property 5 of \cite[Proposition 4.1]{Laga-E6paper};
	Property 7 follows from Lemma \ref{lemma: smooth locus of invariant map};
	Property 8 follows from an argument similar to Property 4 of \cite[Proposition 4.1]{Laga-E6paper};
	Property 9 follows from Propositions \ref{proposition: Kostant section properties} and \ref{proposition: action map is surjective over B1};
	Property 10 follows from Proposition \ref{proposition: isomorphism 2-torsion and Kostant section centraliser thorne} and Theorem \ref{theorem: isomorphism centraliser 2-torsion over B1};
	Property 11 follows from Theorem \ref{theorem: summary compactified jacobian};
	finally Property 12 follows from Theorem \ref{theorem: universal torsor is Zariski trivial}.
\end{proof}


For the remainder of the paper, we fix an admissible integer $N$ and continue to write $S=  \Spec \Z[1/N]$. 
Moreover, to simplify notation, we will drop the subscript $\underline{(\,)}_S$ and write $G, V, B, J, C\dots$ for $\intbigG_S, \intbigV_S, \intbigB_S, \underline{J}_S, \intbigcurve_S \dots$. 

Using these properties, we can extend our previous results to $S$-algebras rather than $\Q$-algebras. 
We mention in particular the following two, which follow from Properties \ref{item: image action map kostant section}, \ref{item: isomorphism centraliser and two torsion} and \ref{item: torsor for orbit stuff}.
(Recall from \S\ref{subsection: notation} the definition of $\HH^1$.)

\begin{proposition}[Analogue of Corollary \ref{corollary: AIT for orbits of V}]\label{proposition: spread out orbit parametrization galois}
Let $R$ be an $S$-algebra and $b\in B^{\rs}(R)$. Then we have a natural bijection of pointed sets:
		\begin{align}
		G(R) \backslash V_b(R) \simeq  \ker\left(\HH^1(R,J_b[2]) \rightarrow \HH^1(R,G)\right).
		\end{align}	
\end{proposition}

\begin{proposition}[Analogue of Theorem \ref{theorem: inject 2-descent orbits}]\label{proposition: inject 2-descent orbits spreading out}
Let $R$ be a local $S$-algebra and $b\in B^{\rs}(R)$. 	
Then there is an injective map 
\begin{align*}
\eta_b\colon J_b(R)/2J_b(R) \hookrightarrow G(R) \backslash V_b(R)
\end{align*}
compatible with base change on $R$.
\end{proposition}

We are now ready to state the main theorem of this chapter.
Write $\sh{E}_p$ for the set of all $b\in \intbigB(\Z_p)$ that lie in $B^{\rs}(\Q_p)$.
It consists of those elements of $\intbigB(\Z_p)$ of nonzero discriminant.
(This is different to $\intbigB^{\rs}(\Z_p)$, which consists of those $b\in \intbigB(\Z_p)$ whose discriminant is a \emph{unit} in $\Z_p$ and for which integral orbits are already constructed in Proposition \ref{proposition: inject 2-descent orbits spreading out}.)

\begin{theorem}\label{theorem: integral representatives exist}
Let $p$ be a prime not dividing $N$. 
Then for any $b\in \sh{E}_p$, the image of the map 
\begin{align*}
\eta_b\colon J_b(\Q_p)/2J_b(\Q_p) \rightarrow G(\Q_p) \backslash V_b(\Q_p)
\end{align*}
from Theorem \ref{theorem: inject 2-descent orbits} is contained in the image of the map $V(\Z_p) \rightarrow G(\Q_p) \backslash V(\Q_p)$.
\end{theorem}
The proof of Theorem \ref{theorem: integral representatives exist} will be given at the end of \S\ref{subsection: proof of integral representatives general case}. 
The essential ingredients are Lemma \ref{lemma: purity for the stack M}, Lemma \ref{lemma: squarefree discriminant implies curve is regular} and the good geometric properties of the compactified Jacobian summarized in Theorem \ref{theorem: summary compactified jacobian}.
We refer to the start of \S\ref{subsection: proof of integral representatives general case} for a general overview of the proof strategy.

\subsection{Some stacks}\label{subsection: some stacks}



For technical purposes related to the proof of Theorem \ref{theorem: integral representatives exist}, we need to introduce some stacks relevant to our setup.
This can be seen as an attempt to `geometrise' the set of $G$-orbits of $V$, and allows for more flexibility in glueing and descent arguments.
Hopefully we soothe the reader by mentioning that we will not need any serious properties of stacks, and we mainly think of them as collections of groupoids where one can glue objects suitably.
All stacks introduced in this paper are considered in the \'etale topology. 
Recall from \S\ref{subsection: spreading out} that we have fixed an admissible integer $N$ and we have set $S = \Z[1/N]$.

\begin{definition}\label{definition: classifying stack of G}
Let $\Bun{G} = [\Spec S/G]$ be the \define{classifying stack} of $G$.
By definition, for any $S$-scheme $X$ the groupoid $\Bun{G}(X)$ has as objects $G$-torsors over $X$.
Morphisms are given by isomorphisms of $G$-torsors.
\end{definition}

\begin{definition}\label{definition: quotient stack V/G}
Let $\mathcal{M} = \left[ G\backslash V \right]$ be the quotient stack of $V$ by the natural $G$-action on $V$.
By definition, for any $S$-scheme $X$ an object of $\mathcal{M}(X)$ consists of a $G$-torsor $T\rightarrow X$ together with a $G$-equivariant morphism $\phi\colon T \rightarrow V$. 
A morphism between two objects $(T,\phi)$ and $(T',\phi')$ consists of an isomorphism $\alpha \colon T \rightarrow T'$ of $G$-torsors satisfying $\phi' \circ \alpha = \phi$.
\end{definition}

Finally, recall that $Z\rightarrow B$ denotes the centraliser of the Kostant section $\kappa$, an extension of the group scheme of Definition \ref{definition: universal centralisers} to $S$.
Consider the quotient stack $\left[ B /Z \right]\rightarrow B$, where $Z$ acts trivially on $B$. For any $B$-scheme $X$, an $X$-point of $\left[B/Z\right](X)$ consists of a $Z$-torsor on $X$. 

These stacks come with a few natural maps between them:
\begin{itemize}
    \item $\mathcal{M}\rightarrow \Bun{G}$: sends a pair $(T,\phi)$ to the $G$-torsor $T$.
    \item $\mathcal{M} \rightarrow B$: sends a pair $(T\xrightarrow{\alpha} X, T\xrightarrow{\phi} V)$ to the unique morphism $X\xrightarrow{f} B$ fitting in the commutative diagram:
    \begin{center} 
    \begin{tikzcd}
    T \arrow[r, "\phi"] \arrow[d,"\alpha"] & V \arrow[d, "\pi"] \\
    X \arrow[r,"f"] & B 
    \end{tikzcd}
    \end{center}
    (Here $\pi$ denotes the invariant map, and the existence and uniqueness of $f$ follows from \'etale descent.)
    We will often regard $\mathcal{M}$ as a stack over $B$. 
    In particular, if $b\in B(X)$ is an $X$-point we write $\mathcal{M}_b$ for the pullback of $\mathcal{M}$ along this point; it is isomorphic to $\left[ G \backslash V_b \right]$.
    \item $V\rightarrow \mathcal{M}$: sends an $X$-point $X\xrightarrow{v} V$ to $(G\times X, \phi_v)$, where $\phi_v\colon G\times X\rightarrow V$ sends $(g,x)$ to $g\cdot v(x)$.
    \item There is a substack $\left[B/Z\right] \hookrightarrow \mathcal{M}$ obtained by `twisting' the Kostant section. For any $B$-scheme $X$, its image consists of those elements of $\mathcal{M}$ that are \'etale locally conjugate to $\kappa_b$ (or rather to its image under $V\rightarrow \mathcal{M}$.)
\end{itemize}

If $\mathcal{G}$ is a groupoid, we write $\pi_0\mathcal{G}$ for its set\footnote{Assuming it is a set, which will always be the case in this paper.} of isomorphism classes.
\begin{lemma}\label{lemma: points on quotient stack with trivial torsor are just orbits}
Let $b$ be an $X$-point of $B$. The map $V_b(X) \rightarrow \mathcal{M}_b(X)$ induces a bijection between the $G(X)$-orbits of $V_b(X)$ and elements of $\pi_0(\mathcal{M}_b(X))$ that map to the trivial element in $\pi_0(\Bun{G}(X))$. 
\end{lemma}
\begin{proof}
This follows formally from the definitions. 
Indeed, if $v,v'\in V_b(X)$ give rise to isomorphic elements $(G\times X, \phi_v)$ and $(G\times X, \phi_{v'})$ in $\mathcal{M}_b(X)$, then there exists an isomorphism $G \times X \xrightarrow{\sim} G\times X$ of $G$-torsors mapping $\phi_v$ to $\phi_{v'}$. 
Such an isomorphism is defined by multiplying an element of $G(X)$, so $v$ and $v'$ are $G(X)$-conjugate.
The argument can be reversed, so we obtain an injection $G(X) \backslash V_b(X)\hookrightarrow \pi_0(\mathcal{M}_b(X))$. 
Since an object $(T,\phi)$ of $\mathcal{M}_b(X)$ is isomorphic to $(G\times X, \phi_v)$ for some $v\in V_b(X)$ if and only if $T$ is the trivial torsor, we conclude.
\end{proof}

\begin{example}\label{example: Zp-points quotient stack}
    Suppose that $X = \Spec \Z_p$ where $p$ is coprime to $N$ and $b\in B(\Z_p)$.
    Since every $G$-torsor over $\Spec \Z_p$ is trivial (by \cite[III.3.11(a)]{milne-etalecohomology} and Lang's theorem), Lemma \ref{lemma: points on quotient stack with trivial torsor are just orbits} gives a bijection between $G(\Z_p) \backslash V_b(\Z_p)$ and $\pi_0(\mathcal{M}_b(\Z_p))$.
    Therefore when we want to show that an orbit $x\in G(\Q_p) \backslash V_b(\Q_p)$ has an integral representative, it suffices to extend it to a $\Z_p$-point of quotient stack $\mathcal{M}_b$.
    This will be used in the proof of Theorem \ref{theorem: integral representatives exist}.
\end{example}

The next lemma can be interpreted as a categorical version of Corollary \ref{corollary: AIT for orbits of V}.

\begin{lemma}\label{lemma: stacks over Brs are the same}
The inclusion $[B/Z] \hookrightarrow \mathcal{M}$ induces an isomorphism of stacks $[B^{\rs}/Z^{\rs}] \simeq \mathcal{M}^{\rs}$ over $B^{\rs}$.
\end{lemma}
\begin{proof}
    It suffices to prove that for any $B$-scheme $X$ and $b\in B^{\rs}(X)$, every two objects in $\mathcal{M}_b(X)$ are \'etale locally isomorphic. (For then every object will be \'etale locally isomorphic to the Kostant section.)
    By passing to an \'etale extension, we may assume that these objects map to the trivial element in $\pi_0(\Bun{G}(X))$.
    It therefore suffices to prove that every two elements of $V^{\rs}_b(X)$ are \'etale locally $G(X)$-conjugate. 
    This is true, since $G\times B^{\rs} \rightarrow V^{\rs}$ is smooth and surjective so has sections \'etale locally; see Part \ref{item: image action map kostant section} of \S\ref{subsection: spreading out}.
\end{proof}

Let $\mathcal{M}^{\reg}\subset \mathcal{M}$ be the open substack whose $X$-points consist of those objects $(T,\phi)$ of $\mathcal{M}(X)$ such that $\phi$ lands in the locus of regular elements $V^{\reg}$, and all morphisms between them.
Note that the map $[B/Z]\rightarrow \mathcal{M}$ factors through $\mathcal{M}^{\reg}$ by (a spreading out of) Part 2 of Proposition \ref{proposition: Kostant section properties}. 

\begin{lemma}\label{lemma: stacks over B1 are the same}
The inclusion $[B/Z] \hookrightarrow \mathcal{M}^{\reg}$ induces an isomorphism of stacks $[B^1/Z^1] \simeq \mathcal{M}^{\reg}|_{B^1}$ over $B^1$.
\end{lemma}
\begin{proof}
By the same reasoning as the proof of Lemma \ref{lemma: stacks over Brs are the same}, it suffices to prove that $G\times B^1\rightarrow V^{1,reg}$ is smooth and surjective. 
This follows from Part \ref{item: image action map kostant section} of \S\ref{subsection: spreading out}, which is a spreading out of Propositions \ref{proposition: Kostant section properties} and \ref{proposition: action map is surjective over B1}.
\end{proof}

The next lemma is a purity result for the stack $\mathcal{M}$, and will be a crucial ingredient in the proof of Theorem \ref{theorem: integral representatives exist}.

\begin{lemma}\label{lemma: purity for the stack M}
Let $X$ be a regular integral $2$-dimensional scheme, let $U\subset X$ be an open subscheme whose complement is finite and let $b\in B(X)$. 
Then the restriction $\mathcal{M}_b(X) \rightarrow \mathcal{M}_{b|_U}(U)$ is an equivalence of categories.
\end{lemma}
\begin{proof}
We will use the following fact \cite[Lemme 2.1(iii)]{ColliotTheleneSansuc-Fibresquadratiques} repeatedly: if $Y$ is an affine $X$-scheme of finite type, then restriction of sections $Y(X)\rightarrow Y(U)$ is bijective.
To prove essential surjectivity, let $(T_U, T_U\xrightarrow{\phi_U} V_b)$ be an object of $\mathcal{M}_{b|_U}(U)$.
By \cite[Th\'eoreme 6.13]{ColliotTheleneSansuc-Fibresquadratiques}, the $G$-torsor $T_U\rightarrow U$ extends to a $G$-torsor $T$ on $X$.
By the fact above applied to $Y = V_b$, $\phi_U$ uniquely extends to a morphism $\phi \colon T\rightarrow V_b$.
The uniqueness of $\phi$ guarantees that $\phi$ is $G$-equivariant.
Since the scheme of isomorphisms $\Isom_{\mathcal{M}}(\mathcal{A},\mathcal{A}')$ between two objects of $\mathcal{M}_b(X)$ is $X$-affine, fully faithfulness follows from the fact applied to this isomorphism scheme.
\end{proof}

\subsection{Orbits of square-free discriminant}\label{subsection: square-free discriminant case}





In this subsection we study orbits of square-free discriminant, which will be useful in the proof of Theorem \ref{theorem: integral representatives exist} and in Chapter \ref{section: the average size of the 2-Selmer group}.
For the remainder of \S\ref{subsection: square-free discriminant case}, we fix a discrete valuation ring $R$ with fraction field $K$, uniformiser $\pi$, residue field $k$ and normalised discrete valuation $\ord_K: K^{\times} \twoheadrightarrow \Z$.
We assume that the integer $N$ fixed in \S\ref{subsection: spreading out} is a unit in $R$. 

\begin{lemma}\label{lemma: squarefree discriminant implies special fibre in D1}
Suppose that $b\in B(R)$ satisfies $\ord_K(\Delta(b)) = 1$. Then $b_k \in D^1(k)$.
\end{lemma}
\begin{proof}
We may assume that $R$ is complete.
Since $\Delta(b)$ reduces to $0\in k$, we have $b_k \in D$. 
Since $\Delta(b)$ is a uniformiser of $R$, the quotient of the regular ring $\HH^0(B_R,\O_{B_R}) = R[p_{d_1},\dots,p_{d_r}]$ by the maximal ideal $(p_{d_1}-p_{d_1}(b),\dots,p_{d_r}-p_{d_r}(b),\Delta)$ is isomorphic to $k$.
Therefore the elements $\{p_{d_1}-p_{d_1}(b),\dots,p_{d_r}-p_{d_r}(b),\Delta\}$ form a regular system of parameters at $b_k$.
Hence $D_R = \Spec R[\{p_{d_i}\}] / (\Delta)$ is regular at $b_k$.
We conclude that $b_k$ is a regular point of $D_R$. 
To prove the lemma, it suffices to prove that the regular locus of $D_R$ coincides with the smooth locus $D^1_R$ of $D_R\rightarrow \Spec R$. 

Indeed, let $Z\subset D$ be the nonsmooth locus of $D\rightarrow \Spec S$, which coincides with the nonregular locus of $D$ and is $S$-flat by Part \ref{item: spreading out discriminant} of \S\ref{subsection: spreading out}.
Let $Z'$ be the nonregular locus of $D_R$, which is closed by the excellence of $R$.
Since taking the smooth locus commutes with base change, $Z_R$ agrees with the nonsmooth locus of $D_R\rightarrow \Spec R$. 
It therefore suffices to show that $Z_R = Z'$. 
Since every smooth point of $D_R\rightarrow \Spec R$ is regular, $Z' \subset Z_R$.
To prove the opposite inclusion, let $L$ be the prime subfield of $K$, in other words the residue field of the image of $\Spec K \rightarrow \Spec S$.
Since $L$ is perfect, $Z_L$ is the nonregular locus of $D_L$ and every point of $D_L \setminus Z_L$ is geometrically regular.
It follows that $Z_K$ is the nonregular locus of $D_K$ and so $Z_K = Z'_K$. 
Since $Z$ is $S$-flat by assumption, $Z_R$ is $R$-flat so $Z_K$ is dense in $Z_R$. 
Since $Z'$ contains the closure of $Z_K'=Z_K$ which is $Z_R$, we conclude that $Z_R\subset Z'$.
\end{proof}

\begin{lemma}\label{lemma: squarefree discriminant implies Vb is regular}
Suppose that $b\in B(R)$ satisfies $\ord_K(\Delta(b)) =1$.
Then the scheme $V_b$ is regular. 
\end{lemma}
\begin{proof}
The idea of the proof is to reduce the statement to $\liesl_2$; this will be achieved by a sequence of standard but somewhat technical reduction steps.
Since regularity can be checked after \'etale extensions and completion, we may assume that $R$ is complete and $k$ is separably closed.
By Lemma \ref{lemma: squarefree discriminant implies special fibre in D1}, $b_k \in D^1(k)$.
Let $v\in V_b(k)$ be a semisimple element. 
Then the centraliser $Z_H(v)$ is a reductive group of semisimple rank one (Part \ref{item: characterisation nodal curves D1} of \S\ref{subsection: spreading out}).

We claim that there exists a lift $\tilde{v} \in V(R)$ of $v$ such that $\tilde{v}_K \in V(K)\subset \lieh(K)$ is semisimple and such that the group scheme $Z_H(\tilde{v})\rightarrow \Spec R$ is smooth with connected fibres.
Indeed, let $\mathfrak{c}\subset V_k$ be a Cartan subspace containing $v$ (here we use the extension of Vinberg theory to positive characteristic of \cite{Levy-Vinbergtheoryposchar}). 
Let $x \in V(R)$ be a lift of some regular semisimple element in $\mathfrak{c}$. Then $x_K$ is regular semisimple (this being an open condition) and its centraliser $\tilde{\mathfrak{c}} \coloneqq \mathfrak{z}_{\lieh}(x)\subset V_R$ is a Cartan subspace lifting $\mathfrak{c}$. 
Since $k$ is separably closed, $\mathfrak{c}$ is a split Cartan subalgebra; since $R$ is complete the same is true for $\tilde{\mathfrak{c}}$.
We may therefore choose an element $\tilde{v} \in \tilde{\mathfrak{c}}$ lifting $v$ that vanishes on the same roots of $\mathfrak{c}$ as $v$; this $\tilde{v}$ will satisfy the desired properties.
The smoothness of the centraliser $L \coloneqq Z_H(\tilde{v}) \rightarrow \Spec R$ follows from \cite[Theorem 1.1(1)]{Cotner-centralizersofreductivegroupschemes}.
The connectedness of the fibres follows from Lemma \ref{lemma: centraliser semisimple element}, whose proof continues to hold if the characteristic of $k$ is not a torsion prime for $\lieh$, which is weaker than our assumption that the order of the Weyl group is invertible in $k$ (Part 1 of \S\ref{subsection: spreading out}).

The involution $\theta\colon \lieh \rightarrow \lieh$ restricts to a stable involution of the Lie algebra $\mathfrak{l}$ of $L$ by \cite[Lemma 2.5]{Thorne-thesis}. 
We claim that the morphism of $R$-schemes $G \times \mathfrak{l}^{\theta = -1} \rightarrow V_R, (g,x) \mapsto g\cdot x$ is smooth.
Since the domain and target are $R$-flat, it suffices to check this $R$-fibrewise \cite[(I.7.4)]{DeligneRapoport-Schemasdemodulesdecourbes}.
This then follows from \cite[Proposition 4.5]{Thorne-thesis} (noting that $X^1 = \mathfrak{l}^{\theta=-1}$ in this case), whose proof continues to hold when the characteristic of $k$ does not divide the order of the Weyl group of $H$.

Let $\mathfrak{l} \xrightarrow{\pi_L} B_L \coloneqq \mathfrak{l} \GIT L$ be the GIT quotient and let $\phi\colon B_L \rightarrow B$ the map induced by the inclusion $\mathfrak{l} \subset \lieh$.
Since $\phi$ is \'etale at $\pi_L(v)\in B_L(k)$ (Lemma \ref{lemma: etale local structure GIT quotient}), the $R$-point $b\in B(R)$ uniquely lifts to an $R$-point $b_L \in B_L(R)$ satisfying $b_{L,k} = \pi_L(v)$.
Since $b_L$ is open in the fibre $\phi^{-1}(b)$, $\mathfrak{l}_{b_L} \coloneqq \pi_L^{-1}(b_L)$ is an open subscheme of $\mathfrak{l} \cap \lieh_b$.
Using the previous paragraph, this implies that the action map $m\colon G \times \mathfrak{l}^{\theta = -1}_{b_L} \rightarrow V_b$ is smooth. 
We claim that $m$ is also surjective. 
By Part \ref{item: image action map kostant section} of \S\ref{subsection: spreading out}, the image of $m$ contains the set $V_b^{\reg}$ of regular elements (in the sense of Lie theory).
The complement $V_b\setminus V_b^{\reg}$ consists of the semisimple elements of the special fibre $V_{b,k}$. 
Since all such semisimple elements are $G(\bar{k})$-conjugate and since the image of $m$ contains $v$, we conclude that $m$ is surjective.

Since regularity is a smooth-local property \cite[Tag \href{https://stacks.math.columbia.edu/tag/036D}{036D}]{stacksproject}, the smooth surjective morphism $m$ shows that it suffices to prove that $\mathfrak{l}^{\theta = -1}_{b_L}$ is a regular scheme. 
We now make $\mathfrak{l}$ more explicit.
Since $2$ is invertible in $R$ and since $L$ is reductive of semisimple rank $1$, $\mathfrak{l} = Z(\mathfrak{l}) \oplus \mathfrak{l}^{der}$, where $Z(\mathfrak{l})$ and $\mathfrak{l}^{der}$ are the centre and the derived subalgebra of $\mathfrak{l}$ respectively. Since $k$ is separably closed, $\mathfrak{l}^{der} \simeq \liesl_{2,R}$. 
We claim that any two stable involutions on $\mathfrak{l}$ are \'etale locally conjugate.
Indeed, the subscheme of elements of $L$ mapping one such involution to another is smooth, so to prove that it has sections \'etale locally we merely have to show it surjects on $\Spec R$, which follows from a spreading out of Lemma \ref{lemma: stable involutions over alg closed field are conjugate}. 
(See \cite[Proposition 5.6]{Laga-F4paper} for the proof of a similar statement.)
Therefore we may assume that in the decomposition $\mathfrak{l} \simeq Z(\mathfrak{l}) \oplus \liesl_{2,R}$, $\theta$ corresponds to the standard stable involution $\Ad((1,-1))$ of $\liesl_{2,R}$ and to $-1$ on $Z(\mathfrak{l})$.
Moreover if $\Delta_L$ denotes the discriminant polynomial of $\mathfrak{l}$ then $\Delta_L(b_L)$ equals $\Delta(b)$ up to a unit in $R$ (Lemma \ref{lemma: etale local structure GIT quotient}).
We may now calculate that $\mathfrak{l}^{\theta= -1}_{b_L}$ is isomorphic to the scheme $(xy = \Delta_L(b_L))$. 
This scheme is regular since $\Delta(b)$ and hence $\Delta_L(b_L)$ is a uniformiser of $R$.
\end{proof}

\begin{lemma}\label{lemma: squarefree discriminant implies curve is regular}
Let $b\in B(R)$ with $\ord_K(\Delta(b)) = 1$. Then $C_b$ is regular and its geometric special fibre is integral and has a unique singularity, which is a node. 
Moreover, the group scheme $J^1_b\rightarrow \Spec R$ (where $J^1$ is introduced in \S\ref{subsection: the family of curves over B 1}) is the N\'eron model of its generic fibre.
\end{lemma}
\begin{proof}
By Lemma \ref{lemma: squarefree discriminant implies Vb is regular} the scheme $V_b$ is regular. 
Moreover since $C^{\circ}$ is (the spreading out of) a transverse slice of the $G$-action on $V$, the map $G\times C_b^{\circ} \rightarrow V_b$ is smooth (Part \ref{item: action map subregular slice is smooth} of \S\ref{subsection: spreading out}). 
Since regularity is smooth-local, it follows that $C_b^{\circ}$ is a regular scheme. Since $C_b$ is smooth in a Zariski neighbourhood of the marked points $\infty_{i,b}$, the scheme $C_b$ is regular too.
We have $b_k\in D^1(k)$ by Lemma \ref{lemma: squarefree discriminant implies special fibre in D1}.
Therefore the special fibre $C_{b_k}$ is geometrically integral and has a unique nodal singularity (Part \ref{item: characterisation nodal curves D1} of \S\ref{subsection: spreading out}).
The claim about $J^1_b$ follows from the regularity of $C_b$ and a result of Raynaud \cite[\S9.5, Theorem 1]{BLR-NeronModels}.
\end{proof}

The next theorem is not necessary for the proof of Theorem \ref{theorem: integral representatives exist}, but completely determines the integral orbits in the case of square-free discriminant and will be useful in \S\ref{subsection: estimates reducibility stabilizers} and Chapter \ref{section: the average size of the 2-Selmer group}.

\begin{theorem}\label{theorem: summary orbits of square-free discriminant}
	Let $R$ be a discrete valuation ring in which $N$ is a unit. Let $K = \Frac R$ and let $\ord_K: K^{\times} \twoheadrightarrow \Z$ be the normalised discrete valuation. 
	Let $b\in B(R)$ and suppose that $\ord_K \Delta(b)\leq 1$.
	Then:
	\begin{enumerate}
		\item If $x\in V_b(R)$, then $Z_{G}(x)(K) = Z_{G}(x)(R)$. 
		\item The natural map $\alpha\colon G(R)\backslash V_b(R) \rightarrow G(K)\backslash V_b(K)$ is injective and its image contains $\eta_b\left(J_b(K)/2J_b(K)\right)$. 
		\item If furthermore $R$ is complete and has finite residue field then the image of $\alpha$ equals $\eta_b\left(J_b(K)/2J_b(K)\right)$. 
	\end{enumerate}
\end{theorem}

\begin{proof}

If $\ord_K \Delta(b)=0$, $J_b$ is smooth and proper over $R$.
Since $Z_{G}(x)$ is finite \'etale over $R$, the first part follows.
By Proposition \ref{proposition: spread out orbit parametrization galois} and Lemma \ref{lemma: injective H^1 for quasifinite etale gp scheme} below, $\alpha$ is injective.
Proposition \ref{proposition: inject 2-descent orbits spreading out} and the equality $J_b(K)=J_b(R)$ implies that $\eta_b\colon J_b(K)/2J_b(K)\rightarrow \bigG(K)\backslash \bigV_b(K)$ factors through $\bigG(R)\backslash \bigV_b(R)$, so the second part follows.
If $R$ is complete and the residue field $k$ is finite, the pointed sets $\HH^1(R, \bigG)$ and $\HH^1(R,J_b)$ are trivial by \cite[III.3.11(a)]{milne-etalecohomology} and Lang's theorem.
The third part then follows from the fact that the $2$-descent map $ J_b(R)/2J_b(R)\rightarrow \HH^1(R,J_b[2])$ is an isomorphism.

We now assume that $\ord_K \Delta(b)=1$.
Then $b_k \in D^1(k)$ by Lemma \ref{lemma: squarefree discriminant implies special fibre in D1}.
By Lemma \ref{lemma: squarefree discriminant implies curve is regular}, $C_b/R$ is regular, has geometrically integral fibres and its special fibre has a unique nodal singularity.
By the same lemma the group scheme $J_b^1/R$ introduced in \S\ref{subsection: the family of curves over B 1} is the N\'eron model of its generic fibre.
Moreover we have an isomorphism $Z^1_b \simeq J^1_b[2]$ of quasi-finite \'etale group schemes over $R$ by Theorem \ref{theorem: isomorphism centraliser 2-torsion over B1} (or rather its spreading out, Part \ref{item: isomorphism centraliser and two torsion} of \S\ref{subsection: spreading out}).

By Lemmas \ref{lemma: squarefree discriminant implies Vb is regular} and \ref{lemma: smooth locus of invariant map} (and the spreading out of the latter, Part \ref{item: smooth locus of V-> B} of \S\ref{subsection: spreading out}) the scheme $V_b$ is regular and the smooth locus of the morphism $V_b \rightarrow \Spec R$ coincides with the locus $V_b^{\reg}$ of regular elements of $V_b$ (this time in the sense of Lie theory).
Since a section of a morphism between regular schemes lands in the smooth locus \cite[\S3.1, Proposition 2]{BLR-NeronModels}, we see that $V_b(R) = V_b^{\reg}(R)$.
By Part \ref{item: image action map kostant section} of \S\ref{subsection: spreading out}, the morphism $G \times \Spec R\rightarrow V_b^{\reg}, (g,b) \mapsto g \cdot \kappa_b$ is a torsor under the group scheme $Z^1_b$ from \S\ref{subsection: representation theory over B1}.
By Lemma \ref{lemma: AIT full generality} we obtain a bijection of pointed sets
\begin{align}\label{equation: orbit parametrisation square-free discriminant}
G(R) \backslash V_b(R) = G(R) \backslash V_b^{\reg}(R) \simeq \ker(\HH^1(R,J_b^1[2])\rightarrow \HH^1(R,G)).
\end{align}
We now prove the first part of the theorem. Since $x \in V_b^{\reg}(R)$ is \'etale locally $G$-conjugate to $\kappa_b$ by the previous paragraph, we may assume that $x = \kappa_b$.
But then $Z_G(\kappa_b) = Z_b^1 \simeq J_b^1[2]$ and $J_b^1$ satisfies the N\'eron mapping property, so $J_b^1[2](R) = J_b^1[2](K)$.

To prove the remaining parts, note that the map $\HH^1(R, J_b^1[2])\rightarrow \HH^1(K, J_b^1[2]))$ is injective (Lemma \ref{lemma: injective H^1 for quasifinite etale gp scheme} below), so by \eqref{equation: orbit parametrisation square-free discriminant} the map $\bigG(R)\backslash \bigV_b(R) \rightarrow \bigG(K)\backslash \bigV_b(K)$ is injective too.
To show that the image of $\bigG(R)\backslash \bigV_b(R)  \rightarrow \bigG(K)\backslash \bigV_b(K)$ contains $\eta_b\left(J_b(K)/2J_b(K)\right)$, note that we have an exact sequence of smooth group schemes
\begin{align*}
0 \rightarrow J^1_b[2]\rightarrow J^1_b \xrightarrow{\times 2} J^1_b \rightarrow 0,
\end{align*}
since $J^1_b$ has connected fibres. This implies the existence of a commutative diagram:
\begin{center}
	\begin{tikzcd}
		 J^1_b(R)/2J^1_b(R) \arrow[d] \arrow[r , "="] & J_b(K)/2J_b(K) \arrow[d] \\
		{\HH^1(R,J^1_b[2])} \arrow[r]      & {\HH^1(K,J_b[2])}                        
	\end{tikzcd}
	\end{center}
It therefore suffices to prove that every element in the image of the map $J^1_b(R)/2J^1_b(R) \rightarrow  \HH^1(R,J^1_b[2])$ has trivial image in $\HH^1(R,\bigG)$. 
This is true, since the pointed kernel of the map $\HH^1(R,\bigG) \rightarrow \HH^1(K,\bigG)$ is trivial (Proposition \ref{proposition: known cases Grothendieck Serre conjecture}).

If $R$ has finite residue field then \cite[III.3.11(a)]{milne-etalecohomology} and Lang's theorem imply that $\HH^1(R,\bigG) = \{1\}$. 
In this case the $\bigG(R)$-orbits of $\bigV_b(R)$ are in bijection with $\HH^1(R,J^1_b[2])$ by \eqref{equation: orbit parametrisation square-free discriminant}.
The triviality of $\HH^1(R,J^1_b)$ (again by Lang's theorem) shows that $\HH^1(R,J^1_b[2])$ is in bijection with $J^1_b(R)/2J^1_b(R) = J_b(K)/2J_b(K)$. 
This proves Part 3, completing the proof of the proposition.  
\end{proof}

\begin{lemma}\label{lemma: injective H^1 for quasifinite etale gp scheme}
    Let $\Gamma$ be a quasi-finite \'etale commutative group scheme over $\Spec R$. Suppose that $\Gamma$ is a N\'eron model of its generic fibre: for every \'etale extension $R \rightarrow R'$ of discrete valuation rings, we have $\Gamma(R') = \Gamma(\Frac R')$. 
	Then the map of \'etale cohomology groups $\HH^1(R,\Gamma) \rightarrow \HH^1(K,\Gamma)$ is injective.
\end{lemma}
\begin{proof}
	Let $j\colon \Spec K \rightarrow \Spec R$ denote the natural inclusion. Then the N\'eron mapping property translates into the equality of \'etale sheaves $j_*j^*\Gamma = \Gamma$. The map $\HH^1(R,\Gamma) \rightarrow \HH^1(K,\Gamma)$ is therefore injective because it is the first term in the five-term exact sequence associated with the Leray spectral sequence $\HH^p(R,\mathrm{R}^qj_*j^*\Gamma) \Rightarrow \HH^{p+q}(K,\Gamma)$.
\end{proof}

\subsection{Proof of Theorem \ref{theorem: integral representatives exist}}\label{subsection: proof of integral representatives general case}

In this section we use the results from \S\ref{subsection: some stacks} and \S\ref{subsection: square-free discriminant case} to complete the proof of Theorem \ref{theorem: integral representatives exist}. 

We start by summarizing the broad strategy. 
Let $p$ be a prime not dividing $N$, $b\in B(\Z_p) \cap B^{\rs}(\Q_p)$, $P \in J_b(\Q_p)$ and $\eta_b(P) \in G(\Q_p) \backslash V_b(\Q_p)$ the orbit constructed in Theorem \ref{theorem: inject 2-descent orbits}.
Let $\lambda \colon \Spec \Q_p \rightarrow \mathcal{M}_b = [G\backslash V_b]$ be the $\Q_p$-point of $\mathcal{M}_b$ corresponding to $\eta_b(P)$ under Lemma \ref{lemma: points on quotient stack with trivial torsor are just orbits}.
We wish to show that $\eta_b(P)$ has a representative in $V_b(\Z_p)$. By Example \ref{example: Zp-points quotient stack}, this is equivalent to showing that $\lambda$ extends to a morphism $\tilde{\lambda} \colon \Spec \Z_p \rightarrow \mathcal{M}_b$.
Instead of constructing such an $\tilde{\lambda}$ directly, we will build a $2$-dimensional scheme around $\Spec \Z_p$, construct a map to $\mathcal{M}$ on a large open subscheme of this scheme and use Lemma \ref{lemma: purity for the stack M} to extend this map to the whole scheme, which can then be specialized to $\Spec \Z_p$ giving the desired $\tilde{\lambda}$.
More precisely:
\begin{enumerate}
    \item Using Bertini theorems and properties of the compactified Jacobian, we construct a $2$-dimensional regular integral scheme $\mathcal{X}$ together with a morphism $x\colon \Spec \Z_p \rightarrow \mathcal{X}$ and a morphism $\tilde{b}\colon \mathcal{X} \rightarrow B$ such that $\tilde{b}$ extends $b$ in the sense that $\tilde{b}\circ x = b$, and such that $\Delta(\tilde{b})$ is square-free in some sense. (See Corollary \ref{corollary: bertini theorem for compactified jacobian} for a precise statement.)
    \item Using the results of \S\ref{subsection: square-free discriminant case} (more specifically Lemma \ref{lemma: squarefree discriminant implies curve is regular}), we find an open subset $U\subset \mathcal{X}$ such that $x_{\Q_p}\in \mathcal{X}(\Q_p)$ lies in $U(\Q_p)$, $\mathcal{X} \setminus U$ is a finite set of closed points and such that $\lambda\colon \Spec \Q_p\rightarrow \mathcal{M}_b$ extends to a morphism $U\rightarrow \mathcal{M}_{\tilde{b}}$.
    \item By the purity Lemma \ref{lemma: purity for the stack M}, the latter morphism extends to a morphism $\mathcal{X} \rightarrow \mathcal{M}_{\tilde{b}}$. Precomposing with $x\colon \Spec \Z_p\rightarrow \mathcal{X}$ gives a morphism $\tilde{\lambda}\colon \Spec \Z_p\rightarrow \mathcal{M}_b$ extending $\lambda$, as desired.
\end{enumerate}

Before we carry out this proof strategy precisely, we need to state the following Bertini type theorem over $\Z_p$, proved in \cite[Proposition 4.22]{Laga-E6paper}.

\begin{proposition}\label{proposition: Bertini type theorem}
	Let $p$ be a prime number. Let $\mathcal{Y} \rightarrow \Z_p$ be a smooth, quasiprojective morphism of relative dimension $d\geq 1$ with geometrically integral fibres. 
	Let $\mathcal{D} \subset \mathcal{Y}$ be an effective Cartier divisor. Assume that $\mathcal{Y}_{\F_p}$ is not contained in $\mathcal{D}$ (i.e. $\mathcal{D}$ is horizontal) and that $\mathcal{D}_{\Q_p}$ is reduced.
	Let $P\in \mathcal{Y}(\Z_p)$ be a section such that $P_{\Q_p} \not\in \mathcal{D}_{\Q_p}$. 
	Then there exists a closed subscheme $\mathcal{X} \hookrightarrow \mathcal{Y}$ containing the image of $P$ satisfying the following properties.
	\begin{itemize}
		\item $\mathcal{X} \rightarrow \Z_p$ is smooth of relative dimension $1$ with geometrically integral fibres.
		\item $\mathcal{X}_{\F_p}$ is not contained in $\mathcal{D}$ and the (scheme-theoretic) intersection $\mathcal{X}_{\Q_p} \cap \mathcal{D}_{\Q_p}$ is reduced.
	\end{itemize} 
\end{proposition}

Recall that $\bar{J}$ denotes the compactified Jacobian introduced in \S\ref{section: the compactified jacobian}, which has been spread out in \S\ref{subsection: spreading out} to a scheme over $\Z[1/N]$.

\begin{corollary}\label{corollary: bertini theorem for compactified jacobian}
Let $p$ be a prime not dividing $N$.
Let $b\in B(\Z_p) \cap B^{\rs}(\Q_p)$ and $P \in J_b(\Q_p)$. 
Then there exists a morphism $\mathcal{X}\rightarrow\Z_p$ which is of finite type, smooth of relative dimension $1$ and has geometrically integral fibres, together with a morphism $\mathcal{X} \rightarrow \bar{J}_{\Z_p}$ satisfying the following properties. 
\begin{enumerate}
\item Let $\tilde{b}$ be the composition $ \mathcal{X} \rightarrow \bar{J}_{\Z_p}\rightarrow B_{\Z_p}$. Then the discriminant of $\tilde{b}$, seen as a map $\mathcal{X} \rightarrow \A^1_{\Z_p}$, is square-free on the generic fibre of $\mathcal{X}$ and not identically zero on the special fibre.
\item There exists a section $x\in \mathcal{X}(\Z_p)$ such that the composition $\Spec \Q_p \xrightarrow{x_{\Q_p}} \mathcal{X} \rightarrow \bar{J}_{\Z_p}$ coincides with $P$.
\end{enumerate}
\end{corollary}
\begin{proof}
	We apply Proposition \ref{proposition: Bertini type theorem} with $\mathcal{Y} = \CJac_{\Z_p}$.
	We define $\mathcal{D}$ to be the pullback of the discriminant locus $\{ \Delta = 0 \} \subset \intbigB_{\Z_p}$ under the morphism $\CJac_{\Z_p} \rightarrow \intbigB_{\Z_p}$. Since the latter morphism is proper, we can extend $P \in \Jac_b(\Q_p)$ to an element of $\CJac_b(\Z_p)$, still denoted by $P$. 
	We claim that the triple $(\mathcal{Y}, \mathcal{D},P)$ satisfies the assumptions of Proposition \ref{proposition: Bertini type theorem}. 
	This follows from an argument identical to \cite[Corollary 4.23]{Laga-E6paper}, using Part \ref{item: compactified jacobian spreading out} of \S\ref{subsection: spreading out}.
\end{proof}

\begin{proof}[Proof of Theorem \ref{theorem: integral representatives exist}]
Choose a relative curve $\mathcal{X} \rightarrow \Z_p$, a map $\mathcal{X} \rightarrow \bar{J}_{\Z_p}$ and a section $x\in \mathcal{X}(\Z_p)$ satisfying the conclusions of Corollary \ref{corollary: bertini theorem for compactified jacobian}, and let $\tilde{b}$ be the composition $\mathcal{X} \rightarrow \bar{J}_{\Z_p} \rightarrow B_{\Z_p}$.
Recall that $J^1$ is an open subscheme of $\bar{J}$; let $\mathcal{X}^1$ denote the open subscheme of $\mathcal{X}$ landing in $J^1_{\Z_p}$.

We claim that the complement of $\mathcal{X}^1$ in $\mathcal{X}$ is a union of finitely many closed points.
Indeed, by Lemma \ref{lemma: squarefree discriminant implies curve is regular} and the fact that the discriminant of $\tilde{b}_{\Q_p}$ is square-free, the group scheme $J^1_{\tilde{b}_{\Q_p}}\rightarrow \mathcal{X}_{\Q_p}$ is a N\'eron model of its generic fibre.
By the N\'eron mapping property, the section $\mathcal{X}_{\Q_p} \rightarrow \bar{J}^1_{\Q_p}$ must land in $J^1_{\Q_p}$.
Since the discriminant of $\mathcal{X}$ is nonzero on the special fibre, it follows that $\mathcal{X}_{\F_p}^1$ is nonempty.
Combining the last two sentences and the fact that $\mathcal{X}_{\F_p}$ is irreducible proves the claim.

To carry out the second step in the proof strategy sketched in the beginning of this section, we construct a morphism $\mathcal{X}^1 \rightarrow \mathcal{M}_{\tilde{b}}$ such that precomposing this morphism with $x_{\Q_p}\colon \Spec \Q_p \rightarrow \mathcal{X}$ equals the morphism $\lambda \colon \Spec \Q_p \rightarrow \mathcal{M}_b$ that corresponds to the orbit $\eta_b(P)$ under Lemma \ref{lemma: points on quotient stack with trivial torsor are just orbits}.
The multiplication-by-$2$ map $J^1 \xrightarrow{\times 2} J^1$ on the semiabelian scheme $J^1\rightarrow B^1$ is a $J^1[2]$-torsor, and pulling back this torsor along the morphism $\mathcal{X}^1\rightarrow J^1$ defines a $J^1[2]$-torsor $T\rightarrow \mathcal{X}^1$. 
Using the isomorphism $J^1[2]\simeq Z^1$ of Theorem \ref{theorem: isomorphism centraliser 2-torsion over B1} (and its spreading out version of Property \ref{item: isomorphism centraliser and two torsion} in \S\ref{subsection: spreading out}), we obtain a $Z^1$-torsor $T'\rightarrow \mathcal{X}^1$.
In the language of \S\ref{subsection: some stacks}, this torsor determines a morphism $\mathcal{X}^1\rightarrow [B/Z]$.
Composing this morphism with the inclusion $[B/Z]\hookrightarrow \mathcal{M}$ (given by twisting the Kostant section and described in \S\ref{subsection: some stacks}), we obtain a morphism $\mathcal{X}^1\rightarrow \mathcal{M}_{\tilde{b}}$.
Specializing this morphism at $x_{\Q_p}$ corresponds to the orbit $\eta_b(P)$ under Lemma \ref{lemma: points on quotient stack with trivial torsor are just orbits}, by the explicit construction of $\eta_b$ in Theorem \ref{theorem: inject 2-descent orbits}.

Finally by Lemma \ref{lemma: purity for the stack M} --- and this is the key point --- the above morphism $\mathcal{X}^1\rightarrow \mathcal{M}_{\tilde{b}}$ extends (uniquely) to a morphism $\mathcal{X} \rightarrow \mathcal{M}_{\tilde{b}}$.
Precomposing with $x\colon \Spec \Z_p\rightarrow  \mathcal{X}$ defines a morphism $\tilde{\lambda}\colon \Spec \Z_p\rightarrow \mathcal{M}_b$ whose $\Q_p$-fibre corresponds to $\eta_b(P)$ under Lemma \ref{lemma: points on quotient stack with trivial torsor are just orbits}. 
By Example \ref{example: Zp-points quotient stack}, $\tilde{\lambda}$ actually arises from an element of $V_b(\Z_p)$.
We conclude that $\eta_b(P)$ has a representative in $V_b(\Z_p)$, completing the proof.
\end{proof}

\subsection{Orbits over \texorpdfstring{$\Z$}{Z}}\label{subsection: orbits over Z}

Recall that $\sh{E}_p = \intbigB(\Z_p) \cap \bigB^{\rs}(\Q_p)$ for all $p$.
Define $\sh{E} \coloneqq \intbigB(\Z) \cap \bigB^{\rs}(\Q)$.
We state the following corollary, whose proof is completely analogous to the proof of \cite[Corollary 5.8]{Thorne-Romano-E8} and uses the fact that $\intbigG$ has class number $1$ (Proposition \ref{proposition: G has class number one}).
\begin{corollary}\label{corollary: weak global integral representatives}
	Let $b_0 \in \sh{E}$. Then for each prime $p$ dividing $N$ we can find an open compact neighbourhood $W_p$ of $b_0$ in $\sh{E}_p$ and an integer $n_p\geq 0$ with the following property. Let $M = \prod_{p | N} p^{n_p}$. Then for all $b\in \sh{E} \cap \left(\prod_{p| N} W_p \right)$ and for all $y \in \Sel_2(\Jac_{M\cdot b})$, the orbit $\eta_{M\cdot b}(y) \in \bigG(\Q) \backslash \bigV_{M\cdot b}(\Q)$ contains an element of $\intbigV_{M\cdot b}(\Z)$. 
\end{corollary} 

This statement about integral representatives will be strong enough to obtain the main theorems in \S\ref{section: geometry-of-numbers}.

\section{Geometry-of-numbers}\label{section: geometry-of-numbers}


In this chapter we will apply the counting techniques of Bhargava to provide estimates for the integral orbits of bounded height in the representation $(\intbigG,\intbigV)$. 
We keep the notation from the previous chapters and continue to assume that $H$ is not of type $A_1$.

\subsection{Heights} \label{subsection: heights}

Recall that $\intbigB = \Spec \Z[p_{d_1},\dots,p_{d_r}]$ and that $\pi: \intbigV \rightarrow \intbigB$ denotes the morphism of taking invariants.
For any $b\in \bigB(\Real)$ we define the \define{height} of $b$ by the formula
$$\height(b) \coloneqq \sup_{1\leq i\leq r} |p_i(b)|^{1/i}.$$
We define $\height(v) = \height(\bigpi(v))$ for any $v\in \bigV(\Real)$.
We have $\height(\lambda\cdot b) = |\lambda|\height(b)$ for all $\lambda\in \Real$ and $b\in \bigB(\Real)$. 
If $A$ is a subset of $\bigV(\Real)$ or $\bigB(\Real)$ and $X\in \Real_{>0}$ we write $A_{<X}\subset A$ for the subset of elements of height $<X$. 
For every such $X$, the set $\intbigB(\Z)_{<X}$ is finite.

The next lemma records a numerological fact, which implies that $\intbigB(\Z)_{<X}$ has order of magnitude $X^{\dim V}$.

\begin{lemma}\label{lemma: sum of invariants equals the dimension of V}
We have $d_1+ \dots +d_r = \dim_{\Q}V$.
\end{lemma}
\begin{proof}
Recall from \S\ref{subsection: a stable Z/2Z-grading} that $\Phi_H$ denotes a root system of $H$. 
We prove the two equalities $d_1+ \dots +d_r  = \frac{1}{2}\#\Phi_H+\rank H = \dim V$.
The first one is classical, see \cite[Corollary 10.2.4]{Carter-SimpleGroupsLieType1972}; the second one follows from \cite[Lemma 2.21]{Thorne-thesis} applied to $x = 0$.
\end{proof}

\subsection{Measures on \texorpdfstring{$G$}{G}}\label{subsection: measures on G}

Let $\omega_{\bigG}$ be a generator for the $\Q$-vector space of left-invariant top differential forms on $\intbigG$ over $\Q$. 
It is well-defined up to an element of $\Q^{\times}$ and it determines Haar measures $dg$ on $\bigG(\Real)$ and $\bigG(\Q_p)$ for each prime $p$.

Recall from \S\ref{subsection: A family of curves} that $m$ denotes the number of marked points of the family of curves $C\rightarrow B$.

\begin{proposition}\label{proposition: tamagawa number}
	The product $\vol\left(\intbigG(\Z)\backslash \intbigG(\Real) \right) \cdot \prod_p \vol\left(\intbigG(\Z_p)\right)$ converges absolutely and equals $2^{m}$, the Tamagawa number of $\bigG$.
\end{proposition}
\begin{proof}
	Proposition \ref{proposition: G has class number one} implies that the product equals the Tamagawa number $\tau(\bigG)$ of $G$. 
	By Proposition \ref{proposition: simply connected cover of $G$}, the group $G$ is semisimple and its fundamental group has order $2 \#\pi_0(H^{\theta})$; let $G_{sc} \rightarrow G$ be its simply connected cover.
	The proof of Proposition \ref{proposition: orbits corresponding to marked points are reducible} (more precisely the isomorphism between \eqref{equation: sequence 2-torsion open curve} and \eqref{equation: sequence centralisers of kostant section in H^theta}) shows that $\#\pi_0(H^{\theta})$ has order $2^{m-1}$.
	Now use the identities $\tau(G)=2^m\tau(G_{sc})$ \cite[Theorem 2.1.1]{Ono-relativetheorytamagawa} and $\tau(G_{sc})=1$ \cite{Langlands-volumefunddomainchevalley}.
\end{proof}

We study the measure $dg$ on $\bigG(\Real)$ using the Iwasawa decomposition, after introducing some notation.
Recall from \S\ref{subsection: a stable Z/2Z-grading} that we have fixed a maximal torus $T\subset H$ with set of roots $\Phi_H$.
Moreover we have fixed a Borel subgroup $P$ containing $T$, which determines a root basis $S_H \subset \Phi_H$ and a set of positive roots $\Phi_H^+$.
Then $T^{\theta}$ is a maximal torus of $G$ and $P^{\theta}$ a Borel subgroup of $G$ \cite[Lemma 5.1]{Richardson-orbitsinvariantsrepresentationsinvolutions}. 
Let $\Phi_G = \Phi(G,T^{\theta})$ be its set of roots, $S_G = \{b_1,\dots,b_k\}$ the corresponding root basis and $\Phi_G^{\pm}$ the subset of positive/negative roots.
Fix, once and for all, a maximal compact subgroup $K\subset G(\Real)$.
If $N$ is the unipotent radical of $P^{\theta}$ we have a decomposition $\bigP^{\theta} = \bigT^{\theta} \bigN \subset \bigG$.
Let $\bigPopp = \bigT \bigNopp \subset \bigG$ be the opposite Borel subgroup. 
Then the natural product maps 
\begin{equation*}
	\bigNopp(\Real)\times T^{\theta}(\Real)^{\circ}\times K \rightarrow \bigG(\Real) ,\; T^{\theta}(\Real)^{\circ} \times \bigNopp(\Real) \times K \rightarrow \bigG(\Real)
\end{equation*} 
are diffeomorphisms. 
If $t\in T^{\theta}(\Real)$, let $\delta_{\bigG}(t) = \prod_{\beta\in \Phi_{\bigG}^-} \beta(t) = \det \Ad(t)|_{\Lie \bigNopp(\Real)}$. 
The following result follows from well-known properties of the Iwasawa decomposition; see \cite[Chapter 3;\S1]{Lang-SL2R}.

\begin{lemma}\label{lemma: Haar measure iwasawa decomposition}
	Let $dt, dn, dk$ be Haar measures on $T^{\theta}(\Real)^{\circ}, \bigNopp(\Real), K$ respectively. Then the assignment 
	\begin{align*}
		f \mapsto \int_{t\in T^{\theta}(\Real)^{\circ} } \int_{n\in \bigNopp(\Real)} \int_{k\in K } f(tnk)\, dk\, dn\, dt = \int_{t\in T^{\theta}(\Real)^{\circ} } \int_{n\in \bigNopp(\Real)} \int_{k\in K } f(ntk)\delta_{\bigG}(t)^{-1} \, dk\, dn\, dt
	\end{align*}
	defines a Haar measure on $\bigG(\Real)$.
\end{lemma}

We now fix Haar measures on the groups $T^{\theta}(\Real)^{\circ}, K$ and $\bigNopp(\Real)$, as follows. 
We give $T^{\theta}(\Real)^{\circ}$ the measure pulled back from the isomorphism $\prod_{\beta\in S_{\bigG}} \beta \colon T^{\theta}(\Real)^{\circ} \rightarrow \Real_{>0}^{\#S_G}$, where $\Real_{>0}$ gets its standard Haar measure $d^{\times} \lambda = d\lambda/\lambda$.
We give $K$ the probability Haar measure. 
Finally we give $\bigNopp(\Real)$ the unique Haar measure $dn$ such that the Haar measure on $\bigG(\Real)$ from Lemma \ref{lemma: Haar measure iwasawa decomposition} coincides with $dg$.

\subsection{Measures on \texorpdfstring{$V$}{V}}\label{subsection: measures on V}

Let $\omega_{\bigV}$ be a generator of the free rank one $\Z$-module of left-invariant top differential forms on $\intbigV$.
Then $\omega_{\bigV}$ is uniquely determined up to sign and it determines Haar measures $dv$ on $\bigV(\Real)$ and $\bigV(\Q_p)$ for every prime $p$.
We define the top form $\omega_{\bigB} = dp_{d_1}\wedge\cdots \wedge dp_{d_r}$ on $\intbigB$.
It defines measures $db$ on $\bigB(\Real)$ and $\bigB(\Q_p)$ for every prime $p$. 

\begin{lemma}\label{lemma: relations different forms on V,G,B}
There exists a unique rational number $W_0\in \Q^{\times}$ with the following property.
Let $k/\Q$ be a field extension, let $\mathfrak{c}$ a Cartan subalgebra of $\bigh_k$ contained in $\bigV_k$, and let $\mu_{\mathfrak{c}}\colon \bigG_k \times \mathfrak{c} \rightarrow \bigV_k$ be the action map. 
Then $\mu^*_{\mathfrak{c}}\omega_{\bigV} = W_0 \omega_{\bigG}\wedge \pi|_{\mathfrak{c}}^*\omega_{\bigB}$.
\end{lemma}
\begin{proof}
    The proof is identical to that of \cite[Proposition 2.13]{Thorne-E6paper}.
    Here we use the fact that the sum of the invariants equals the dimension of the representation: $d_1 + \cdots d_r =  \dim_{\Q} \bigV$ (Lemma \ref{lemma: sum of invariants equals the dimension of V}).
\end{proof}

\begin{lemma}\label{lemma: the constants W0 and W}
	Let $W_0\in \Q^{\times}$ be the constant of Lemma \ref{lemma: relations different forms on V,G,B}. Then:
	\begin{enumerate}
		\item Let $\intbigV(\Z_p)^{\rs} \coloneqq \intbigV(\Z_p)\cap \bigV^{\rs}(\Q_p)$ and define a function $m_p: \intbigV(\Z_p)^{\rs} \rightarrow \Real_{\geq 0}$ by the formula
		\begin{equation}
		m_p(v) = \sum_{v' \in \intbigG(\Z_p)\backslash\left( \bigG(\Q_p)\cdot v\cap \intbigV(\Z_p) \right)} \frac{\#Z_{\intbigG}(v')(\Q_p)  }{\#Z_{\intbigG}(v')(\Z_p) } .
		\end{equation}
		Then $m_p$ is locally constant. 
		\item Let $\intbigB(\Z_p)^{\rs} \coloneqq \intbigB(\Z_p)\cap \bigB^{\rs}(\Q_p)$ and let $\psi_p: \intbigV(\Z_p)^{\rs}  \rightarrow \Real_{\geq 0}$ be a bounded, locally constant function which satisfies $\psi_p(v) = \psi_p(v')$ when $v,v'\in \intbigV(\Z_p)^{\rs}$ are conjugate under the action of $\bigG(\Q_p)$. 
		Then we have the formula 
		\begin{equation}
		\int_{v\in \intbigV(\Z_p)^{\rs}} \psi_p(v) \mathrm{d} v = |W_0|_p \vol\left(\intbigG(\Z_p)\right) \int_{b\in \intbigB(\Z_p)^{\rs}} \sum_{v\in \bigG(\Q_p)\backslash \intbigV_b(\Z_p) } \frac{m_p(v)\psi_p(v)  }{\# Z_{\intbigG }(v)(\Q_p)} \mathrm{d} b .
		\end{equation}
		
	\end{enumerate}
\end{lemma}
\begin{proof}
	The proof is identical to that of \cite[Proposition 3.3]{Romano-Thorne-ArithmeticofsingularitiestypeE}, using Lemma \ref{lemma: relations different forms on V,G,B}. 
\end{proof}

\subsection{Fundamental sets}\label{subsection: fundamental sets}


Let $K\subset \bigG(\Real)$ be the maximal compact subgroup fixed in \S\ref{subsection: measures on G}.
For any $c\in \Real_{>0}$, define $T_c \coloneqq \{t\in T^{\theta}(\Real)^{\circ} \mid \forall \beta\in S_{\bigG} ,\, \beta(t) \leq c   \}$.
A \define{Siegel set} is, by definition, any subset $\Siegel_{\omega, c} \coloneqq \omega \cdot T_c \cdot K$, where $\omega \subset \bigNopp(\Real)$ is a compact subset and $c>0$.

\begin{proposition}\label{proposition: properties Siegel set}
	\begin{enumerate}
		\item For every $\omega\subset \bigNopp(\Real)$ and $c>0$, the set $$\{\gamma\in \intbigG(\Z) \mid \gamma \cdot \Siegel_{\omega,c} \cap \Siegel_{\omega,c}\neq \emptyset \}$$ is finite. 
		\item We can choose $\omega\subset \bigNopp(\Real)$ and $c>0$ such that $\intbigG(\Z) \cdot \Siegel_{\omega,c} = \bigG(\Real)$.
	\end{enumerate}
\end{proposition}
\begin{proof}
	The first part follows from the Siegel property \cite[Corollaire 15.3]{Borel-introductiongroupesarithmetiques}.
	By \cite[Theorem 4.15]{PlatonovRapinchuk-Alggroupsandnumbertheory}, the second part is reduced to proving that $\bigG(\Q)=\bigP(\Q)\cdot \intbigG(\Z)$, which follows from \cite[\S6, Lemma 1(b)]{Borel-densityMaximalityarithmetic}.
\end{proof}

Fix $\omega \subset \bigNopp(\Real)$ and $c>0$ so that $\Siegel_{\omega,c}$ satisfies the conclusions of Proposition \ref{proposition: properties Siegel set}.
By enlarging $\omega$, we may assume that $\Siegel_{\omega,c}$ is semialgebraic. 
We drop the subscripts and for the remainder of \S\ref{section: geometry-of-numbers} we write $\Siegel$ for this fixed Siegel set. 
The set $\Siegel$ will serve as a fundamental domain for the action of $\intbigG(\Z)$ on $\bigG(\Real)$. 

A $\intbigG(\Z)$-coset of $\bigG(\Real)$ may be represented more than once in $\Siegel$, but by keeping track of the multiplicities this will not cause any problems.
The surjective map $\varphi\colon \Siegel \rightarrow \intbigG(\Z)\backslash \bigG(\Real)$ has finite fibres and if $g\in \Siegel$ we define $\mu(g) \coloneqq \# \varphi^{-1}(\varphi(g))$.
The function $\mu \colon \Siegel \rightarrow \mathbb{N}$ is uniformly bounded by $\mu_{\max} \coloneqq \# \{\gamma \in \intbigG(\Z) \mid \gamma \Siegel \cap \Siegel \neq \emptyset \}  $ and has semialgebraic fibres.
By pushing forward measures via $\varphi$, we obtain the formula
\begin{equation}\label{equation: weighted volume Siegel set}
	\int_{g\in \Siegel} \mu(g)^{-1} \,dg = \vol\left(\intbigG(\Z) \backslash \bigG(\Real) \right).
\end{equation}

We now construct special subsets of $\bigV^{\rs}(\Real)$ which serve as our fundamental domains for the action of $\bigG(\Real)$ on $\bigV^{\rs}(\Real)$. 
By the same reasoning as in \cite[\S 2.9]{Thorne-E6paper}, we can find open subsets $L_1,\dots, L_k$ of $\{b\in \bigB^{\rs}(\Real)\mid \height(b)=1 \}$ and sections $s_i \colon L_i \rightarrow \bigV(\Real)$ of the map $\pi\colon \bigV \rightarrow \bigB $ satisfying the following properties:
\begin{itemize}
	\item For each $i$, $L_i$ is connected and semialgebraic and $s_i$ is a semialgebraic map with bounded image. 
	\item Set $\Lambda = \Real_{>0}$. Then we have an equality 
	\begin{equation}
		\bigV^{\rs}(\Real) = \bigcup_{i=1}^k \bigG(\Real) \cdot \Lambda \cdot s_i(L_i). 
	\end{equation}
\end{itemize}

If $v\in s_i(L_i)$ let $r_i = \# Z_{\bigG}(v)(\Real)$; this integer is independent of the choice of $v$. 
We record the following change-of-measure formula, which follows from Lemma \ref{lemma: relations different forms on V,G,B}.

\begin{lemma}\label{lemma: change of variables section}
Let $f \colon \bigV(\Real) \rightarrow \mathbb{C}$ be a continuous function of compact support and $i\in \{1,\dots,k\}$.
Let $G_0\subset \bigG(\Real)$ be a measurable subset and let $m_{\infty}(v)$ be the cardinality of the fibre of the map $G_0 \times\Lambda\times  L_i \rightarrow \bigV(\Real), (g,\lambda,l)\mapsto g\cdot\lambda\cdot  s_i(l)$ above $v\in \bigV(\Real)$. 
Then 
\begin{align*}
\int_{v\in G_0 \cdot \Lambda \cdot s_i(L_i)} f(v)m_{\infty}(v) \,dv = |W_0| \int_{b\in \Lambda\cdot L_i}\int_{g\in G_0} f(g\cdot s_i(b)) \,dg \,db,
\end{align*}
where $W_0\in \Q^{\times}$ is the scalar of Lemma \ref{lemma: relations different forms on V,G,B}.
\end{lemma}

\subsection{Counting integral orbits of \texorpdfstring{$\bigV$}{V}}\label{subsection: counting integral orbits in V}

For any $\intbigG(\Z)$-invariant subset $A\subset \intbigV(\Z)$, define 
$$N(A,X) \coloneqq \sum_{v\in \intbigG(\Z)\backslash A_{<X}} \frac{1}{\# Z_{\intbigG}(v)(\Z)}.$$
(Recall that $A_{<X}$ denotes the set of elements of $A$ of height $<X$.)
Let $k$ be a field of characteristic not dividing $N$.
We say an element $v\in V(k)$ with $b=\bigpi(v)$ is:
\begin{itemize}
	\item \define{$k$-reducible} if $\Delta(b)=0$ or if it is $G(k)$-conjugate to a Kostant section, and \define{$k$-irreducible} otherwise.
	\item \define{$k$-soluble} if $\Delta(b)\neq 0$ and $v$ lies in the image of the map $\eta_b: J_b(k)/2J_b(k) \rightarrow G(k)\backslash V_b(k)$ of Proposition \ref{proposition: inject 2-descent orbits spreading out}.
\end{itemize}

For any $A\subset \intbigV(\Z)$, write $A^{irr}\subset A$ for the subset of $\Q$-irreducible elements. 
Write $\bigV(\Real)^{sol} \subset \bigV(\Real)$ for the subset of $\Real$-soluble elements. 
Write $g$ for the common arithmetic genus of the curves $C \rightarrow B$.

\begin{theorem}\label{theorem: counting R-soluble elements, no congruence}
	We have
	\begin{displaymath}
	N(\intbigV(\Z)^{irr} \cap \bigV(\Real)^{sol},X) = \frac{|W_0|}{2^g}\vol\left(\intbigG(\Z)\backslash \bigG(\Real)\right) \vol\left(\bigB(\Real)_{<X}  \right)+ o\left(X^{\dim V}\right),
	\end{displaymath}
	where $W_0\in \Q^{\times}$ is the scalar of Lemma \ref{lemma: relations different forms on V,G,B}.
\end{theorem}

\begin{remark}
    It is possible to obtain a power saving error term in Theorem \ref{theorem: counting R-soluble elements, no congruence} using a finer analysis and the methods of \cite{ShankarTsimerman-countingpowersaving}, but we have not pursued this here.
\end{remark}

We first explain how to reduce Theorem \ref{theorem: counting R-soluble elements, no congruence} to Proposition \ref{prop: counting sections}.
Recall that there exists $\G_m$-actions on $\bigV$ and $\bigB$ such that the morphism $\bigpi: \bigV \rightarrow \bigB$ is $\G_m$-equivariant and that we write $\Lambda = \Real_{>0}$. 
By an argument identical to \cite[Lemma 5.5]{Laga-E6paper}, the subset $\bigV(\Real)^{sol}\subset \bigV^{\rs}(\Real)$ is open and closed in the Euclidean topology.
Therefore by discarding some of the subsets $L_1,\dots,L_k$ of \S\ref{subsection: fundamental sets}, we may write 
$\bigV(\Real)^{sol} = \bigcup_{i\in J} \bigG(\Real)\cdot \Lambda \cdot s_i(L_i)$
for some $J\subset \{1,\dots,k\}$.
Moreover for every $b\in \bigB^{\rs}(\Real)$ we have equalities $$\#\left(\bigG(\Real)\backslash \bigV_b(\Real)^{sol}\right)/\#Z_{\bigG}(\kappa_b)(\Real)= \#\left(J_b(\Real)/2J_b(\Real)\right)/\#J_b[2](\Real)=1/2^g,    $$
where the first follows from the definition of $\Real$-solubility and Proposition \ref{proposition: isomorphism 2-torsion and Kostant section centraliser thorne}, and the second is a general fact about real abelian varieties.
Therefore by the inclusion-exclusion principle, to prove Theorem \ref{theorem: counting R-soluble elements, no congruence} it suffices to prove the following proposition. (See \cite[\S5.2]{Laga-E6paper} for more details concerning this step.)

For any nonempty subset $I$ of $\{1,\dots,k\}$, write $L_I=\pi\left( \cap_{i\in I} \bigG(\Real)\cdot s_i(L_i) \right)$.
Fix an $i\in I$, write $s_I$ for the restriction of $s_i$ to $L_I$ and write $r_I=r_i$. (The section $s_I$ may depend on $i$ but the number $r_I$ does not if $L_I$ is nonempty.)

\begin{proposition}\label{prop: counting sections}
    In the above notation, let $(L, s, r)$ be $(L_I, s_I, r_I)$ for some $I\subset \{1,\dots, k\}$.
	Then
	$$N((\bigG(\Real)\cdot\Lambda\cdot s(L)) \cap \intbigV(\Z)^{irr},X) = \frac{|W_0|}{r}\vol\left(\intbigG(\Z) \backslash \bigG(\Real) \right)\vol((\Lambda\cdot L)_{<X}) +o\left(X^{\dim V} \right).$$

\end{proposition}

So to prove Theorem \ref{theorem: counting R-soluble elements, no congruence} it remains to prove Proposition \ref{prop: counting sections}. For the latter we will follow the general orbit-counting techniques established by Bhargava, Shankar and Gross \cite{BS-2selmerellcurves, Bhargava-Gross-hyperellcurves} closely.
The only notable differences are that we work with a Siegel set instead of a true fundamental domain and that we have to carry out a case-by-case analysis for cutting off the cusp in \S\ref{subsection: cutting off cusp}. 
For the remainder of \S\ref{section: geometry-of-numbers} we fix a triple $(L,s,r)$ as above with $L\neq \emptyset$.

\subsection{First reductions}

We first reduce Proposition \ref{prop: counting sections} to estimating the number of (weighted) lattice points in a region of $\bigV(\Real)$. 
Recall that $\Siegel$ denotes the Siegel set fixed in \S\ref{subsection: fundamental sets} which comes with a multiplicity function $\mu\colon \Siegel \rightarrow \mathbb{N}$.
Because $\intbigG(\Z)\cdot \Siegel = \bigG(\Real)$, every element of $\bigG(\Real)\cdot \Lambda \cdot s(L)$ is $\intbigG(\Z)$-equivalent to an element of $\Siegel \cdot \Lambda \cdot s(L)$.
In fact, we can be more precise about how often a $\intbigG(\Z)$-orbit will be represented in $\Siegel \cdot \Lambda\cdot s(L)$. Let $\nu \colon \Siegel \cdot \Lambda\cdot s(L) \rightarrow \Real_{>0}$ be the `weight' function defined by 
\begin{equation}\label{equation: weight function nu}
	x\mapsto \nu(x) \coloneqq \sum_{\substack{g\in \Siegel \\ x\in g\cdot \Lambda\cdot s(L) }} \mu(g)^{-1}.
\end{equation}
Then $\nu$ takes only finitely many values and has semialgebraic fibres. 
We now claim that if every element of $\Siegel \cdot \Lambda \cdot s(L)$ is weighted by $\nu$, then the $\intbigG(\Z)$-orbit of an element $x\in \bigG(\Real) \cdot \Lambda \cdot s(L)$ is represented exactly $\#Z_{\bigG}(x)(\Real)/\#Z_{\intbigG}(x)(\Z)$ times. 
More precisely, for any $x \in \bigG(\Real) \cdot \Lambda \cdot s(L)$ we have 
\begin{equation}\label{equation: multiplicity orbit in fund domain}
	\sum_{x' \in \intbigG(\Z) \cdot x \cap \Siegel \cdot \Lambda\cdot s(L) } \nu(x') = \frac{\#Z_{\bigG}(x)(\Real)}{\#Z_{\intbigG}(x)(\Z)}.
\end{equation}
This follows from an argument similar to \cite[p. 202]{BS-2selmerellcurves} by additionally keeping track of the multiplicity function $\mu$.


In conclusion, for any $\intbigG(\Z)$-invariant subset $A\subset \intbigV(\Z) \cap \bigG(\Real) \cdot \Lambda\cdot s(L)$ we have 
\begin{equation}\label{equation: N(S,a) in terms of lattice points fund domain}
	N(A,X) = \frac{1}{r} \#\left[A \cap (\Siegel\cdot \Lambda\cdot s(L))_{<X}  \right],
\end{equation}
with the caveat that elements on the right-hand side are weighted by $\nu$. 
(Recall that $r = \#Z_{\bigG}(v)(\Real)$ for some $v\in s(L) $.)

\subsection{Averaging and counting lattice points}\label{subsection: averaging and counting lattice points}

We consider an averaged version of (\ref{equation: N(S,a) in terms of lattice points fund domain}) and obtain a useful expression for $N(A,X)$ (Lemma \ref{lemma: bhargavas trick}) using a trick due to Bhargava. 
Then we use this expression to count orbits lying in the `main body' of $\bigV$ using geometry-of-numbers techniques, see Proposition \ref{proposition: counting lattice points main body}.

Fix a compact, semialgebraic subset $G_0 \subset \bigG(\Real) \times \Lambda$ of nonempty interior, that in addition satisfies $K\cdot G_0 = G_0$, $\vol(G_0) = 1$ and the projection of $G_0$ onto $\Lambda$ is contained in $[1,K_0]$ for some $K_0>1$. 
Moreover we suppose that $G_0$ is of the form $G_0'\times [1,K_0]$ where $G_0'$ is a subset of $\bigG(\Real)$. 
Equation (\ref{equation: N(S,a) in terms of lattice points fund domain}) still holds when $L$ is replaced by $hL$ for any $h\in \bigG(\Real)$, by the same argument as above.
Thus for any $\intbigG(\Z)$-invariant $A\subset \intbigV(\Z) \cap \bigG(\Real) \cdot \Lambda\cdot s(L)$ we obtain
\begin{equation}\label{equation: counting bigV using latttice points}
	N(A,X) = \frac{1}{r} \int_{h\in G_0} \# \left[ A\cap (\Siegel\cdot\Lambda \cdot hs(L))_{<X} \right] \, dh,
\end{equation}
again with the caveat that elements on the right are weighted by a function similar to \eqref{equation: weight function nu}.
We use Equation (\ref{equation: counting bigV using latttice points}) to \emph{define} $N(A,X)$ for any subset $A\subset \intbigV(\Z) \cap \bigG(\Real) \cdot \Lambda\cdot s(L)$ which is not necessarily $\intbigG(\Z)$-invariant. 
We can rewrite this integral using the decomposition $\Siegel = \omega \cdot T_c \cdot K$, Lemma \ref{lemma: Haar measure iwasawa decomposition} and an argument similar to \cite[\S2.3]{BS-2selmerellcurves}, which we omit.

\begin{lemma}\label{lemma: bhargavas trick}
	Given $X\geq 1$, $n\in \bigNopp(\Real)$, $t\in T^{\theta}(\Real)$ and $\lambda \in \Lambda$, define $B(n,t,\lambda,X) \coloneqq (nt\lambda G_0 \cdot s(L))_{<X}$.
	Then for any subset $A\subset \intbigV(\Z) \cap (\bigG(\Real)\cdot \Lambda\cdot s(L) ) $ we have
	\begin{equation}\label{equation: lemma bhargavas trick}
		N(A,X) = \frac{1}{r}\int_{\lambda = K_0^{-1}}^{X} \int_{t\in T_c} \int_{n\in \omega} \#\left[A \cap B(n,t,\lambda,X)  \right]\mu(nt)^{-1} \delta_{\bigG}(t)^{-1} \, dn \, dt \, d^{\times}\lambda,
	\end{equation}
	where an element $v\in A\cap B(n,t,\lambda,X)$ on the right-hand side is counted with weight $\#\{h\in G_0 \mid v\in nt\lambda h \cdot s(L)) \}$.

\end{lemma}

Before estimating the integrand of (\ref{equation: lemma bhargavas trick}) by counting lattice points in the bounded regions $B(n,t,\lambda,X)$, we first need to handle the so-called cuspidal region after introducing some notation. 

Let $\Phi_V$ be the set of weights of the $T^{\theta}$-action on $V$. 
Any $v\in \bigV(\Q)$ can be decomposed as $\sum v_{a}$ where $v_{a}$ lies in the weight space corresponding to $a \in \Phi_{\bigV}$.  
For a subset $M\subset \Phi_V$, let $\bigV(M)\subset \bigV$ be the subspace of elements $v$ with $v_{a} = 0$ for all $a \in M$.   
Define $S(M) \coloneqq \bigV(M)(\Q)\cap \intbigV(\Z)$.

Let $a_0 \in X^*(T^{\theta})$ denote the restriction of the highest root $\alpha_0 \in \Phi_{\bigH}$ to $T^{\theta}$.
It turns out that $a_0\in \Phi_V$:
if $H$ is not of type $A_{2n}$, this follows from the fact that the Coxeter number of $H$ is even so the root height of $\alpha_0$ with respect to $S_H$ is odd; if $H$ is of type $A_{2n}$, this can be checked explicitly.

We define $S(\{a_0\})$ (sometimes written $S(a_0)$) as the \define{cuspidal region} and $\intbigV(\Z)\setminus S(a_0)$ as the \define{main body} of $\bigV$. 
The next proposition, proved in \S\ref{subsection: cutting off cusp}, says that the number of irreducible elements in the cuspidal region is negligible. 
\begin{proposition}\label{proposition: cutting off cusp}
	There exists $\delta >0$ such that $N(S(a_0)^{irr},X) = O(X^{\dim V-\delta})$. 
\end{proposition}

\begin{remark}
    The method of proof shows that an explicit $\delta$ can in principle be obtained in each case. For example, we may take $\delta = 1+\varepsilon$ for every $\varepsilon>0$ when $H$ is of type $A_{2n}$ by \cite[Proposition 10.5]{Bhargava-Gross-hyperellcurves}.
\end{remark}

Having dealt with the cuspidal region, we may now count lattice points in the main body using the following proposition \cite[Theorem 1.3]{BarroeroWidmer-lattice}, which strengthens a well-known result of Davenport \cite{Davenport-onaresultofLipschitz}.

\begin{proposition}\label{proposition: count lattice points barroero}
	Let $m,n\geq 1$ be integers, and let $Z\subset \Real^{m+n}$ be a semialgebraic subset. 
	For $T\in \Real^m$, let $Z_T = \{x\in \Real^n\mid (T,x) \in Z\}$, and suppose that all such subsets $Z_T$ are bounded.
	Then for any unipotent upper-triangular matrix $u\in \GL_n(\Real)$, we have
	\begin{align*}
		\#(Z_T \cap u\Z^n) = \vol(Z_T)+O(\max\{1,\vol(Z_{T,j}\}),
	\end{align*}
	where $Z_{T,j}$ runs over all orthogonal projections of $Z_T$ to all $j$-dimensional coordinate hyperplanes $(1\leq j \leq n-1)$. 
	Moreover, the implied constant depends only on $Z$. 
\end{proposition}

\begin{proposition}\label{proposition: counting lattice points main body}
	Let $A = \intbigV(\Z)\cap (\bigG(\Real)\cdot \Lambda\cdot s(L))$.
	Then 
	\begin{equation*}
		N(A \setminus S(a_0) ,X) = \frac{|W_0|}{r}\vol\left(\intbigG(\Z) \backslash \bigG(\Real) \right)\vol((\Lambda\cdot L)_{<X})+o(X^{\dim V}).
	\end{equation*}
\end{proposition}
\begin{proof}
	Choose generators for the weight space $\intbigV_{a}$ (as a finite free $\Z$-module) for every $a\in \Phi_{\bigV}$ and let $\left\Vert \cdot \right\Vert$ denote the supremum norm of $\bigV(\Real)$ with respect to this choice of basis.
	Since the set $\omega\cdot G_0 \cdot s(L)$ is bounded, we can choose a constant $J>0$ such that $\left \Vert v \right\Vert \leq J$ for all $v\in \omega\cdot G_0 \cdot s(L)$.
	Let $F(n,t,\lambda,X) = \left\{ v\in B(n,t,\lambda,X)\mid v_{a_0} \neq 0 \right\}$.
	If $F(n,t,\lambda,X) \cap \intbigV(\Z) \neq \emptyset$, there exists an element $v\in B(n,t,\lambda,X)$ such that $\left\Vert v_{a_0}\right\Vert\geq 1$, hence $\lambda a_0(t)\geq  1/J$.
	
	We wish to estimate $\#[(A\setminus S(a_0)) \cap B(n,t,\lambda,X)]=\#[\intbigV(\Z) \cap F(n,t,\lambda,X)]$ for all $t\in T_c,n\in \omega,\lambda \geq K_0^{-1}$ and $X$ using Proposition \ref{proposition: count lattice points barroero}. 
	An element $v\in F(n,t,\lambda,X)$ has weight $\#\{h\in G_0 \mid v\in nt\lambda h \cdot s(L)) \}$, and $F(n,t,\lambda,X)$ is partitioned into finitely many bounded semialgebraic subsets of constant weight.
	Moreover we have an equality of (weighted) volumes $\vol(F(n,t,\lambda,X))=\vol(B(n,t,\lambda,X))$.
	Since $t$ and $\lambda$ stretch the elements of $\intbigV_a$ by a factor $a(t)$ and $\lambda$ respectively, for any $M\subset \Phi_{\bigV}$ the volume of the projection of $F(n,t,\lambda,X)$ to $\bigV(M)(\Real)$ is bounded above by $O\left(\lambda^{\dim V-\#M}\prod_{a \in \Phi_{\bigV} \setminus M} a(t)\right)$. 
	Since $\Phi_V$ is closed under inversion, we have $\prod_{a\in \Phi_{\bigV}}a(t) = 1$.
	Moreover since $a_0$ is the highest weight of the representation $V$, we know that for every $a\in \Phi_V$ we can write $a = a_0 - \sum_{\beta\in S_G} n_{\beta} \beta$ for some nonnegative rationals $n_{\beta}$.
	Since $t\in T_c$ by assumption, we have $a(t) \geq c^{-\sum n_{\beta}} a_0(t)$. It follows that 
	\begin{align*}
	\lambda^{\dim V-\#M}\prod_{a \in \Phi_{\bigV} \setminus M} a(t)  = \lambda^{\dim V-\#M}\prod_{a \in  M} a(t)^{-1} \ll  \lambda^{\dim V - \#M} a_0(t)^{-\#M}.
	\end{align*}
	
	Putting the results from the previous paragraph together, we conclude by Proposition \ref{proposition: count lattice points barroero} that the number of weighted elements of $[(A\setminus S(a_0)) \cap B(n,t,\lambda,X)]$ is given by:
	\begin{displaymath}
	\begin{cases}
		0 & \text{if } \lambda a_0(t)< 1/J, \\
		\vol(B(n,t,\lambda,X))+O(\lambda^{\dim V-1}a_0(t)^{-1}) & \text{otherwise.}
	\end{cases}
	\end{displaymath}
	By that same proposition, the implied constant in this estimate is independent of $t\in T_c,n\in \omega,\lambda \geq K_0^{-1}$ and $X$.
	Therefore by Lemma \ref{lemma: bhargavas trick} $N(A \setminus S(a_0) ,X) $ equals
	\begin{align}\label{equation: integral estimate for number points in main body}
	\frac{1}{r}\int_{\lambda=K_0^{-1}}^{X} \int_{ t\in T_c,a_0(t)\geq 1/\lambda J } \int_{n\in \omega} \left(\vol(B(n,t,\lambda,X))+O(\lambda^{\dim V-1}a_0(t)^{-1}) \right)\mu(nt)^{-1} \delta_{\bigG}(t)^{-1} \, dn \, dt \, d^{\times}\lambda.
	\end{align}
	We first show that the integral of the second summand is $o(X^{\dim V})$. 
	One easily reduces to showing that 
	\begin{align}\label{equation: first equation tail estimate}
	\int_{\lambda=K_0^{-1}}^{X} \int_{ t\in T_c,a_0(t)\geq 1/\lambda J } \lambda^{\dim V-1}a_0(t)^{-1} \delta_{\bigG}(t)^{-1}\, dt \, d^{\times}\lambda = o(X^{\dim V}).
	\end{align}
	Write $S_G = \{\beta_1,\dots,\beta_k\}$ and identify $T^{\theta}(\Real)^{\circ}$ with $\Real_{>0}^k$ using the isomorphism $t\mapsto (\beta_i(t))$.
	Write $a_0 = \sum h_i \beta_i$ and $\sum_{\Phi_G^+}\beta = \sum \delta_i \beta_i$ with $h_i,\delta_i\in \Q$.
	Since the coefficients $\delta_i$ are strictly positive, there exists an $0<\epsilon<1$ such that $\delta_i-\epsilon h_i >0$ for all $i$.
	Since $\lambda^{1-\epsilon}a_0(t)^{1-\epsilon}\gg 1$ on $\{t\in T_c\mid a_0(t) \geq 1/\lambda J\}$, it follows that
	\begin{align}\label{equation: second equation tail estimate}
	\int_{\lambda=K_0^{-1}}^{X} \int_{ t\in T_c,a_0(t)\geq 1/\lambda J } \lambda^{\dim V-1}a_0(t)^{-1} \delta_{\bigG}(t)^{-1}\, dt \, d^{\times}\lambda 
	\ll \int_{\lambda=K_0^{-1}}^{X} \lambda^{\dim V-\epsilon}\int_{ t\in T_c,a_0(t)\geq 1/\lambda J } a_0(t)^{-\epsilon} \delta_{\bigG}(t)^{-1}\, dt \, d^{\times}\lambda.
	\end{align}
	Since the exponents of $t_i$ in $a_0(t)^{-\epsilon} \delta_{\bigG}(t)^{-1}$ are strictly positive, the inner integral of the right-hand side of \eqref{equation: second equation tail estimate} is bounded independently of $\lambda$.
	It follows that the right-hand side of \eqref{equation: second equation tail estimate} is $\ll \int_{\lambda = K_0^{-1}}^{X}\lambda^{\dim V-\epsilon} d^{\times} \lambda  = O(X^{\dim V-\epsilon})$, as claimed.

	On the other hand, the integral of the first summand in \eqref{equation: integral estimate for number points in main body} is
	$$
	\frac{1}{r}  \int_{g\in \Siegel} \vol\left(\left(g\cdot \Lambda\cdot G_0 \cdot s(L)\right)_{<X}\right)\mu(g)^{-1}\,dg   + o(X^{\dim V}),
	$$
	using the fact that $\vol(B(n,t,\lambda,X)) = O(\lambda^{\dim V})$.
	Lemma \ref{lemma: change of variables section} shows that 
	$$
	\vol\left(\left(g\cdot \Lambda\cdot G_0 \cdot s(L)\right)_{<X}\right)= |W_0|\vol\left(\left(\Lambda\cdot L\right)_{<X}\right)\vol(G_0)=|W_0|\vol\left(\left(\Lambda\cdot L\right)_{<X}\right).
	$$
	The proposition follows from Formula (\ref{equation: weighted volume Siegel set}).
\end{proof}

\subsection{End of the proof of Proposition \ref{prop: counting sections}}

The following proposition is proven in \S\ref{subsection: estimates reducibility stabilizers}.

\begin{proposition}\label{proposition: estimates red}
	Let $\bigV^{red}$ denote the subset of $\Q$-reducible elements $v\in \intbigV(\Z)$ with $v\not\in S(a_0)$. 
	Then $N(\bigV^{red},X) = o(X^{\dim V})$. 
\end{proposition}

We now finish the proof of Proposition \ref{prop: counting sections}.
Again let $A = \intbigV(\Z)\cap (\bigG(\Real)\cdot \Lambda\cdot s(L))$. 
Then 
\begin{equation*}
	N(A^{irr},X) = N(A^{irr} \setminus S(a_0),X)+N(S(a_0)^{irr},X)
\end{equation*}
The second term on the right-hand side is $o(X^{\dim V})$ by Proposition \ref{proposition: cutting off cusp}, and $N(A^{irr} \setminus S(a_0),X)= N(A \setminus S(a_0),X)+o(X^{\dim V})$ by Proposition \ref{proposition: estimates red}.
Using Proposition \ref{proposition: counting lattice points main body}, we obtain
\begin{align*}
	N(A^{irr},X) = \frac{|W_0|}{r}\vol\left(\intbigG(\Z) \backslash \bigG(\Real) \right)\vol((\Lambda\cdot L)_{<X})+o(X^{\dim V}).
\end{align*}
This completes the proof of Proposition \ref{prop: counting sections}, hence also that of Theorem \ref{theorem: counting R-soluble elements, no congruence}.

\subsection{Congruence conditions}\label{subsection: congruence conditions}

We now introduce a weighted version of Theorem \ref{theorem: counting R-soluble elements, no congruence}. If $w\colon \intbigV(\Z) \rightarrow \Real$ is a function and $A\subset \intbigV(\Z)$ is a $\intbigG(\Z)$-invariant subset we define 
\begin{equation}\label{definition N with congruence conditions}
N_w(A,X) \coloneqq \sum_{\substack{v\in \intbigG(\Z)\backslash A \\ \height(v)<X}} \frac{w(v)}{\# Z_{\intbigG}(v)(\Z)}.
\end{equation}
We say a function $w$ is \define{defined by finitely many congruence conditions} if $w$ is obtained from pulling back a function $\bar{w} \colon \intbigV(\Z/M\Z) \rightarrow \Real$ along the projection $\intbigV(\Z) \rightarrow \intbigV(\Z/M\Z)$ for some $M \geq 1$. 
For such a function write $\mu_w$ for the average of $\bar{w}$ where we put the uniform measure on $\intbigV(\Z/M\Z)$. 
The following theorem follows immediately from the proof of Theorem \ref{theorem: counting R-soluble elements, no congruence}, compare \cite[\S2.5]{BS-2selmerellcurves}.

\begin{theorem}\label{theorem: counting finitely many congruence}
	Let $w\colon \intbigV(\Z) \rightarrow \Real $ be defined by finitely many congruence conditions. Then 
	\begin{displaymath}
	N_w(\intbigV(\Z)^{irr} \cap \bigV(\Real)^{sol},X) = \mu_w \frac{|W_0|}{2^g}\vol\left(\intbigG(\Z)\backslash \bigG(\Real)\right) \vol\left(\bigB(\Real)_{<X}  \right)+ o\left(X^{\dim V}\right),
	\end{displaymath}
	where $W_0\in \Q^{\times}$ is the scalar of Lemma \ref{lemma: relations different forms on V,G,B}.
\end{theorem}

Next we will consider infinitely many congruence conditions. 
Suppose we are given for each prime $p$ a $\intbigG(\Z_p)$-invariant function $w_p: \intbigV(\Z_p) \rightarrow [0,1]$ with the following properties:
\begin{itemize}
	\item The function $w_p$ is locally constant outside the closed subset $\{v\in \intbigV(\Z_p) \mid \Delta(v) = 0\} \subset \intbigV(\Z_p)$. 
	\item For $p$ sufficiently large, we have $w_p(v) = 1$ for all $v \in \intbigV(\Z_p)$ such that $p^2 \nmid \Delta(v)$. 
\end{itemize}
In this case we can define a function $w: \intbigV(\Z) \rightarrow [0,1]$ by the formula $w(v) = \prod_{p} w_p(v)$ if $\Delta(v) \neq 0$ and $w(v) = 0$ otherwise. Call a function $w: \intbigV(\Z) \rightarrow [0,1]$ defined by this procedure \define{acceptable}.

\begin{theorem}\label{theorem: counting infinitely many congruence conditions}
	Let $w: \intbigV(\Z) \rightarrow [0,1]$ be an acceptable function. Then
	\begin{displaymath}
	N_w(\intbigV(\Z)^{irr}\cap \bigV^{sol}(\Real) ,X) \leq \frac{|W_0|}{2^g} \left(\prod_p \int_{\intbigV(\Z_p)} w_p(v) \mathrm{d} v \right)  \vol\left(\intbigG(\Z) \backslash \bigG(\Real) \right) \vol\left(\bigB(\Real)_{<X}  \right) + o(X^{\dim V}). 
	\end{displaymath}
\end{theorem}
\begin{proof}
	This inequality follows from Theorem \ref{theorem: counting finitely many congruence}; the proof is identical to the first part of the proof of \cite[Theorem 2.21]{BS-2selmerellcurves}. 
\end{proof}

To obtain a lower bound in Theorem \ref{theorem: counting infinitely many congruence conditions} when infinitely many congruence conditions are imposed, one needs a uniformity estimate that bounds the number of irreducible $\intbigG(\Z)$-orbits whose discriminant is divisible by the square of a large prime. 
The following conjecture is the direct analogue of \cite[Theorem 2.13]{BS-2selmerellcurves}.

\begin{conjecture}\label{conjecture: uniformity estimate}
For a prime $p$, let $\mathcal{W}_p(V)$ denote the subset of $v\in \intbigV(\Z)^{irr}$ such that $p^2 \mid \Delta(v)$. 
Then for any $M>0$, we have 
\begin{align*}
\lim_{X\rightarrow +\infty} \frac{ N(\cup_{p>M} \mathcal{W}_p(V) ,X) }{X^{\dim V}} = O\left(\frac{1}{\log M}\right),
\end{align*}
where the implied constant is independent of $M$.
\end{conjecture}

Conjecture \ref{conjecture: uniformity estimate} is related to computing the density of square-free values of polynomials.
See \cite{Bhargava-geometricsievesquarefree} for some remarks about similar questions, for known results and why these uniformity estimates seem difficult in general.
By an identical proof to that of \cite[Theorem 2.21]{BS-2selmerellcurves}, we obtain:

\begin{proposition}\label{proposition: counting infinitely many congruence assuming a uniformity estimate}
Assume that Conjecture \ref{conjecture: uniformity estimate} holds for $(G,V)$. 
Let $w: \intbigV(\Z) \rightarrow [0,1]$ be an acceptable function. Then
	\begin{displaymath}
	N_w(\intbigV(\Z)^{irr}\cap \bigV^{sol}(\Real) ,X) = \frac{|W_0|}{2^g} \left(\prod_p \int_{\intbigV(\Z_p)} w_p(v) \mathrm{d} v \right)  \vol\left(\intbigG(\Z) \backslash \bigG(\Real) \right) \vol\left(\bigB(\Real)_{<X}  \right) + o(X^{\dim V}). 
	\end{displaymath}
\end{proposition}

\subsection{Estimates on reducibility and stabilisers}\label{subsection: estimates reducibility stabilizers}

In this subsection we give the proof of Proposition \ref{proposition: estimates red} and the following proposition, which will be useful in \S\ref{section: the average size of the 2-Selmer group}. 
\begin{proposition}\label{proposition: estimates bigstab}
	Let $\bigV^{bigstab}$ denote the subset of $\Q$-irreducible elements $v\in \intbigV(\Z)$ with $\#Z_{\bigG}(v)(\Q)>1$. 
	Then $N(\bigV^{bigstab},X) = o(X^{\dim V})$.
\end{proposition}

By the same reasoning as \cite[\S10.7]{Bhargava-Gross-hyperellcurves} it will suffice to prove Lemma \ref{lemma: red and bigstab mod p} below, after having introduced some notation. 

Let $N$ be the integer of \S\ref{subsection: integral structures} and let $p$ be a prime not dividing $N$. We define $\bigV_p^{red}\subset \bigV(\Z_p)$ to be the set of vectors whose reduction mod $p$ is $\F_p$-reducible.
We define $\bigV_p^{bigstab} \subset \bigV(\Z_p)$ to be the set of vectors $v\in \bigV(\Z_p)$ such that $p | \Delta(v)$ or the image of $v$ in $\bigV(\F_p)$ has nontrivial stabiliser in $\bigG(\F_p)$. 

\begin{lemma}\label{lemma: red and bigstab mod p}
	We have 
	$$\lim_{Y\rightarrow +\infty} \prod_{N<p<Y} \int_{\bigV_p^{red}} \,dv = 0,$$
	and similarly
	$$\lim_{Y\rightarrow +\infty} \prod_{N<p<Y} \int_{\bigV_p^{bigstab}} \,dv = 0.$$ 
\end{lemma}
\begin{proof}
	The proof is very similar to the proof of \cite[Proposition 6.9]{Thorne-Romano-E8} using the root lattice calculations of \S\ref{subsection: monodromy}.
	We first treat the case of $V_p^{bigstab}$.
	Let $p$ be a prime not dividing $N$. We have the formula
	$$ \int_{V_p^{bigstab}} dv = \frac{1}{\#\bigV(\F_p)}\# \{v\in \bigV(\F_p) \mid \Delta(v) = 0 \text{ or } Z_{\bigG}(v)(\F_p) \neq 1 \}   .$$
	Since $\{ \Delta = 0 \}$ is a hypersurface we have
	\begin{align}\label{equation: hypersurface mod p points}
	\frac{1}{\#\bigV(\F_p)}\# \{v\in \bigV(\F_p) \mid \Delta(v) = 0\} = O(p^{-1}). 
	\end{align}
	If $v\in \bigV^{\rs}(\F_p)$ then $\#Z_{\bigG}(v)(\F_p)$ depends only on $\pi(v)$ by (the $\Z[1/N]$-analogue of) Lemma \ref{lemma: Z and A are descents of regular centralisers}.
	Therefore if $b\in B^{\rs}(\F_p)$, Proposition \ref{proposition: spread out orbit parametrization galois} and Lang's theorem imply that $\#\bigV_b(\F_p)$ is partitioned into $\#\HH^1(\F_p,J_b[2])$ many orbits under $G(\F_p)$, each of size $\#G(\F_p)/\#J_b[2](\F_p)$. 
	Since $\#J_b[2](\F_p) =\#(J_b(\F_p)/2 J_b(\F_p)) =  \# \HH^1(\F_p,J_b[2])$, we have $\#V^{\rs}(\F_p) = \#G(\F_p) \# B^{\rs}(\F_p)$.

	So to prove the lemma in case of $V_p^{bigstab}$ it suffices to prove that there exists a $0< \delta <1$ such that
	$$
	\frac{1}{\# \bigB^{\rs}(\F_p)}\#\{b\in \bigB^{\rs}(\F_p) \mid J_b[2](\F_p) \neq 1 \} \rightarrow \delta
	$$
	as $p \rightarrow +\infty$.
	We will achieve this using the results of \cite[\S9.3]{Serre-lecturesonNx(p)}. 
	Recall from \S\ref{subsection: a stable Z/2Z-grading} that $\bigT$ is a split maximal torus of $\bigH$ with Lie algebra $\liet$ and Weyl group $W$.
	These objects spread out to objects $\intbigT, \intbigH,\intbigt$ over $\Z$. In \S\ref{subsection: monodromy} we have defined a $W$-torsor $f\colon \liet^{\rs}\rightarrow \bigB^{\rs}$ which extends to a $W$-torsor $\intbigt_S^{\rs} \rightarrow \intbigB_S^{\rs}$, still denoted by $f$. The group scheme $\Jac[2] \rightarrow \intbigB^{\rs}_S$ is trivialised along $f$ and the monodromy action is given by the natural action of $W$ on $N_L$ using the same logic and notation as Proposition \ref{proposition: monodromy of J[2]}.
	Let $C\subset W$ be the subset of elements of $W$ which fix some nonzero element of $N_L$.
	Then \cite[Proposition 9.15]{Serre-lecturesonNx(p)} implies that
	\begin{align}\label{equation: result on Serre lectures on NX(p)}
	\frac{1}{\# \bigB^{\rs}(\F_p)}\#\{b\in \bigB^{\rs}(\F_p) \mid J_b[2](\F_p) \neq 1 \}  = \frac{\#C}{\#W}+O(p^{-1/2}).
	\end{align}
	Since $C\neq W$ by Part 3 of Proposition \ref{proposition: consequences of monodromy}, we conclude the proof of the lemma in this case.
	
	We now treat the case $V_p^{red}$. 
	Again by \eqref{equation: hypersurface mod p points} it suffices to prove that there exists a nonnegative $\delta <1$ such that 
	\begin{align}\label{equation: proof of Vred 1}
	\frac{1}{\#\bigV^{\rs}(\F_p)}\# \{v\in \bigV^{\rs}(\F_p) \mid v \text{ is } \F_p\text{-reducible} \} < \delta 
	\end{align}
	for all sufficiently large $p$. 
	By the first paragraph of the proof of this lemma, there are exactly $\#J_b[2]$ orbits of $V_b(\F_p)$ for all $b\in B^{\rs}(\F_p)$, each of size $\#G(\F_p)/\#J_b[2](\F_p)$. 
	Therefore \eqref{equation: proof of Vred 1} equals
	\begin{align}\label{equation: proof of Vred 2}
	\frac{1}{\#B^{\rs}(\F_p)} \sum_{b\in B^{\rs}(\F_p)} \frac{\#\{\F_p\text{-reducible orbits in }V_b(\F_p) \} }{\#J_b[2](\F_p)}
	\end{align}
	Each summand in \eqref{equation: proof of Vred 2} belongs to the set $\{1/2^{2g},2/2^{2g},\dots,(2^{2g}-1)/2^{2g},1\}$; let $\eta_p$ be the proportion of $b\in B^{\rs}(\F_p)$ for which this summand equals $1$.
	The the quantity in \eqref{equation: proof of Vred 2} is $\leq \eta_p + (1-\eta_p)(2^{2g}-1)/2^{2g} = 1+ (\eta_p-1)/2^{2g}$. 
	By \eqref{equation: result on Serre lectures on NX(p)}, $\eta_p \rightarrow \eta \coloneqq 1-\#C/\#W$ as $p$ tends to infinity. 
	Since $1\in C$, we see that $\eta <1$ and hence $1+ (\eta-1)/2^{2g}<1$, completing the proof of the lemma.
\end{proof}

\begin{proof}[Proof of Proposition \ref{proposition: estimates red}]
We first claim that if $v\in \intbigV(\Z)$ with $b=\bigpi(v)$ is $\Q$-reducible, then for each prime $p$ not dividing $N$ the reduction of $v$ in $\bigV(\F_p)$ is $\F_p$-reducible.
Indeed, either $\Delta(b)=0$ in $\F_p$ (in which case $v$ is $\F_p$-reducible), or $p\nmid \Delta(b)$ and $v$ is $G(\Q)$-conjugate to $\kappa'_b$ for some Kostant section $\kappa'$. 
In the latter case Part 2 of Theorem \ref{theorem: summary orbits of square-free discriminant} implies that $v$ is $G(\Z_p)$-conjugate to $\kappa'_b$, so their reductions are $G(\F_p)$-conjugate, proving the claim.
By a congruence version of Proposition \ref{proposition: counting lattice points main body}, for every subset $L\subset \bigB(\Real)$ considered in Proposition \ref{prop: counting sections} and for every $Y>0$ we obtain the estimate:
\begin{equation*}
    N(\bigV^{red}\cap \bigG(\Real)\cdot \Lambda \cdot s(L),X)\leq C \left(\prod_{N<p<Y} \int_{\bigV_p^{red}} \,dv\right)\cdot X^{\dim V} +o(X^{\dim V}),
\end{equation*}
where $C>0$ is a constant independent of $Y$.
By Lemma \ref{lemma: red and bigstab mod p}, the product of the integrals converges to zero as $Y$ tends to infinity, so $ N(\bigV^{red}\cap \bigG(\Real)\cdot \Lambda \cdot s(L),X)=o(X^{\dim V})$.
Since this holds for every such subset $L$, the proof is complete.
\end{proof}

\begin{proof}[Proof of Proposition \ref{proposition: estimates bigstab}]
Note that we have not used Theorem \ref{theorem: counting R-soluble elements, no congruence} in the proof of Proposition \ref{proposition: estimates red}, but we may use it now to prove Proposition \ref{proposition: estimates bigstab}. 
Again the reduction of an element of $\bigV^{bigstab}$ modulo $p$ lands in $\bigV_p^{bigstab}$ if $p$ does not divide $N$, by Part 1 of Theorem \ref{theorem: summary orbits of square-free discriminant}.
Since $\lim_{X\rightarrow +\infty} N(\bigV^{bigstab},X)/X^{\dim V}$ is $O(\prod_{N<p<Y}\int_{\bigV_p^{red}} dv) $ by Theorem \ref{theorem: counting finitely many congruence} and the product of the integrals converges to zero by Lemma \ref{lemma: red and bigstab mod p}, the proof is complete.
\end{proof}

\subsection{Cutting off the cusp}\label{subsection: cutting off cusp}

In this section we consider the only remaining unproved assertion of this chapter, namely Proposition \ref{proposition: cutting off cusp}.
This is the only substantial part of this paper where we rely on previous papers treating specific cases.
Case $A_{2g}$ is treated in \cite[Proposition 10.5]{Bhargava-Gross-hyperellcurves}; Case $A_{2g+1}$ ($g\geq 1$) is \cite[Proposition 21]{ShankarWang-hypermarkednonweierstrass}; Case $D_{2g+1}$ ($g\geq 2$) is \cite[Proposition 7.6]{Shankar-2selmerhypermarkedpoints}; Case $E_6$ is \cite[Proposition 3.6]{Thorne-E6paper}; Case $E_7,E_8$ is \cite[Proposition 4.5]{Romano-Thorne-ArithmeticofsingularitiestypeE}. 
Note that these authors sometimes use a power of the height that we use.
It remains to consider the case where $H$ is of type $D_{2n}$ and $n\geq 2$. 
We first reduce the statement to a combinatorial result, after introducing some notation. 
This reduction step is valid for any $H$, and we do not yet assume that $H$ is of type $D_{2n}$.

Recall that every element $a\in X^*(T^{\theta})\otimes \Q$ can be uniquely written as $\sum_{i=1}^k n_i(a) \beta_i$ for some rational numbers $n_i(a)$.
We define a partial ordering on $X^*(T^{\theta})\otimes \Q$ by declaring that $a\geq b$ if $n_i(a-b)\geq 0$ for all $i=1,\dots,k$.
By restriction, this induces a partial ordering on $\Phi_V$.
The restriction of the highest root $a_0\in \Phi_V$ is the unique maximal element with respect to this partial ordering.

If $(M_0, M_1)$ is a pair of disjoint subsets of $\Phi_{\bigV}$ we define $S(M_0,M_1) \coloneqq \{v\in \intbigV(\Z) \mid \forall {a}\in M_0, v_{a}=0; \forall {a}\in M_1, v_{a} \neq 0  \}$.
Let $\mathcal{C}$ be the collection of nonempty subsets $M_0\subset \Phi_{\bigV}$ with the property that if $a \in M_0$ and $b \geq a$ then $b \in M_0$. 
Given a subset $M_0 \in\mathcal{C}$ we define $\lambda(M_0) \coloneqq \{ a\in \Phi_{\bigV}\setminus M_0 \mid M_0 \cup\{a\} \in \mathcal{C}\}$, i.e. the set of maximal elements of $\Phi_{\bigV}\setminus M_0$. 

By definition of $\mathcal{C}$ and $\lambda$ we see that $S(\{a_0\}) = \cup_{M_0\in \mathcal{C}} S(M_0,\lambda(M_0))$. 
Therefore to prove Proposition \ref{proposition: cutting off cusp}, it suffices to prove that for each $M_0\in \mathcal{C}$, either $S(M_0,\lambda(M_0))^{irr}=\emptyset$ or $N(S(M_0,\lambda(M_0)),X)=O(X^{\dim V-\epsilon})$ for some $\epsilon >0$. 
By the same logic as \cite[Proposition 3.6 and \S5]{Thorne-E6paper} (itself based on a trick due to Bhargava), the latter estimate holds if there exists a subset $M_1\subset \Phi_{\bigV}\setminus M_0$ and a function $f\colon M_1 \rightarrow \Real_{\geq 0}$ with $\sum_{a\in M_1} f(a) <\#M_0$ such that
\begin{align*}
	\sum_{\beta \in \Phi_{G}^+ }\beta - \sum_{a \in M_0} a  +\sum_{a\in M_1}f(a)a
\end{align*}
has strictly positive coordinates with respect to the basis $S_{G}$. 
It will thus suffice to prove the following combinatorial proposition, which is the analogue of \cite[Proposition 29]{Bhargava-Gross-hyperellcurves}.

\begin{proposition}\label{proposition: combinatorial cutting off the cusp}
Let $M_0 \in \mathcal{C}$ be a subset such that $\bigV(M_0)(\Q)$ contains $\Q$-irreducible elements. 
Then there exists a subset $M_1\subset \Phi_{\bigV} \setminus M_0$ and a function $f\colon M_1 \rightarrow \Real_{\geq 0}$ satisfying the following conditions: 
\begin{itemize}
		\item We have $\sum_{a \in M_1} f(a) < \# M_0$.
		\item For each $i = 1,\dots, k$ we have $\sum_{\beta \in \Phi_{G}^+} n_i(\beta)- \sum_{a \in M_0} n_i(a) + \sum_{a \in M_1} f(a) n_i(a) >0$. 
\end{itemize}
\end{proposition}

We will prove Proposition \ref{proposition: combinatorial cutting off the cusp} in the remaining case $D_{2n}$ in Appendix \ref{section: appendix cutting off cusp D2n}.

\section{The average size of the 2-Selmer group} \label{section: the average size of the 2-Selmer group}

\subsection{An upper bound}

In this chapter we prove Theorem \ref{theorem: intro average size 2-Selmer} stated in the introduction. 
Recall that we write $\sh{E}$ for the set of elements $b\in \intbigB(\Z)$ of nonzero discriminant. 
We recall that we have defined a height function $\height$ for $\sh{E}$ in \S\ref{subsection: heights}. We say a subset $\mathcal{F}\subset \sh{E}$ is defined by \define{finitely many congruence conditions} if $\mathcal{F}$ is the preimage of a subset of $\intbigB(\Z/n\Z)$ under the reduction map $\sh{E} \rightarrow \intbigB(\Z/n\Z)$ for some $n\geq 1$.

\begin{theorem}\label{theorem: main theorem}
	Let $\mathcal{F}\subset \sh{E}$ be a subset defined by finitely many congruence conditions.
	Let $m$ be the number of marked points.
	Then we have 
	\begin{equation*}
	\limsup_{X\rightarrow +\infty} \frac{ \sum_{b\in \mathcal{F},\; \height(b)<X }\# \Sel_2\Jac_b  }{\#  \{b \in \mathcal{F} \mid \height(b) < X \}}	\leq 3\cdot 2^{m-1}.
	\end{equation*}
\end{theorem}

The proof is along the same lines as the discussion in \cite[\S7]{Thorne-Romano-E8}. 

We first prove a `local' result. Recall that $\sh{E}_p$ denotes the set of elements $b\in \intbigB(\Z_p)$ of nonzero discriminant.
Define $\mathcal{F}_p$ as the closure of $\mathcal{F}$ in $\sh{E}_p$, equivalently $\mathcal{F}_p$ is the preimage in $\sh{E}_p$ of a subset of $\intbigB(\Z/n\Z)$ that defines $\mathcal{F}$.

For every $b\in B^{\rs}(\Q)$, consider the subgroup of $J_b(\Q)/2J_b(\Q)$ generated by differences of the marked points $\{\infty_1-\infty_2,\dots,\infty_1-\infty_m\}$. 
The image of this subgroup under the map $J_b(\Q)/2J_b(\Q) \hookrightarrow \Sel_2 J_b$ is by definition the subgroup $\Sel_2^{triv}J_b$ of `marked' elements.
Its complement $\Sel_2^{\text{\sout{triv}}}J_b$ is the subset of `nonmarked' elements.

\begin{proposition}\label{proposition: local result of main theorem}
	Let $b_0 \in \mathcal{F}$. Then we can find for each prime $p$ dividing $N$ an open compact neighbourhood $W_p$ of $b_0$ in $\sh{E}_p$ such that the following condition holds. Let $\mathcal{F}_W = \mathcal{F} \cap \left(\prod_{p | N} W_p \right)$. Then we have 
	\begin{equation*}
	\limsup_{X\rightarrow +\infty} \frac{ \sum_{b\in \mathcal{F}_W,\; \height(b)<X }\# \Sel_2^{\text{\sout{triv}}}J_b  }{\#  \{b \in \mathcal{F}_W \mid \height(b) <X \}}	\leq 2^m.
	\end{equation*}
\end{proposition}
\begin{proof}
	Choose the sets $W_p$ and integers $n_p\geq 0$ for $p| N$ satisfying the conclusion of Corollary \ref{corollary: weak global integral representatives}. 
	If $p$ does not divide $N$, set $W_p = \mathcal{F}_p$ and $n_p = 0$. Let $M \coloneqq \prod_{p} p^{n_p}$. 
	
	For $v\in \intbigV(\Z)$ with $\bigpi(v) = b$, define $w(v) \in \Q_{\geq 0}$ by the following formula:
	\begin{displaymath}
	w(v) = 
	\begin{cases}
	\left( \sum_{v'\in \intbigG(\Z)\backslash \left( \intbigG(\Q)\cdot v \cap \intbigV(\Z) \right)}  \frac{\# Z_{\intbigG}(v')(\Q)}{\# Z_{\intbigG}(v')(\Z)} \right)^{-1} & \text{if }b\in p^{n_p}\cdot W_p \text{ and } \bigG(\Q_p)\cdot v \in \eta_{b}(\Jac_b(\Q_p)/2\Jac_b(\Q_p)) \text{ for all }p, \\
	0 & \text{otherwise.}
	\end{cases}
	\end{displaymath}
	Define $w'(v)$ by the formula $w'(v) = \#Z_{\intbigG}(v)(\Q) w(v)$. Corollaries \ref{corollary: inject 2-Selmer orbits} and \ref{corollary: weak global integral representatives} and Proposition \ref{proposition: orbits corresponding to marked points are reducible} imply that if $b\in M \cdot \mathcal{F}_W$, nonmarked elements in the $2$-Selmer group of $\Jac_b$ correspond bijectively to $\bigG(\Q)$-orbits of $\bigV_b(\Q)$ that intersect $\intbigV(\Z)$ nontrivially, that are $\Q$-irreducible and that are soluble at $\Real$ and $\Q_p$ for all $p$. In other words, in the notation of \eqref{definition N with congruence conditions} we have the formula:
	\begin{equation}\label{equation: selmer count vs orbit count}
	\sum_{\substack{b \in \mathcal{F}_W \\ \height(b) <X}} \#\Sel_2^{\text{\sout{triv}}}J_b
	= \sum_{\substack{b \in M\cdot\mathcal{F}_W \\ \height(b) <M  X}}\#\Sel_2^{\text{\sout{triv}}}J_b
	= N_{w'}(\intbigV(\Z)^{irr}\cap\bigV(\Real)^{sol}  ,M  X).
	\end{equation}
	Proposition \ref{proposition: estimates bigstab} implies that
	\begin{equation}\label{equation: compare w and w'}
	 N_{w'}(\intbigV(\Z)^{irr}\cap \bigV(\Real)^{sol}  ,M  X) =  N_{w}(\intbigV(\Z)^{irr}\cap \bigV(\Real)^{sol},M  X) + o(X^{\dim V}).
	\end{equation}
	It is more convenient to work with $w(v)$ than with $w'(v)$ because $w(v)$ is an acceptable function in the sense of \S\ref{subsection: congruence conditions}. 
	Indeed, for $v\in \intbigV(\Z_p)$ with $\bigpi(v)=b$, define $w_p(v) \in \Q_{\geq 0}$ by the following formula
	\begin{displaymath}
	w_p(v) = 
	\begin{cases}
	\left( \sum_{v'\in \intbigG(\Z_p)\backslash \left( \intbigG(\Q_p)\cdot v \cap \intbigV(\Z_p) \right)}  \frac{\# Z_{\intbigG}(v')(\Q_p)}{\# Z_{\intbigG}(v')(\Z_p)} \right)^{-1} & \text{if }b\in p^{n_p}\cdot W_p \text{ and } \bigG(\Q_p)\cdot v \in \eta_{b}(\Jac_b(\Q_p)/2\Jac_b(\Q_p)   ), \\
	0 & \text{otherwise.}
	\end{cases}
	\end{displaymath}
	Then an argument identical to \cite[Proposition 3.6]{BS-2selmerellcurves} (using that $\intbigG$ has class number $1$ by Proposition \ref{proposition: G has class number one}) shows that $w(v)  =\prod_pw_p(v)$ for all $v\in\intbigV(\Z)$. The remaining properties for $w(v)$ to be acceptable follow from Part 1 of Lemma \ref{lemma: the constants W0 and W} and Theorem \ref{theorem: summary orbits of square-free discriminant}. 
	Using Lemma \ref{lemma: the constants W0 and W} we obtain the formula
	\begin{equation}\label{equation: mass formula w}
	\int_{v\in \intbigV(\Z_p)} w_p(v) d v = |W_0|_p \vol\left(\intbigG(\Z_p) \right) \int_{b \in p^{n_p}\cdot {W_p}} \frac{\#\Jac_b(\Q_p)/2\Jac_b(\Q_p)}{\#\Jac_b[2](\Q_p)}d b.
	\end{equation}
	Using the equality $\#\Jac_b(\Q_p)/2\Jac_b(\Q_p) = |1/2^g|_p  \#\Jac_b[2](\Q_p)$ for all $b\in \sh{E}_p$ (which is a general fact about abelian varieties), we see that the integral on the right-hand side equals $|1/2^g|_p\vol(p^{n_p}\cdot W_p)=|1/2^g|_pp^{-n_p\dim_{\Q}V} \vol(W_p)$.
	Combining the identities (\ref{equation: selmer count vs orbit count}) and (\ref{equation: compare w and w'}) shows that 
	\begin{align*}
	\limsup_{X\rightarrow +\infty} X^{-\dim V} \sum_{\substack{b \in \mathcal{F}_W \\ \height(b) <X}}\#\Sel_2^{\text{\sout{triv}}}J_b
	& = \limsup_{X\rightarrow +\infty} X^{-\dim V}N_w(\intbigV(\Z)^{irr}\cap \bigV(\Real)^{sol} ,M X).
	\end{align*}
	This in turn by Theorem \ref{theorem: counting infinitely many congruence conditions} is less than or equal to
	 \begin{displaymath}
	 \frac{|W_0|}{2^g} \left(\prod_p \int_{\intbigV(\Z_p)} w_p(v) d v \right)  \vol\left(\intbigG(\Z) \backslash \bigG(\Real) \right) 2^{\dim B} M^{\dim V}.
	 \end{displaymath}
	Using (\ref{equation: mass formula w}) this simplifies to
	\begin{displaymath}
	\vol\left(\intbigG(\Z)\backslash \intbigG(\Real) \right) \prod_p \vol\left(\intbigG(\Z_p)\right) 2^{\dim B}\prod_{p} \vol(W_p).
	\end{displaymath}
	On the other hand, an elementary point count shows that 
	\begin{displaymath}
	\lim_{X\rightarrow +\infty} \frac{\#  \{b \in \mathcal{F}_W \mid ht(b) <X \}}{X^{\dim V}} = 2^{\dim B }\prod_p \vol(W_p).
	\end{displaymath}
	We conclude that 
	\begin{displaymath}
	\limsup_{X\rightarrow +\infty} \frac{ \sum_{b\in \mathcal{F}_W,\; ht(b)<X }\#\Sel_2^{\text{\sout{triv}}}J_b   }{\#  \{b \in \mathcal{F}_W \mid ht(b) <X \}}	\leq \vol\left(\intbigG(\Z)\backslash \intbigG(\Real) \right) \cdot \prod_p \vol\left(\intbigG(\Z_p)\right).
	\end{displaymath}
	Since the Tamagawa number of $\intbigG$ is $2^m$ (Proposition \ref{proposition: tamagawa number}), the proposition follows. 
\end{proof}

To deduce Theorem \ref{theorem: main theorem} from Proposition \ref{proposition: local result of main theorem}, choose for each $i\geq1$ sets $W_{p,i} \subset\sh{E}_p$ (for $p$ dividing $N$) such that if $W_i = \sh{E} \cap \left( \prod_{p | N} W_{p,i} \right)$, then $W_i$ satisfies the conclusion of Proposition \ref{proposition: local result of main theorem} and we have a countable partition $\mathcal{F} =\mathcal{F}_{W_1}\sqcup\mathcal{F}_{W_1} \sqcup \cdots$. 
By an argument identical to the proof of Theorem 7.1 in \cite{Thorne-Romano-E8}, we see that for any $\varepsilon >0$, there exists $k\geq 1$ such that 
\begin{displaymath}
\limsup_{X\rightarrow +\infty}  \frac{ \sum_{\substack{b \in \sqcup_{i\geq k} \mathcal{F}_{W_i} , \height(b) <X  }} \#\Sel_2^{\text{\sout{triv}}}J_b    }{ \# \{b \in\mathcal{F} \mid \height(b) <X  \}  }<\varepsilon.
\end{displaymath}
This implies that 
\begin{align*}
\limsup_{X\rightarrow +\infty}  \frac{ \sum_{\substack{b\in\mathcal{F} , \height(b) <X  }} \#\Sel_2^{\text{\sout{triv}}}J_b     }{ \# \{b \in\mathcal{F}\mid \height(b) <X  \}  } &\leq 2^m \limsup_{X\rightarrow +\infty}\frac{\# \{b \in \sqcup_{i<k} \mathcal{F}_{W_i} \mid \height(b) <X \}  }{ \# \{b \in \mathcal{F} \mid \height(b) <X  \}  } +\varepsilon \\
&\leq 2^m+\varepsilon.  
\end{align*}
Since the above inequality is true for any $\varepsilon >0$, it is true for $\varepsilon=0$.
Since the subgroup $\Sel^{triv}J_b$ has size at most $2^{m-1}$, we conclude the proof of Theorem \ref{theorem: main theorem}.

\begin{remark}
A small modification of the above argument shows that Theorem \ref{theorem: main theorem} remains valid when $\mathcal{F}\subset \sh{E}$ is the subset of so-called `minimal' elements, namely those elements $b\in \mathcal{E}$ with $N^{-1} \cdot b \not\in \intbigB(\Z)$ for all integers $N\geq 1$. 
\end{remark}

\subsection{A conditional lower bound}\label{subsection: a conditional lower bound}

We show that the upper bound in Theorem \ref{theorem: main theorem} is sharp if we assume Conjecture \ref{conjecture: uniformity estimate}.
We first need to establish (unconditionally) a lower bound for the subgroup of marked elements $\Sel^{triv} J_b$.

\begin{proposition}\label{proposition: average size trivial selmer group}
Let $\mathcal{F}\subset \sh{E}$ be a subset defined by finitely many congruence conditions.
Then the limit
\begin{align*}
\lim_{X\rightarrow +\infty} \frac{ \#\{b\in \mathcal{F} \mid \height(b)<X, \#\Sel_2(J_b)^{triv} = 2^{m-1} \}  }{\#\{b\in \mathcal{F} \mid \height(b)<X \} }
\end{align*}
exists and equals $1$.
\end{proposition}
\begin{proof}
Let $b\in \mathcal{F}$ and consider the maximal torus $Z_H(\kappa_b)$ of $H$.
By Proposition \ref{proposition: orbits corresponding to marked points are reducible}, $ \#\Sel_2(J_b)^{triv} = 2^{m-1}$ if and only if the map $Z_G(\kappa_b) \rightarrow Z_{H^{\theta}}(\kappa_b)$ is surjective on $\Q$-points.
The Galois action on $Z_H(\kappa_b)$ induces a homomorphism $\Gal(\Q^s\mid \Q) \rightarrow W$ by Proposition \ref{proposition: monodromy of J[2]}, where $W$ is the Weyl group of the split torus $T\subset H$ with character group $L$.
If this homomorphism is surjective, then $Z_{H^{\theta}}(\kappa_b)(\Q) = T[2]^W = (L^{\vee}/2L^{\vee})^W = \{0\}$  by Part 1 of Proposition \ref{proposition: consequences of monodromy}, so $Z_G(\kappa_b) \rightarrow Z_{H^{\theta}}(\kappa_b)$ is automatically surjective on $\Q$-points.

It therefore suffices to prove that the limit
\begin{align*}
\lim_{X\rightarrow +\infty} \frac{ \#\{b\in \mathcal{F} \mid \height(b)<X, \Gal(\Q^s\mid \Q) \rightarrow W \text{ surjective}  \}  }{\#\{b\in \mathcal{F} \mid \height(b)<X \} }
\end{align*}
exists and equals $1$.
This follows from a version of Hilbert's irreducibility theorem; see \cite[Theorem 2.1]{Cohen-distributionofgaloisgroupshilbertirreducibility}, adapted as in \cite[\S5, Notes (iii)]{Cohen-distributionofgaloisgroupshilbertirreducibility} to account for the fact that the coordinates of $B$ have unequal weights.
\end{proof}

\begin{theorem}
Assume that Conjecture \ref{conjecture: uniformity estimate} holds for $(G,V)$.
Let $\mathcal{F}\subset \sh{E}$ be a subset defined by finitely many congruence conditions.
Then the limit
\begin{equation*}
	\lim_{X\rightarrow +\infty} \frac{ \sum_{b\in \mathcal{F},\; \height(b)<X }\# \Sel_2^{ \text{\sout{triv}}  } \Jac_b  }{\#  \{b \in \mathcal{F} \mid \height(b) < X \}}
\end{equation*}
exists and equals $2^{m}$.
Moreover, the average size of the $2$-Selmer group $\Sel_2 J_b$ exists and equals $3\cdot 2^{m-1}$.
\end{theorem}
\begin{proof}
The proof of the first statement is identical to the proof of Theorem \ref{theorem: main theorem}, using Proposition \ref{proposition: counting infinitely many congruence assuming a uniformity estimate} instead of Theorem \ref{theorem: counting infinitely many congruence conditions}.
The second statement follows from the first and Proposition \ref{proposition: average size trivial selmer group}.
\end{proof}

\appendix

\section{Cutting off the cusp for \texorpdfstring{$D_{2n}$}{D2n}}
\label{section: appendix cutting off cusp D2n}

In this appendix chapter we prove Proposition \ref{proposition: cutting off cusp} in the case that $H$ is of type $D_{2n}$ for all $n \geq 2$.
The methods employed here are fairly standard but somewhat intricate, and are sometimes inspired by \cite[\S7.2.1]{Shankar-2selmerhypermarkedpoints}.
In \S\ref{subsection: recollections on even orthogonal groups} we recall some results and notation on groups of type $D_n$. 
In \S\ref{subsection: explicit model of stable involution of D2n} we make the representation $(G,V)$ in the case $D_{2n}$ and some related objects explicit.
In \S\ref{subsection: reducibility conditions} we establish sufficient conditions for a vector $v\in V(\Q)$ to be $\Q$-reducible.
In \S\ref{subsection: bounding the remaining cusp integrals} we finish the proof of Proposition \ref{proposition: cutting off cusp}.

\subsection{Recollections on even orthogonal groups}\label{subsection: recollections on even orthogonal groups}

Let $n\geq 2$ be an integer.
Let $W$ be a $2n$-dimensional $\Q$-vector space with basis $\mathcal{B} =  \{e_1,\dots,e_n,e_n^*,\dots,e_1^*\}$. Let $b$ be the symmetric bilinear form with the property that $b(e_i,e_j) = b(e_i^*,e_j^*) = 0$ and $b(e_i,e_j^*) = \delta_{ij}$ for all $1\leq i,j \leq n$. 
For every linear map $f\colon W\rightarrow W$, there is a unique adjoint linear map $f^*\colon W\rightarrow W$ satisfying $b(fv,w) = b(v,f^*w)$ for all $v,w\in W$.
We define the $\Q$-algebraic group $H \coloneqq \SO(W,b) = \{g\in \SL(W) \mid gg^* = 1\}$.
Then $\lieh \coloneqq \Lie H$ can be naturally identified with $\{f\in \End(W) \mid f+f^* = 0\}$.
Below we will make various aspects of the semisimple group $H$ explicit.

Using $\mathcal{B}$ to represent an element $f \colon W \rightarrow W$ as a $(2n) \times (2n)$-matrix $A$, $f^*$ corresponds to reflecting $A$ along its antidiagonal. 
Consider the maximal torus $T = \{\text{diag}(t_1,\dots,t_n,t_n^{-1}, \dots, t_1^{-1})\} \subset H$.
Its character group $X^*(T)$ is freely generated by the characters $(t_1,\dots) \mapsto t_i$ with $1\leq i\leq n$, and we abusively denote these characters by $t_i \in X^*(T)$ too.

\paragraph{Root system}
The roots of $\lieh$ with respect to $T$ are given by 
\begin{align}
    \Phi_H = \{ \pm t_i \pm t_j \mid 1\leq i\neq j\leq n\}\subset X^*(T).
\end{align}
The standard upper triangular Borel subgroup of $\GL_{2n}$ (with respect to the basis $\mathcal{B}$ of $W$) determines a root basis of $\Phi_H$, given by
\begin{align}\label{equation: D2n basis SH}
    S_H =\{t_1-t_2,\dots,t_{n-1}-t_n,t_{n-1}+t_n\}.
\end{align}
We denote the elements of this root basis by $\alpha_1,\dots,\alpha_n$.
The highest root of $\Phi_H$ with respect to $S_H$ is $t_1+t_2$, which has height $2n-3$.

\paragraph{Weyl group}
The Weyl group $W_H$ is isomorphic to $S_n\rtimes (\Z/2\Z)^{n-1}$. 
Explicitly, elements of $W_H$ correspond to pairs $(\sigma,(\epsilon_i))$, where $\sigma\in S_n$ is a permutation and $\epsilon_i \in \{\pm 1\}$ is a sign for each $1\leq i\leq n$, with the property that $\prod_i \epsilon_i = 1$.
An element $(\sigma,(\epsilon_i))$ acts on $X^*(T)$ via the rule $t_i \mapsto \epsilon_i t_{\sigma(i)}$.
In particular, $-1\in W_H$ if and only if $n$ is even.

\paragraph{Stable involution} 
To describe the stable involution of $H$ in \S\ref{subsection: explicit model of stable involution of D2n} in the case that $n$ is even, we determine the elements $s\in T$ with the property that $\alpha(s)=-1$ for every simple root $\alpha \in S_H$.
Such an $s$ is not uniquely determined since $H$ is not adjoint, but it is uniquely determined up to multiplication by the element $(-1,\dots,-1)\in T$ of the centre of $H$.
Using the description \eqref{equation: D2n basis SH} we see that $s$ is of the form
\begin{align}\label{equation: stable involution Dn case}
    \pm (1,-1,1,-1,\dots,(-1)^n).
\end{align}


\paragraph{Change of variables}
We record the following computation which will be useful in \S\ref{subsection: bounding the remaining cusp integrals}:
\begin{align}\label{equation: change of variables Dn bases}
    \begin{cases}
    t_i  = \alpha_i+\dots+\alpha_{n-2} +\frac{1}{2}(\alpha_{n-1}+\alpha_n) ,&\, (1\leq i\leq n-2)\\
    t_{n-1} = \frac{1}{2}(\alpha_{n-1}+\alpha_{n}), &\\
    t_n  = \frac{1}{2}(-\alpha_{n-1}+\alpha_{n}) . &
    \end{cases}
\end{align}

\paragraph{Sum of positive roots}
The sum of the positive roots of $\Phi_H$ with respect to $S_H$ is
\begin{align}
    \sum_{\alpha \in \Phi_H^{+}} \alpha &= 2(n-1)t_1+ 2(n-2) t_2 + \dots + 2 t_{n-1}\\
    &= \sum_{k=1}^{n-2} k(2n-k-1)\alpha_k + \frac{n(n-1)}{2}(\alpha_{n-1}+\alpha_n).
\end{align}

\paragraph{Discriminant}
Let $\liet \coloneqq \Lie T$ and write an element of $\liet$ as $\text{diag}(t_1,\dots,t_n,-t_n,\dots,-t_1)$.
Let $\Delta\in \Q[\lieh]^H$ be the discriminant polynomial of $H$, defined as the image of $\prod_{\alpha\in \Phi_H} \alpha$ under the Chevalley isomorphism $\Q[\mathfrak{t}]^{W_H}\rightarrow \Q[\lieh]^{H}$ of Proposition \ref{proposition: classical invariant theory of H on lieh}.
Let $\Pff \in \Q[\lieh]^H$ be the image of the product $\prod_{i=1}^n t_i \in k[\liet]^{W_H}$ under the Chevalley isomorphism.
If we write $t_i \coloneqq t_{2n+1-i}$ for $n+1\leq i\leq 2n$, then we may compute that 
\begin{align}
    \prod_{1\leq i <j\leq 2n} (t_i-t_j)^2 = \prod_{\alpha\in \Phi_H} \alpha(t)^2 \prod_{i=1}^nt_i^2 .
\end{align}
It follows that if we write $\chi_v$ for the characteristic polynomial of a square matrix $v$ and $\disc(\chi_v)$ for its discriminant, we have the identity 
\begin{align}\label{equation: discriminant of Dn in terms of A2n}
    \disc(\chi_v) = \Delta(v)^2\Pff(v)^2
\end{align}
for every $v\in \lieh$.

\subsection{An explicit model for the split stable involution of \texorpdfstring{$D_{2n}$}{D{2n}}}
\label{subsection: explicit model of stable involution of D2n}

Let $n\geq 2$ be an integer. 
Let $W_1$ be the $\Q$-vector space with basis $\{e_1,\dots,e_n,e_n^*,\dots,e_1^*\}$, and let $b_1$ be the symmetric bilinear form with the property that $b_1(e_i,e_j) = b_1(e_i^*,e_j^*) = 0$ and $b_1(e_i,e_j^*) = \delta_{ij}$ for all $1\leq i,j \leq n$. 
Let $W_2$ be the $\Q$-vector space with basis $\{f_1,\dots,f_n,f_n^*,\dots,f_1^*\}$, and let $b_2$ be the bilinear form of $W_2$ constructed similarly to $b_1$ with $e_i$ and $e_i^*$ replaced by $f_i$ and $f_i^*$.
Let $(W,b) \coloneqq (W_1,b_1) \oplus (W_2,b_2)$.
Let $H' \coloneqq \SO(W,b)$, let $H$ be the quotient of $H'$ by its centre of order $2$ and let $\lieh \coloneqq \Lie H = \Lie H'$.
With respect to the basis 
\begin{align}\label{equation: first basis of SO4n}
\{e_1,\dots,e_n,e_n^*,\dots,e_1^*,f_1,\dots,f_n,f_n^*,\dots,f_1^*\},
\end{align}
the adjoint of a $(4n)\times (4n)$-block matrix 
\begin{align}
    \begin{pmatrix}
    A & B \\
    C & D
    \end{pmatrix}
\end{align}
with respect to $b$ is given by 
\begin{align}
    \begin{pmatrix}
    A^* & C^* \\
    B^* & D^*
    \end{pmatrix}.
\end{align}
Here if $X$ is a $(2n)\times (2n)$-matrix we write $X^*$ for its reflection around the antidiagonal.
It follows that in this basis $\lieh$ is given by 
\begin{align}
    \left\{ \begin{pmatrix} B & A \\ -A^* & C \end{pmatrix} \mid B^* = -B, C^* = -C \right\}.
\end{align}

\paragraph{Stable involution}
The ordered basis
\begin{align}\label{equation: second basis SO4n}
    \{e_1,f_1,\dots,e_n,f_n,f_n^*,e_n^*,\dots,f_1^*,e_1^*\}
\end{align}
of $W$ determines a maximal torus and root basis of $H$ as in \S\ref{subsection: recollections on even orthogonal groups}.
Let $\theta$ be the involution of $H$ constructed using the recipe in \S\ref{subsection: a stable Z/2Z-grading} with respect to this root basis.
Since $-1$ is contained in the Weyl group of $H$ (as observed in \S\ref{subsection: recollections on even orthogonal groups}), $\theta$ is inner.
The description of \eqref{equation: stable involution Dn case} shows that with respect to the first basis \eqref{equation: first basis of SO4n} of $W$, $\theta$ is given by conjugating by the element $s' = \text{diag}(1,\dots,1,-1,\dots,-1)$, where the first $2n$ entries are $1$'s and the last $2n$ entries are $-1$'s.
Using this description, it is easy to see that
\begin{align*}
    \lieg &\coloneqq \lieh^{\theta} = \left\{ \begin{pmatrix} B & 0 \\ 0 & C \end{pmatrix} \mid B^* = -B, C^* = -C \right\}, \\
    V &\coloneqq \lieh^{\theta=-1} = \left\{ \begin{pmatrix} 0 & A \\ -A^* & 0 \end{pmatrix} \mid A \in \text{Mat}_{2n,2n} \right\}.
\end{align*}
Moreover $G \coloneqq (H^{\theta})^{\circ}$ is isomorphic to $(\SO(W_1)\times \SO(W_2))/\Delta(\mu_2)$, where $\Delta(\mu_2)$ denotes the image of the diagonal inclusion of $\mu_2$ into the centre $\mu_2\times \mu_2$ of $\SO(W_1)\times \SO(W_2)$.
Using these identifications, we see that the map $$\begin{pmatrix} 0 & A \\ -A^* & 0 \end{pmatrix}  \mapsto A$$ establishes a bijection between $V$ and the representation $\Hom(W_2,W_1)$, where $(g,h)\in \SO(W_1)\times \SO(W_2)$ acts on $f\colon W_2\rightarrow W_1$ via $g\circ f\circ h^{-1}$.
In terms of matrices, the action is given by $(g,h)\cdot A = gAh^{-1}$.
We will typically view an element of $V(\Q)$ as a $(2n) \times (2n)$-matrix $A$ or a linear operator $f\colon W_2\rightarrow W_1$.

\paragraph{Roots}
Let $T'$ be the maximal torus $\text{diag}(t_1,\dots,t_n,t_n^{-1},\dots,t_1^{-1},s_1,\dots,s_n,s_n^{-1},\dots,s_1^{-1})$ of $H'$ (again using the basis \eqref{equation: first basis of SO4n}), and let $T$ be its image in $H$.
Then $T$ is a maximal torus of $H$ and $G$; let $\Phi_H$ and $\Phi_G$ be the corresponding sets of roots. 
Let $W_H = N_H(T)/T$ and $W_G = N_G(T)/T$ be the respective Weyl groups.
The basis \eqref{equation: second basis SO4n} determines a set of positive roots $\Phi_H^+$ of $H$ (as in \S\ref{subsection: recollections on even orthogonal groups}) and by restriction a set of positive roots $\Phi_G^+$ of $G$.
The corresponding simple roots are given by: 
\begin{align*}
    S_H &= \{t_1-s_1,s_1-t_2,\dots,s_{n-1}-t_n,t_n-s_n,t_n+s_n\}, \\
    S_G &= \{t_1-t_2,\dots,t_{n-1}-t_n,t_{n-1}+t_n\} \cup \{s_1-s_2,\dots,s_{n-1}-s_n,s_{n-1}+s_n\}. 
\end{align*}
We label the elements of $S_G$ by $ \{\beta_1,\dots,\beta_{n-1},\beta_n\} \cup \{\gamma_1,\dots,\gamma_{n-1},\gamma_n\}$.
We have $\Phi_H = \Phi_G \sqcup \Phi_V$ and $\Phi_V = \{\pm t_i \pm s_j \mid 1\leq i,j\leq n\}$.

\paragraph{Component group}
Let $s$ be the image of $s'$ (defined below \eqref{equation: second basis SO4n}) in $T(\Q)$.
Lemma \ref{lemma: component group using weyl groups} shows that the inclusion $N_{H^{\theta}}(T)\hookrightarrow H^{\theta}$ induces an isomorphism $Z_{W_H}(s)/W_G\simeq H^{\theta}/G$.
In fact, let 
\begin{align}
    \Omega \coloneqq \{ w\in W_H \mid w(S_G) = S_G\}.
\end{align}
Then using the description of the Weyl group of $H$ and $G$ from \S\ref{subsection: recollections on even orthogonal groups} we see that $Z_{W_H}(s) = W_G \rtimes \Omega$ and $\Omega \simeq \Z/2\times \Z/2$.
Explicit generators of $\Omega$ are given by $\omega_1,\omega_2$, where
\begin{align*}
\omega_1 \colon& t_i \leftrightarrow s_i, \\
\omega_2\colon  & 
\begin{cases} s_i \mapsto s_i , t_i \mapsto t_i, & (1\leq i \leq n-1) \\
 t_n \mapsto -t_n , s_n \mapsto -s_n. &
\end{cases}
\end{align*}

\paragraph{Weights}
Using the description of elements of $V$ as $(2n)\times (2n)$ matrices, we organise the weights $\Phi_V$ using the position of their eigenspaces:
\begin{align}\label{equation: big weight matrix}
\left(
\begin{array}{c c ;{2pt/2pt} c | c;{2pt/2pt} c c}
t_1-s_1 & \cdots  & t_1-s_n & t_1+s_n & \cdots & t_1+s_1 \\
\vdots & \ddots &  \vdots & \vdots & \ddots & \vdots \\
\hdashline[2pt/2pt]
t_n-s_1 & \cdots & t_n-s_n & t_n+s_n & \cdots & t_n+s_1 \\
\hline
-t_n-s_1 & \cdots & -t_n-s_n & -t_n+s_n & \cdots & -t_n+s_1 \\
\hdashline[2pt/2pt]
\vdots & \ddots& \vdots & \vdots & \ddots & \vdots \\
-t_1-s_1 & \cdots  & -t_1-s_n & -t_1+s_n & \cdots & -t_1+s_1 \\
\end{array}   
\right)
\end{align}
The group $\Omega$ acts on the set of weights $\Phi_V$ as follows: $\omega_1$ flips the elements of $\Phi_V$ along the antidiagonal of \eqref{equation: big weight matrix}, and $\omega_2$ swaps the two middle rows and the two middle columns.

\paragraph{The partial ordering}
Recall from \S\ref{subsection: cutting off cusp} that we have defined a partial ordering on $X^*(T')$ by declaring that $a\geq b$ if and only if $a-b$ has nonnegative coordinates with respect to the basis $S_G$.
Note that this partial ordering is preserved by the action of $\Omega$ on $X^*(T')$.

We describe the induced partial ordering on the subset $\Phi_V$ using the organisation \eqref{equation: big weight matrix}.
We first consider the restriction of the partial ordering to the rows and columns of \eqref{equation: big weight matrix}.
Let $1\leq i\leq 2n$ and write $t_i \coloneqq -t_{2n+1-i}, s_i \coloneqq -s_{2n+1-i}$ if $i\geq n+1$.
The Hasse diagram of the partial ordering restricted to the weights of row $i$ is given by:
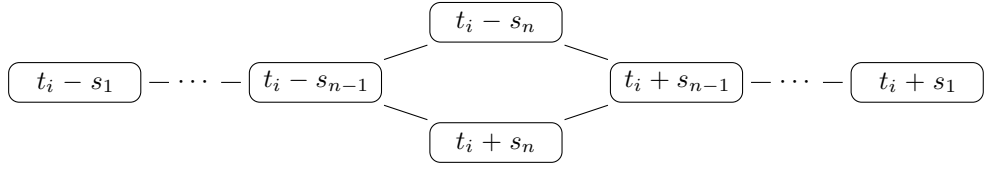
\begin{center}
\begin{tikzpicture}
[scale = 0.8,
block/.style ={rectangle, draw, text width=5em,align=center, rounded corners, minimum height=1.5em}
]
	\node[block] (1) at (-1,0) {$t_i-s_{1}$} ;
	\node (2) at (1,0) {$\cdots$};
	\node[block] (3) at (3,0) {$t_i-s_{n-1}$};
	\node[block] (4) at (6,1) {$t_i-s_n$};
	\node[block] (5) at (6,-1) {$t_i+s_n$};
	\node[block] (6) at (9,0) {$t_i+s_{n-1}$};
	\node (7) at (11,0) {$\cdots$};
	\node[block] (8) at (13,0) {$t_i+s_1$};
	\draw[shorten <= 3pt, shorten >=3pt] (1) -- (2);
	\draw[shorten <= 3pt, shorten >=3pt] (2) -- (3);
	\draw[shorten <= 3pt, shorten >=3pt] (3) -- (4);
	\draw[shorten <= 3pt, shorten >=3pt] (3) -- (5);
	\draw[shorten <= 3pt, shorten >=3pt] (4) -- (6);
	\draw[shorten <= 3pt, shorten >=3pt] (5) -- (6);
	\draw[shorten <= 3pt, shorten >=3pt] (6) -- (7);
	\draw[shorten <= 3pt, shorten >=3pt] (7) -- (8);
\end{tikzpicture}
\end{center}
(In this diagram, $a\leq b$ if and only if $b$ is to the right of $a$.)
The Hasse diagram of column $2n+1-i$ is given by swapping the roles of $s_j$ and $t_j$ in the above diagram for every $1\leq j\leq 2n$.
The partial ordering $\Phi_V$ is the one generated by the relations between two elements lying in the same row or column.

For example, $t_1+s_1$ is the maximal element of $\Phi_V$, and the restriction of the partial ordering to the four $n\times n$ blocks is given by: $a\leq b$ if and only if $b$ is to the top right of $a$.

\paragraph{Regular nilpotent element}
For $\alpha \in \Phi_V$, let $X_{\alpha}$ be the $(2n)\times (2n)$-matrix with coefficient $1$ at the entry corresponding to $\alpha$ using \eqref{equation: big weight matrix} and zeroes elsewhere.
The element $E \coloneqq \sum_{\alpha \in S_H} X_{\alpha}$ is a regular nilpotent element of $V(\Q)$ and gives rise to an $\liesl_2$-triple $(E,X,F)$ and Kostant section $\kappa  = E+ \mathfrak{z}_{\lieh}(F)$, see \S\ref{subsection: Kostant sections}.
Since elements of $\mathfrak{z}_{\lieh}(F)$ are supported on $\Phi_V \cap \Phi_H^{-}$ (where $\Phi^-_H = \Phi_H\setminus \Phi_H^+$), every element of $\kappa$ is of the form
\begin{align}\label{equation: kostant section matrix D2n}
\left(
\begin{array}{c c c c | c c c c}
1 & 0 & \cdots & 0 & 0 & \cdots & \cdots & 0\\
* & \ddots & \ddots & \vdots & \vdots & \ddots & \ddots & \vdots \\
\vdots & \ddots & \ddots & \vdots & 0 & \ddots & \ddots & \vdots\\
* & \cdots & * & 1 & 1 & 0 & \cdots & 0 \\ \hline
* & \cdots & \cdots & * & * & 1 & \cdots & 0\\
\vdots & \ddots & \ddots & \vdots & \vdots & \ddots & \ddots & \vdots \\
\vdots & \ddots & \ddots & \vdots & \vdots & \ddots & \ddots & 1 \\
* & \cdots & \cdots & * & * & \cdots & \cdots & * \\ 
\end{array}   
\right)
\end{align}
If $\omega\in \Omega$, then $E_{\omega} \coloneqq \sum_{\alpha \in \omega(S_H)} X_{\alpha}$ is again regular nilpotent and gives rise to a Kostant section $\kappa_{\omega}$.
Then $\{\kappa_{\omega} \mid \omega\in \Omega \}$ is a full set of representatives of $G(\Q)$-orbits of Kostant sections.

\subsection{Reducibility conditions}\label{subsection: reducibility conditions}

Recall that if $A$ is a $(2n)\times (2n)$-matrix then $A^*$ denotes its reflection along the antidiagonal.

\begin{proposition}\label{proposition: discriminant in SO4n case}
Let $k/\Q$ be a field and $A\in V(k)$.
The following are equivalent:
\begin{enumerate}
    \item $A$ is a regular semisimple element of $V(k)$;
    \item $AA^*$ is a regular semisimple $(2n) \times (2n)$-matrix (in other words, the characteristic polynomial of $AA^*$ has distinct roots in $k^s$);
    \item $A^*A$ is a regular semisimple $(2n)\times (2n)$-matrix.
\end{enumerate}
\end{proposition}
\begin{proof}
Let $\Delta\in \Q[V]^G$ be the discriminant polynomial of $\lieh\simeq \so_{4n}$ restricted to $V$.
Let $D$ be the block matrix $\begin{pmatrix} 0 & A \\ -A^* & 0 \end{pmatrix}$.
If $C$ is a square matrix, write $\chi_C\in k[X]$ for its characteristic polynomial. 
If $f$ is a polynomial, write $\disc(f)$ for its discriminant in the usual sense. 
The identity \eqref{equation: discriminant of Dn in terms of A2n} implies that $\disc(\chi_D) = \Delta(A)^2\cdot \Pff(D)^2$.
We have $\Pff(D) = \pm \det(A)$ since both square to $\det(D)$, so
\begin{align}\label{equation: disc zero Dn 1}
\disc(\chi_D) =\Delta(A)^2 \cdot \det(A)^2.
\end{align}
On the other hand, if $f(X) = g(X^2)$ for some polynomial $g\in k[X]$, then it is elementary to check that $\disc(f) = \pm \disc(g)^2 f(0)$.
Moreover, by calculating determinants of block matrices we have $\chi_D(X) = \chi_{-AA^*}(X^2)$.
Therefore
\begin{align}\label{equation: disc zero Dn 2}
    \disc(\chi_D) = \pm \disc(\chi_{-AA^*})^2\cdot \det(A)^2.
\end{align}
Both identities \eqref{equation: disc zero Dn 1} and \eqref{equation: disc zero Dn 2} hold in $\Q[V][X]$, i.e. they hold when the coefficients of $A$ are interpreted as variables.
Since $\det \in \Q[V]$ is not identically zero, it follows that $\Delta(A) = \pm \disc(\chi_{-AA^*})=\pm \disc(\chi_{AA^*})$. Since $\chi_{AA^*} = \chi_{A^*A}$ we also have $\Delta(A) = \pm \disc(\chi_{A^*A})$.
Since $A$ is a regular semisimple element of $V(k)$ if and only if $\Delta(A)\neq 0$, the proposition follows.
\end{proof}

In the next lemma, we organise the set $\Phi_V$ using the matrix \eqref{equation: big weight matrix}, and we recall from \S\ref{subsection: averaging and counting lattice points} that for a subset $M\subset \Phi_V$ we have defined $V(M)$ as the subspace of $v = \sum_{a\in \Phi_V} v_a \in V$ with the property that $v_a=0$ for all $a\in M$.

\begin{corollary}\label{corollary: explicit discriminant zero conditions}
Suppose that a subset $M\subset \Phi_V$ satisfies at least one of the following conditions:
\begin{enumerate}
    \item $M$ contains a top right $i\times (2n+1-i)$ block for some $1\leq i \leq 2n $;
    \item $M$ contains the top right $i\times j$ and $j\times i$ blocks for some $i,j \geq 1$ satisfying $i+j = 2n$.
\end{enumerate}
Then every element of $V(M)(\Q)$ is not regular semisimple.
\end{corollary}
\begin{proof}
\begin{enumerate}
    \item We may suppose (using the fact that $A\mapsto A^*$ preserves regular semisimplicity) that $i\leq n$.
    Let $X_1 = \text{span}\{e_1,\dots,e_n\} \subset W_1$.
    Then $\GL(X_1)$ embeds inside $\SO(W_1)$, using the map $g\mapsto \begin{pmatrix} g & 0 \\ 0 & (g^*)^{-1} \end{pmatrix}$.
    Suppose that $A\in V(M)(\Q)$.
    Using the $\GL(X_1)$-action to put the top left $i \times (i-1)$ block of $A$ in row echelon form, we may suppose that $M$ contains the top right $1\times 2n$ block; in other words, we may suppose that $i=1$.
    In that case, the matrix $AA^*$ has zeroes on the first row and the last column. 
    This implies that the characteristic polynomial of $AA^*$ is divisible by $X^2$, which implies that $A$ is not regular semisimple by Proposition \ref{proposition: discriminant in SO4n case}.
    \item Assume that $i\leq j$ and let $A\in V(M)(\Q)$.
    Then the matrix $B = AA^*$ is of the form:
    \begin{align}
    \left(
    \begin{array}{c | c| c}
    B_1 & 0 & 0 \\
    \hline
    * & B_2 & 0 \\
    \hline 
    * & * & B_3 
    \end{array}   
    \right).
    \end{align}
    Here $B_1,B_3$ are $i\times i$ matrices and $B_2$ is a $(j-i)\times (j-i)$ matrix (it is possible that $i=j$).
    Recall that $\chi_C$ denotes the characteristic polynomial of a square matrix $C$. 
    Then we have $\chi_{AA^*}  = \chi_{B_1} \chi_{B_2} \chi_{B_3}$.
    Since $B_3 = B_1^*$, the polynomial $\chi_{AA^*} = \chi_{B_1}^2\chi_{B_2}$ has repeated roots.
    By Proposition \ref{proposition: discriminant in SO4n case}, this shows that $A$ is not regular semisimple.
\end{enumerate}
\end{proof}

Recall from \S\ref{subsection: explicit model of stable involution of D2n} that we may interpret an element $A\in V(\Q)$ as a linear map $W_2\rightarrow W_1$.
Using the perfect pairings $b_i$ on $W_i$ we may thus interpret $A^*$ as a linear map $W_1\rightarrow W_2$, and $AA^*$ as a linear map $W_1\rightarrow W_1$.

\begin{proposition}\label{proposition: Q-reducibility SO4n using maximal isotropic stuff}
Let $k/\Q$ be a field and $A\in V(k)$.
Assume that there exists an $(n-1)$-dimensional subspace $X\subset W_1$ such that $\text{span}\{X,AA^*(X)\}$ is an $n$-dimensional isotropic subspace of $(W_1,b_1)$.
Then $A$ is $k$-reducible.
\end{proposition}
\begin{proof}
If $A$ is not regular semisimple then $A$ is $k$-reducible by definition, so assume that $A$ has invariants $b\in B^{\rs}(k)$. (Recall that $B = V\GIT G$.)
The Grassmannian of $n$-dimensional isotropic subspaces of $(W_1,b_1)$ has two connected components called rulings of $(W_1,b_1)$, and the subspaces spanned by $\{e_1^*,\dots,e_{n}^*\}$ and $\{e_1^*,\dots,e_{n-1}^*,e_n\}$ lie in distinct rulings \cite[\S2.2]{Wang-maximallinearspacescontainedbaselociquadrics}; let $\mathcal{R}$ be the ruling containing the subspace $X_1 \coloneqq \text{span}\{e_1^*,\dots,e_n^*\}$. 
Using the $\Omega$-action we may assume that $\text{span}\{X,AA^*(X)\}$ lies in $\mathcal{R}$.
To prove the proposition, it suffices to prove the claim that $G(k)$ acts simply transitively on the set of pairs $(D,Y)$, where $D\in V_b(k)$ and $Y\subset W_1$ is an $(n-1)$-dimensional subspace such that $\text{span}\{Y,DD^*(Y)\}$ is an $n$-dimensional isotropic subspace of $(W_1,b_1)$ contained in the ruling $\mathcal{R}$.
Indeed, the description of the Kostant section $\kappa_b$ from \eqref{equation: kostant section matrix D2n} shows that $(\kappa_b,X_1)$ is such a pair, so the claim implies that $(A,X)$ and $(\kappa_b,X_1)$ are $G(k)$-conjugate.
The proof of the claim is identical to the proof of \cite[Proposition 4.4]{Shankar-2selmerhypermarkedpoints} using the results of \cite[\S2.2.2]{Wang-maximallinearspacescontainedbaselociquadrics}; we omit the details.
\end{proof}

\begin{corollary}\label{corollary: explicit k-reducibility condition D2n}
Suppose that $M$ contains the top right $(n-1)\times (n+1)$ and $n\times (n-1)$ blocks. Then every element of $V(M)(\Q)$ is $\Q$-reducible.
\end{corollary}
\begin{proof}
If $A\in V(M)(\Q)$, a computation shows that the entries of $AA^*$ in the top right $(n-1)\times n$ and $n\times (n-1)$ blocks are zero. In other words, $AA^*$ looks like:
\begin{align}
\left(
\begin{array}{c c ;{2pt/2pt} c | c;{2pt/2pt} c c}
* & *   & * & 0 & 0 & 0 \\
*  & *  &  * & 0 & 0 & 0 \\
\hdashline[2pt/2pt]
*  & *  & * & * & 0 & 0 \\
\hline
*  & *  & * & *  & *  & *  \\
\hdashline[2pt/2pt]
 * &  * &  *  &  *  &  * & * \\
* & *  & *  & *  & * &* \\
\end{array}   
\right).
\end{align}
It follows that the subspace $X = \text{span}\{e_{n-1}^*,\dots,e_1^*\}$ satisfies the assumptions of Proposition \ref{proposition: Q-reducibility SO4n using maximal isotropic stuff}.
\end{proof}

Recall that $\mathcal{C}$ denotes the collection of subsets $M$ of $\Phi_V$ with the property that for all $a,b\in \Phi_V$ with $b\in M$ and $a\geq b$ it follows that $a\in M$.
Also recall the description of the partial ordering on $\Phi_V$ in \S\ref{subsection: explicit model of stable involution of D2n}.

\begin{proposition}\label{proposition: summary Q-reducibility conditions}
Let $M\in \mathcal{C}$ and suppose that $V(M)(\Q)$ contains $\Q$-irreducible elements.
Then the following properties hold:
\begin{enumerate}
    \item $\{t_i-s_i,s_i-t_i\} \subset \Phi_V\setminus M$ for all $1\leq i\leq n-1$;
    \item $\{t_n-s_n,t_n+s_n,-t_n+s_n,-t_n-s_n\} \subset \Phi_V\setminus M$;
    \item for every $1\leq i\leq n-2$, either $t_i-s_{i+1}$ or $s_i-t_{i+1}$ lies in $\Phi_V\setminus M$;
    \item $\#(\{t_{n-1}-s_n,t_{n-1}+s_n,t_n+s_{n-1},-t_n+s_{n-1}\} \cap M) \leq 2$.
\end{enumerate}
\end{proposition}
\begin{proof}
Note that if $\omega\in \Omega$ then $V(M)(\Q)$ contains $\Q$-irreducible elements if and only if $V(\omega(M))(\Q)$ does. 
The first three parts follow from applying Corollary \ref{corollary: explicit discriminant zero conditions} to $\omega(M)$ for all $\omega\in \Omega$ and properties of the partial ordering of $\Phi_V$.
Part 4 follows from applying Corollary \ref{corollary: explicit k-reducibility condition D2n} to $\omega(M)$ for $\omega\in M$.
\end{proof}

The reader is invited to visualise the conditions of Proposition \ref{proposition: summary Q-reducibility conditions} using the organisation of the weights $\Phi_V$ of \eqref{equation: big weight matrix}.

\subsection{Bounding the remaining cusp integrals}\label{subsection: bounding the remaining cusp integrals}

Let $\mathcal{C}^{good}$ be the subset of $\mathcal{C}$ consisting of those $M\in \mathcal{C}$ that satisfy Condiditions 1--4 of Proposition \ref{proposition: summary Q-reducibility conditions}.
Note that $\omega(M) \in \mathcal{C}^{good}$ if $M\in \mathcal{C}^{good}$ and $\omega\in \Omega$.

\begin{lemma}\label{lemma: every element of SG sum of two elements not in M D2n}
If $M\in\mathcal{C}^{good}$, then every element of $S_G$ is of the form $a_1+a_2$ for some $a_1,a_2 \in \Phi_V \setminus M$.
\end{lemma}
\begin{proof}
Using the $\Omega$-action it suffices to consider $\beta_1,\dots,\beta_{n-1}$.
For $1\leq i\leq n-2$, we have identities
\begin{align*}
    \beta_i = t_i-t_{i+1} &= (t_i-s_i)+\boxed{(s_i-t_{i+1})} \\
    &= \boxed{(t_i-s_{i+1})}+(s_{i+1}-t_{i+1}).
\end{align*}
At least one of the two boxed terms is in $\Phi_V\setminus M$ by Part 3 of Proposition \ref{proposition: summary Q-reducibility conditions}, and the unboxed terms are always in $\Phi_V\setminus M$ by Part 1 of that proposition.
To treat $\beta_{n-1}$, consider the identities
\begin{align*}
    \beta_{n-1} = t_{n-1}-t_n &= \boxed{(t_{n-1}-s_n)}+(s_n-t_n) \\
    &=\boxed{(t_{n-1}+s_n)} + (-s_n-t_n) \\
    &=(t_{n-1}-s_{n-1}) + \boxed{(s_{n-1}-t_n)}.
\end{align*}
One of the three boxed terms must be contained in $\Phi_V\setminus M$ by Part 4 of Proposition \ref{proposition: summary Q-reducibility conditions}, and all the unboxed terms are contained in $\Phi_V\setminus M$ by Parts 1 and 2 of that proposition.
\end{proof}

The discussion in \S\ref{subsection: cutting off cusp} and Proposition \ref{proposition: summary Q-reducibility conditions} show that in order to prove Proposition \ref{proposition: cutting off cusp} when $H$ is of type $D_{2n}$ for all $n\geq 2$, it suffices to prove the following proposition.

\begin{proposition}\label{proposition: cutting off cusp D2n case}
For every $M\in \mathcal{C}^{good}$, there exists a function $f\colon \Phi_V \setminus M \rightarrow \Real_{\geq 0}$ with the following properties:
\begin{enumerate}
    \item $\sum_{a\in \Phi_V\setminus M} f(a) < \#M$;
    \item the vector 
    \begin{align}\label{equation: weight of general cusp datum}
    \sum_{\beta\in \Phi_G^+}\beta- \sum_{a\in M} a+\sum_{a\in \Phi_V\setminus M} f(a)a
    \end{align}
    has strictly positive coefficients with respect to the basis $S_G$.
\end{enumerate}
\end{proposition}
We prove Proposition \ref{proposition: cutting off cusp D2n case} using induction on $n$.
The base case $n=2$ is easy to check explicitly, and also follows from Case 1 of the proof of Proposition \ref{proposition: cutting off cusp D2n induction step} below (which only assumes $n\geq 2$). See \cite[p. 1217]{Thorne-averagesizeelliptictwomarkedfunctionfields} which considers the $D_{4}$ case in detail and also proves this base case.

To perform the induction step, let $\Phi_V^{[1]}$ be the subset of $\Phi_V$ of vectors of the form $\{\pm t_1 \pm s_i\} \cup \{ \pm t_i \pm s_1\}$; in other words, $\Phi_V^{[1]}$ consists of the first and last rows and columns of \eqref{equation: big weight matrix}.
Similarly let $\Phi_G^{[1]}$ for the subset of roots of $\Phi_G$ that have a nonzero coordinate at $\beta_1$ or $\gamma_1$ in the root basis $S_G$. 
Write $\Phi_V = \Phi_V^{[1]} \sqcup \Phi_V^{[n-1]}$ and $\Phi_G = \Phi_G^{[1]} \sqcup \Phi_G^{[n-1]}$.
Then $\Phi_V^{[n-1]}$ and $\Phi_G^{[n-1]}$ arise from the constructions of \S\ref{subsection: explicit model of stable involution of D2n} with $n$ replaced by $n-1$.
Moreover $\sum_{\beta\in \Phi_G^{[1],+}}\beta = (2n-2)t_1+(2n-2)s_1$.
To prove Proposition \ref{proposition: cutting off cusp D2n case}, it therefore suffices to prove the following statement.

\begin{proposition}\label{proposition: cutting off cusp D2n induction step}
Let $n\geq 3$ be an integer, let $M\in \mathcal{C}^{good}$ and write $M^{[1]} \coloneqq M \cap \Phi_V^{[1]}$.
Then there exists a function $f^{[1]}\colon \Phi_V \setminus M \rightarrow \Real_{\geq 0}$ with the following properties:
\begin{enumerate}
    \item $\sum_{a\in \Phi_V\setminus M} f^{[1]}(a) < \#M^{[1]}$;
    \item the vector 
    \begin{align}\label{equation: weight of induction step}
    (2n-2)t_1+(2n-2)s_1 - \sum_{a\in M^{[1]}} a+\sum_{a\in \Phi_V\setminus M} f^{[1]}(a)a
    \end{align}
    has strictly positive coefficients with respect to the basis $S_G$.
\end{enumerate}
\end{proposition}
\begin{proof}
Note that if the proposition is true for $M$, it is also true for $\omega(M)$ for every $\omega \in \Omega$. 
We may therefore replace $M$ by a $\Omega$-conjugate in what follows.
We also note that it suffices to find for each $M\in \mathcal{C}^{good}$ a function $f^{[1]}\colon \Phi_V\setminus M \rightarrow \Real_{\geq 0}$ that satisfies the first property and such that \eqref{equation: weight of induction step} has nonnegative (instead of positive) coefficients with respect to $S_G$.
Indeed, by Lemma \ref{lemma: every element of SG sum of two elements not in M D2n}, every element of $\beta\in S_G$ is a sum $a_1+a_2$ of two elements of $\Phi_V\setminus M$ so by adding to $f^{[1]}$ the function $a_1\mapsto \epsilon, a_2\mapsto \epsilon$ for some very small $\epsilon$, we may ensure that $f^{[1]}$ has strictly positive coefficient at every element of $S_G$. 

We will distinguish three cases, after introducing some notation.
We say $M\in \mathcal{C}^{good}$ is \define{bounded} if there exists a function $f^{[1]}\colon \Phi_V\setminus M\rightarrow \Real_{\geq 0}$ satisfying the conclusions of Proposition \ref{proposition: cutting off cusp D2n induction step}.
If $M\in \mathcal{C}^{good}$ we write $w_1(M) \coloneqq (2n-2)t_1+(2n-2)s_1- \sum_{a\in M^{[1]}} a$.
Recall that if $n+1 \leq i\leq 2n$ then we write $t_i \coloneqq -t_{2n+1-i}$ and $s_i \coloneqq -s_{2n+1-i}$.
We use $O(\geq 0)$ as a shorthand for an element of $X^*(T')$ that has nonnegative coordinates with respect to $S_G$.
We also recall the useful formulae \eqref{equation: change of variables Dn bases}.

\textbf{Case 1.}
Suppose that $M^{[1]} \subset \{t_1-s_n,t_1+s_n,\dots,t_1+s_1,\dots,t_n+s_1,-t_n+s_1\}$.
Let $a$ (respectively $b$) be the number of elements of $M^{[1]}$ contained in the first row (respectively last column). 
Then $1\leq a,b \leq n+1$ and since we may switch the roles of $t_n$ and $-t_n$ and similarly for $\pm s_n$, we may assume that $M^{[1]} = \{t_1+s_a,\dots,t_1+s_1,\dots,t_b+s_1\}$. 
Using the $\Omega$-action we may assume that $a \geq b$.
We have 
\begin{align*}
    w_1(M) &= (2n-2-a)t_1-t_2-\dots-t_b \\
    &+(2n-2-b)s_1-s_2-\dots -s_a.
\end{align*}
If $a,b\leq n-1$, then \eqref{equation: change of variables Dn bases} shows that $w_1(M)$ has positive $S_G$ coefficients, so $M$ is evidently bounded.
We may therefore assume that $a=n$ or $n+1$.
A computation shows that the coefficients of $w_1(M)$ at $\{\beta_1,\dots,\beta_{n-2}\} \cup \{\gamma_1,\dots, \gamma_{n-2}\}$ are nonnegative (in fact at least $n+1-a$), so we focus on the coefficients at $\beta_{n-1},\beta_n,\gamma_{n-1},\gamma_n$.
If $b\leq n-1$, then $w_1(M)  = O(\geq 0)+ \frac{1}{2}(n-a)(\beta_{n-1}+\beta_n)$.
This is negative only when $a=n+1$, so assume that this is the case.
Using Lemma \ref{lemma: every element of SG sum of two elements not in M D2n}, write $\beta_{n-1} = a_1+a_2, \beta_n = a_3+a_4$ for some $a_i\in \Phi_V\setminus M$. 
Choose a function $f^{[1]}\colon \Phi_V\setminus M \rightarrow \Real_{\geq 0}$ such that $f^{[1]}$ supported on $\{a_1,\dots,a_4\}$, such that $\sum_{a\in \Phi_V\setminus M} f^{[1]}(a) a = \frac{1}{2}(\beta_{n-1}+\beta_n)$ and such that $\sum_{a\in \Phi_V\setminus M} f^{[1]}(a)<\#M^{[1]}$. Such a function exists since $2< \# M^{[1]}$ and shows that $M$ is bounded in this case.

It remains to treat the case where $n \leq a,b \leq n+1$.
A calculation shows that 
\begin{align*}
    w_1(M) = \begin{cases} O(\geq 0)-\frac{1}{2}\beta_n-\frac{1}{2}\gamma_n & \text{if }(a,b) = (n,n), \\
    O(\geq 0)-\beta_n & \text{if }(a,b)=(n+1,n), \\
    O(\geq 0)-\frac{1}{2}(\beta_{n-1}+\beta_n+\gamma_{n-1}+\gamma_n) & \text{if }(a,b)=(n+1,n+1).
    \end{cases}
\end{align*}
In each of these cases, we can use Lemma \ref{lemma: every element of SG sum of two elements not in M D2n} to show that $M$ is bounded.

\textbf{Case 2.}
Suppose that $M^{[1]}$ contains an element of the form $t_1-s_a$ with $2\leq a \leq n-1$ and is contained in the set $\{t_1-s_a,\dots,t_1+s_1,\dots,t_n+s_1,-t_n+s_1\}$.
Let $b$ be the number of elements of $M^{[1]}$ contained in the last column, so that $\#M^{[1]} = (n-a)+n+b$.
Then $1\leq b\leq n+1$, and by using the $\Omega$-action we may assume that $M^{[1]} = \{t_1-s_a,\dots,t_1+s_1,\dots,t_b+s_1\}$.
We have 
\begin{align*}
    w_1(M) &= (a-3)t_1-t_2-\dots -t_b \\
    &+ (2n-2-b)s_1-s_2-\dots-s_{a-1}.
\end{align*}
By assumption $t_1-s_{a-1} \in \Phi_V\setminus M$ and $2n-a-b\geq 0$. We compute that 
\begin{align}\label{equation: induction step D2n case 2}
    w_1(M) + (2n-a-b)(t_1-s_{a-1}) &= (2n-b-3)t_1-t_2-\dots -t_b + O(\geq 0) 
\end{align}
All the $S_G$ coefficients of the above expression are nonnegative unless $2n-b-2-b<0$, in other words unless $b\geq n$.
Therefore if $b<n$ the function $f^{[1]}$ mapping $t_1-s_{a-1}$ to $2n-a-b$ and all other elements of $\Phi_V\setminus M$ to zero shows that $M$ is bounded, since $2n-a-b < \#M^{[1]}=2n-a+b$.
If $b=n$, \eqref{equation: induction step D2n case 2} equals $O(\geq 0) - \beta_n$.
Therefore Lemma \ref{lemma: every element of SG sum of two elements not in M D2n} and the inequality $(2n-a-b)+2 < \#M^{[1]}$ show that $M$ is bounded in this case.
If $b=n+1$, \eqref{equation: induction step D2n case 2} equals $O(\geq 0) - \beta_{n-2} - \beta_{n-1} -\beta_n$.
Therefore Lemma \ref{lemma: every element of SG sum of two elements not in M D2n} and the inequality $(2n-a-b)+6 < \#M^{[1]}$ (which holds since $n\geq 3$) show that $M$ is again bounded in this case.

\textbf{Case 3.}
Suppose that $M^{[1]}$ is of the form $\{t_1-s_a,\dots,t_1-s_n,t_1+s_n,\dots,t_1+s_1,\dots, t_n+s_1, -t_n+s_1,\dots,-t_b+s_1\}$ for some $2\leq a,b \leq n-1$.
Using the $\Omega$-action we may assume that $a\geq b$.
A calculation using \eqref{equation: change of variables Dn bases} shows that
\begin{align*}
    w_1(M) &= (a-3)t_1-t_2-\dots-t_{b-1}  \\ &+ (b-3)s_1-s_2-\dots -s_{a-1} \\
    &= (a-3)\beta_1 + \dots + (a-b-1)\beta_{b-1}+\dots + (a-b-1)\beta_{n-2} + \frac{1}{2}\left(a-b-1\right)\left(\beta_{n-1}+\beta_n\right) \\
    &+ (b-3)\gamma_1+\dots + (b-a-1)\gamma_{a-1} +\dots + (b-a-1)\gamma_{n-2} + \frac{1}{2}(b-a-1)(\gamma_{n-1}+\gamma_n).
\end{align*}
By assumption $s_1-t_{b-1}\in \Phi_V\setminus M$ and $a\geq b$, and we compute that  
\begin{align*}
    w_1(M) + (a-b)(s_1-t_{b-1}) &= O(\geq 0) -\beta_{b-1} -\dots -\beta_{n-2} -\frac{1}{2}(\beta_{n-1}+\beta_n)\\
    &+O(\geq 0) -\gamma_{a-1}-\dots-\gamma_{n-2}-\frac{1}{2}(\gamma_{n-1}+\gamma_n) .
\end{align*}
Therefore to prove that $M$ is bounded it suffices to prove (using Lemma \ref{lemma: every element of SG sum of two elements not in M D2n}) that 
\begin{align*}
    (a-b)+2((n-b+1)+(n-a+1)) < \#M^{[1]} = 4n-a-b+1.
\end{align*}
This inequality is equivalent to $2b>3$, which is true since $b\geq 2$.

\textbf{Conclusion.}
Since $M\in \mathcal{C}^{good}$, every $M$ has an $\Omega$-conjugate that falls under one of the above three cases. 
\end{proof}

\bibliographystyle{alpha}

\end{document}